\documentclass[final,3p,times,twoside]{elsarticle}
\usepackage{hyperref}
\usepackage{color}

\usepackage{graphicx}
\usepackage{amsthm,amsmath,amssymb}
\usepackage{amsmath,amssymb}
\usepackage{mathrsfs}
\usepackage{amssymb}
\usepackage{subfigure}
\usepackage{bm}
\usepackage{amsfonts}
\usepackage{amsxtra}
\usepackage{upgreek}
\usepackage{mathtools}
\usepackage{enumerate}
\usepackage{framed}
\usepackage{indentfirst}
\usepackage{pstricks}
\usepackage{tikz}
\usepackage{lipsum}
\usepackage[scr=boondoxo]{mathalfa}
\usepackage{booktabs}
\usepackage{cases}
\usepackage{tcolorbox}
\usepackage{tabularx}
\usepackage{multirow}

\usepackage{algorithm}
\usepackage{algpseudocode}

\usepackage{makecell}

\newtheorem{theorem}{Theorem}
\newtheorem{corollary}{Corollary}
\newtheorem{lemma}{Lemma}
\newtheorem{remark}{Remark}

\usepackage[misc]{ifsym}

\usetikzlibrary{patterns}
\graphicspath{ {figures/} }

\journal{Elsevier}

\begin{document}

\begin{frontmatter}

\title{\textbf{M}ulti\textbf{L}evel \textbf{V}ariational \textbf{M}ulti\textbf{S}cale (ML-VMS) framework for large-scale simulation}

%% Group authors per affiliation: PLEASE ADD ALL AFFLICATIONS FOR LIU ALSO HIDENN-AI, LLC
\author[1]{Lei Zhang}
\author[2]{Jiachen Guo}
\author[3]{Shaoqiang Tang}
\author[4]{Thomas J.R. Hughes}
\author[5,6]{Wing Kam Liu}

\address[1]{School of Engineering Science, University of Chinese Academy of Sciences, Beijing 100049, China}
\address[2]{Theoretical and Applied Mechanics Program, Northwestern University, 2145 Sheridan Road, Evanston, 60201, IL, USA}

\address[3]{School of Mechanics and Engineering Science, Peking University, Beijing 100871, China}

\address[4]{Institute for Computational Engineering and Sciences, The University of Texas at Austin, 201 East 24th Street, Stop C0200, Austin, 78712, TX, USA}

\address[5]{Department of Mechanical Engineering, Northwestern University, 2145 Sheridan Road, Evanston, 60201, IL, USA}

\address[6]{HIDENN-AI, LLC, 1801 Maple Ave, Evanston, 60201, IL, USA}
% \begin{Large}

\begin{abstract}
In this paper, we propose the MultiLevel Variational MultiScale (ML-VMS) method, a novel approach that seamlessly integrates a multilevel mesh strategy into the Variational Multiscale (VMS) framework. A key feature of the ML-VMS method is the use of the Convolutional Hierarchical Deep Neural Network (C-HiDeNN) as the approximation basis, which enables fine-grained control over the trade-off between computational efficiency and interpolation accuracy. The framework employs a coarse mesh throughout the domain, with localized fine meshes placed only in subdomains of high interest, such as those surrounding a source. Solutions at different resolutions are robustly coupled through the variational weak form and interface conditions. Crucially, our method departs from existing VMS-based multilevel approaches by approximating the fine-scale solution directly using the fine-scale basis functions. 
Compared to existing multilevel methods, ML-VMS (1) can couple an arbitrary number of mesh levels across different scales using variational multiscale framework; (2) allows approximating functions with arbitrary orders with linear finite element mesh due to the C-HiDeNN basis; (3) is supported by a rigorous theoretical error analysis; (4) features several tunable hyperparameters (e.g., order $p$, patch size $s$) with a systematic guide for their selection. We first show the theoretical error estimates of ML-VMS. Then through numerical examples, we demonstrate that ML-VMS with the C-HiDeNN takes less computational time than the FEM basis given comparable accuracy.
Furthermore, we incorporate a space-time reduced-order model (ROM) based on C-HiDeNN-Tensor Decomposition (TD) into the ML-VMS framework. For a large-scale single-track laser powder bed fusion (LPBF) transient heat transfer problem that is equivalent to a full-order finite element model with $10^{10}$ spatial degrees of freedom (DoFs), our 3-level ML-VMS C-HiDeNN-TD achieves an approximately 5,000x speedup on a single CPU over a single-level linear FEM-TD ROM.

\end{abstract}


% \end{Large}

\end{frontmatter}

% \linenumbers

% \begin{Large}

\textbf{Highlights}
\begin{itemize}

\item A novel multiscale framework that resolves disparate spatial and temporal scales using linear meshes.

\item An efficient approach with arbitrary interpolation orders by adjusting hyperparameters without increasing DoFs.

\item A space-time formulation provides unconditional stability, solving all time steps at once.

%     \item ML-VMS framework resolves disparate scales with compatible linear meshes.

% \item Achieves higher-order accuracy via hyperparameters without increasing DoFs.

% \item Unconditionally stable space-time model solves for all time steps simultaneously.

    % \item \textbf{A Novel Framework for Multiscale Problems:} We introduce the Multilevel Variational MultiScale (ML-VMS) framework, which utilizes different levels of linear meshes to efficiently resolve disparate length scales, thereby avoiding the mesh incompatibility issues common in other multilevel higher-order methods.

    % \item \textbf{Efficient Higher-Order Approximations without Increasing DoFs:} Our approach allows for arbitrary interpolation orders across all levels simply by adjusting C-HiDeNN hyperparameters, not by increasing degrees of freedom (DoFs). This unique formulation significantly improves computational efficiency and reduces memory usage.

    % \item \textbf{Unconditionally Stable Space-Time Formulation:} Using a C-HiDeNN space-time TD model, the framework avoids explicit time-stepping stability constraints, permitting different time steps across different meshes and simultaneous solution for all time steps all at once.

    % \item \textbf{Extensible to Parameterized Multiphysics Systems:} The ML-VMS framework can be readily extended to complex parameterized multiphysics simulations, creating data-free digital twin models that operate without expensive offline data generation.

\end{itemize}
\section{Introduction}

\begin{figure}[htbp]
\centering
\includegraphics[width=4in]{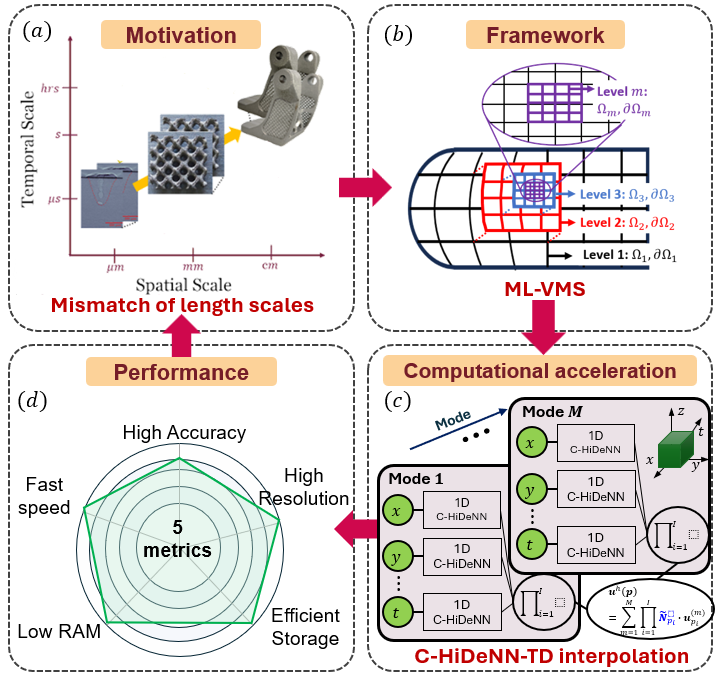}
\caption{MultiLevel Variational MultiScale (ML-VMS) method: (a) capturing disparate length scales in LPBF poses a significant challenge to current available methods; (b) ML-VMS provides an integrated framework for modeling a general $m$ level system; (c) ML-VMS with C-HiDeNN-TD achieves significant computational acceleration as a reduced-order model; (d) The proposed ML-VMS C-HiDeNN-TD achieves superior performance in all 5 performance metrics.}
\label{fig:summary}
\end{figure}

Many challenging engineering problems, such as fatigue stress concentration analysis \cite{kafka2021image, mashayekhi2020fatigue}, additive manufacturing (AM) process simulation and design \cite{debroy2018additive,li2024statistical}, simulating, verifying, and manufacturing complex electronic systems in electronic design automation (EDA) \cite{garyfallou2018large,li2006vector,jiao2007layered}, are usually characterized by multiscale characteristics that pose significant challenges to standard numerical tools to resolve cross-scale physics. For example, in the simulation of the laser powder bed fusion (LPBF) process, it is critical to resolve the fine-scale melt pool physics, which occurs over microseconds and micrometers, while simultaneously capturing the global-scale heat transfer and resulting residual stresses, a process that unfolds over days and meters\cite{debroy2018additive}. As shown in Fig. \ref{fig:summary}(a), capturing such cross-scale physics plays a vital role in determining the microstructure and the quality of the resulting part. Despite significant advancements in computer hardware and the development of GPU parallel computing over the past few decades, simulating the entire computational domain using fine resolutions is not only a waste of computational resources but also remains intractable with conventional simulation tools. Consequently, there is a pressing need for novel modeling approaches that can effectively and efficiently bridge vast scale differences by simultaneously capturing microscale melt-pool dynamics and part-scale transient heat transfer.

A common strategy to address multiscale features is the use of adaptive mesh refinement (AMR) algorithms \cite{cao1999anr, carraturo2019suitably, kollmannsberger2018hierarchical,zhang20193,olleak2020efficient,gan2021benchmark}. In AMR, the mesh surrounding the melt pool is adaptively refined, while a coarser mesh is employed in regions distant from it. This approach can significantly reduce the total number of degrees of freedom (DoFs) compared to a uniformly fine mesh. However, AMR presents challenges in designing meshes with disparate resolutions, which can lead to ill-conditioning. Furthermore, frequent remeshing introduces significant computational overhead \cite{patil2015generalized}. 

An alternative to AMR involves multilevel methods, such as hierarchical bases \cite{yserentant1986multi, bank1988hierarchical} and multigrid techniques \cite{briggs2000multigrid, hackbusch2013multi}, which are characterized by a hierarchy of basis functions across the entire domain. This global approach, however, can be computationally expensive for large-scale problems with highly localized features. To address this limitation, another class of ``multilevel" methods has emerged that superimposes a local fine mesh onto a global coarse mesh to capture fine-scale physics. For example, Viguerie et al. \cite{viguerie2020fat} developed a two-level method inspired by the Fat Boundary Method \cite{maury2001fat} to reconcile localized material properties with the global problem. Carraturo et al. \cite{carraturo2022two} later applied a modified version to simulate part-scale LPBF printing. Henceforth, we adopt this latter definition: "multilevel" refers specifically to methods that combine a global coarse mesh with a local fine mesh, typically using two levels and low-order interpolations like linear FEM \cite{leonor2024go}.

A critical issue in multiscale methods is how to combine solutions with different resolutions seamlessly. The Variational MultiScale (VMS) proposed by Hughes et al. provided a foundational framework \cite{hughes1998variational} for this. Their work employs an overlapping sum decomposition to separate the fine and coarse scales. This strategy effectively captures subgrid terms and improves the overall accuracy of the final output.
Based on VMS, Leonor and Wagner \cite{leonor2024go} proposed the GO-MELT method. In this approach, different levels of meshes of varying resolutions are used to describe the solution: the coarsest mesh captures the global temperature variation, while subdomains with finer meshes track the moving laser source. Solutions at different scales are coupled through the variational weak form and interface conditions. GO-MELT and the existing two-level methods \cite{carraturo2022two, leonor2024go} employ linear finite elements as approximation functions at all levels. However, the requirement that the fine mesh must fully cover the grid points of the coarse mesh in the overlapping domains limits the use of high-order finite elements to achieve higher-order accuracy. Meshes at the boundary of each level must be carefully designed, as high-order finite elements require the introduction of additional internal nodes. Moreover, GO-MELT leverages an explicit time-stepping method and is limited by the small time steps required by the stability condition \cite{leonor2024go}. For a large-scale simulation, this may result in a significantly large number of time steps \cite{kopp2022space}.

To address the difficulty of using higher-order finite elements across different scales in the GO-MELT framework, we adopt the \textbf{C}onvolution \textbf{Hi}erarchical \textbf{De}ep-learning \textbf{N}eural \textbf{N}etwork (C-HiDeNN) as the efficient and accurate interpolation function \cite{lu2023convolution, park2023convolution, li2023convolution} in the \textbf{M}ulti\textbf{L}evel \textbf{V}ariational \textbf{M}ulti\textbf{S}cale framework (ML-VMS) as shown in Fig. \ref{fig:summary} (b). Unlike higher-order Lagrange elements in FEM, which require adding extra nodes to each element, C-HiDeNN leverages neighboring nodes to construct convolutional patch functions. For this reason, C-HiDeNN can use linear finite element meshes while achieving higher polynomial orders and accuracy without increasing the degrees of freedom (DoFs). As a result, C-HiDeNN can achieve the same level of accuracy but with much fewer DoFs compared to linear FEM \cite{park2023convolution}. This feature is particularly advantageous for multiscale modeling, as it enables the straightforward application of different orders of accuracy across various scales.

The high computational cost of ultra large-scale simulations can be substantially alleviated by using model order reduction techniques such as tensor decomposition (TD)\cite{zhang2022hidenn, guo2025interpolating, guo2024convolutional,lu2024extended}. The resulting method, C-HiDeNN-TD, has proven to be effective in substantially accelerating simulations in applications such as additive manufacturing \cite{guo2025tensor} and topology optimization \cite{li2023convolution}. A key feature of C-HiDeNN-TD is its ability to overcome the curse of dimensionality, as its number of unknowns scales linearly with the number of problem dimensions \cite{guo2025tensor}. As shown in Fig. \ref{fig:summary} (c), we integrate space-time C-HiDeNN-TD into the ML-VMS framework. Compared to existing spatiotemporal two-level methods \cite{viguerie2022spatiotemporal} and GO-MELT \cite{leonor2024go}, ML-VMS with space-time C-HiDeNN-TD  formulation automatically realizes multiple levels of temporal discretizations and allows the use of larger time steps than explicit time-stepping schemes (e.g., the forward Euler scheme in GO-MELT). When combined with higher-order temporal elements, this approach substantially reduces the number of steps required for an accurate long-duration analysis. As a result, ML-VMS C-HiDeNN-TD has exceptional performance in terms of speed, accuracy, resolution, memory usage and storage requirement, as summarized in Fig. \ref{fig:summary}(d).

In this paper, we propose the Multilevel Variational MultiScale (ML-VMS) method, which integrates the multilevel mesh into the VMS framework while leveraging the flexibility of C-HiDeNN basis functions and speedup through tensor decomposition. To this aim, we first review the basics of C-HiDeNN interpolation theory in Section \ref{sec:CHIDENN}. The general $m$-level ML-VMS framework and its corresponding solution schemes are introduced in Section \ref{sec:MLVMS}. We also provide a rigorous error analysis of ML-VMS for elliptical problems. We extend ML-VMS with the space-time C-HiDeNN-TD approximation to further speed up the computation in Section \ref{sec:CHIDENNTD}. Then we investigate the performance of ML-VMS numerically in Section \ref{sec:NUMERICAL} and discuss the technical details in Section \ref{sec:discussions}. Finally, in Section \ref{sec:CONCLUSION} we discuss the potential improvements
of the current approach in the near future.
% and highlight the estimation of optimal C-HiDeNN hyperparameters for multilevel interpolations

\section{Review of C-HiDeNN}
\label{sec:CHIDENN}
\subsection{C-HiDeNN interpolations}
\begin{figure}[htbp]
\centering
\includegraphics[width=6in]{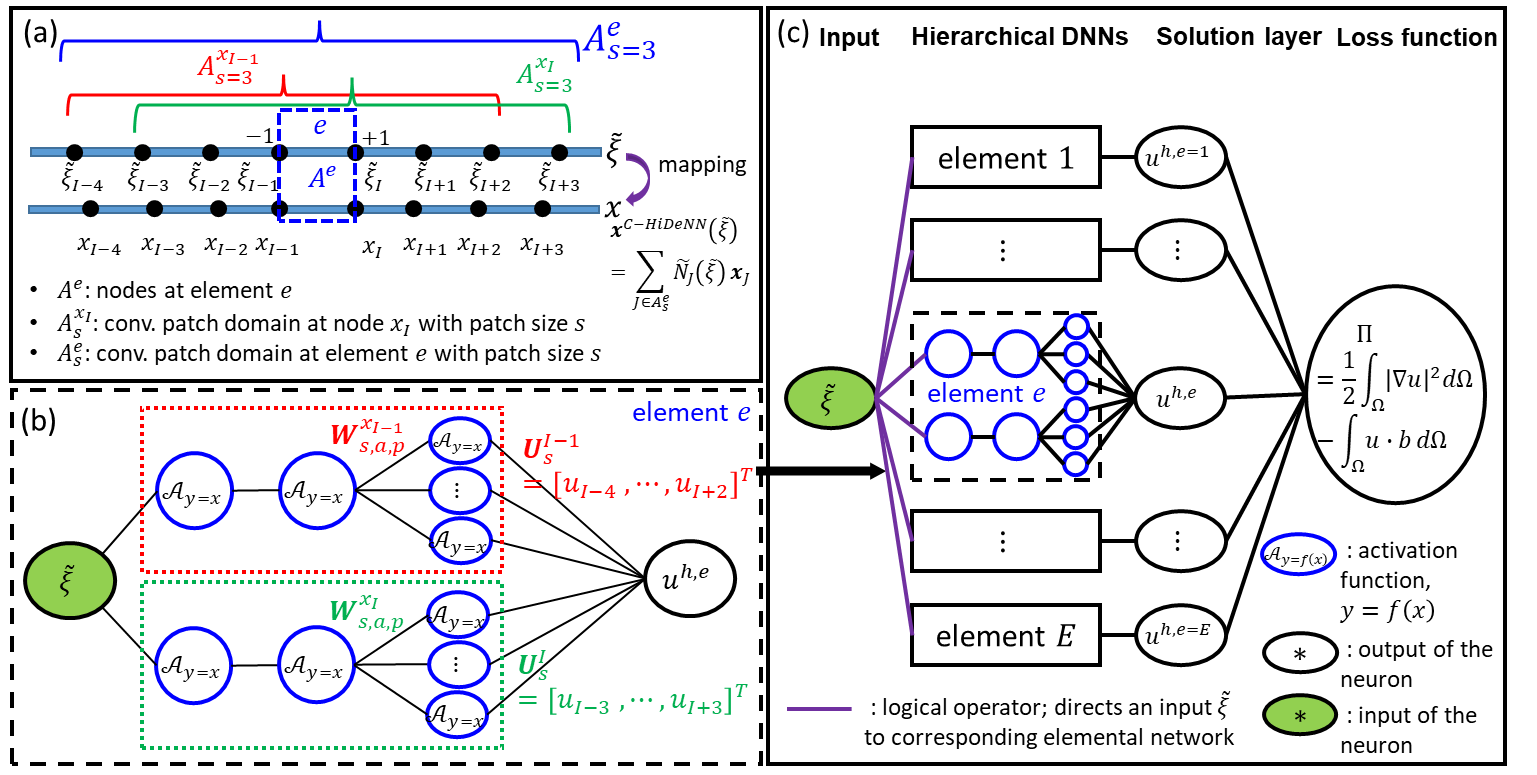}
\caption{Neural network structure of Convolution Hierarchical Deep-learning Neural Network (C-HiDeNN) for solving one-dimensional (1D) Poisson’s equation. This figure is borrowed from \cite{zhang2025multi}.}
\label{fig:C-HiDeNN}
\end{figure}

Convolution Hierarchical Deep-learning Neural Network (C-HiDeNN) \cite{lu2023convolution} is a pruned neural network for high-order interpolation. The neural network structure of C-HiDeNN is illustrated in Figure 1. C-HiDeNN interpolation at coordinate $x$ of element $e$ is expressed as:
\begin{equation}
    u^{h,e}(\bm{x})=\mathcal{I}u(\bm{x})=\sum_{I\in A^e} N_I(\bm{x}) \sum_{K\in A_s^{\bm{x}_I}}W^{x_I}_{s,a,p,K}(\bm{x})u_K = \sum_{J\in A^e_s} \tilde{N}_J(\bm{x})u_J,
\end{equation}
where $\mathcal{I}$ denotes the C-HiDeNN interpolation operator, and $u_J=u(\bm{x}_J)$ is the nodal value of $u(\bm{x})$ at node $J$.
 $N_I(\bm{x})$ can be an arbitrary function with compact support associated with node $I$ in element $e$, satisfying the partition of unity, such as linear finite element (FE) shape functions. In this work, we specifically use linear FE shape functions as $N_I(\bm{x})$. Unlike linear shape functions, which are influenced only by the nearest neighbor nodes, the construction of $W^{x_I}_{s,a,p,K}$ in C-HiDeNN also depends on the nodes in the neighboring elements. The domain of these neighboring elements is referred to as a "convolution patch" with a convolution patch size $s$. $A^e$ denotes the set of nodes in element $e$. $A_s^{\bm{x}_I}$ is the set of nodes within the $s$ layers of elements that surround the center node $I$. The support domain of element $e$ (i.e., $A_s^e$ ) is the union of all nodal convolution patches: $A_s^e=\cup_{I\in A^e} A_s^{\bm{x}_I}$. In the structure of the neural network, convolution patch size $s$ determines the number of neighbors to build the third hidden layer, as shown in Fig. \ref{fig:C-HiDeNN}(b) whose weights are the convolutional functions of the convolution patch $\bm{W}_{s,a,p}^{\bm{x}_I}  (x)=\{W_{s,a,p,K}^{\bm{x}_I} (\bm{x})|K\in A_s^{\bm{x}_I} \}$.  
 
 The convolution patch function $\bm{W}_{s,a,p}^{\bm{x}_I} (\bm{x})$ can be constructed based on meshfree interpolation theories \cite{liu1995reproducing, chen1996reproducing}. A detailed construction procedure based on radial basis functions with dilation parameter $a$ is provided in \ref{appenx:ConvolutionPatch}. These convolution patch functions depend on 3 hyperparameters: the dilation parameter $a$ that determines the influence domain of the radial basis function, the reproducing polynomial order $p$, and the size of the convolution patch $s$.
 
 The nodal shape function $\tilde{N}_J(\bm{x})$ is expressed as the sum of the product of $N_I(\bm{x})$ and $W^{x_I}_{s,a,p,J}(\bm{x})$ over nodes $x_I\in\{x_I|J\in A^{\bm{x}_I}_s\}$, i.e.,
 \begin{equation}
     \tilde{N}_J(\bm{x}) = \sum_{x_I\in\{x_I|J\in A^{\bm{x}_I}_s\}} N_I(\bm{x})\cdot W^{x_I}_{s,a,p,J}(\bm{x}).
 \end{equation}
The C-HiDeNN interpolants can thus be reformulated following the standard procedure of interpolation: the nodal shape function $\tilde{N}_J(\bm{x})$ is multiplied by nodal variable $u_J$, followed by a summation over all nodes. Compared to standard FEM, which uses only the nodes inside an element to approximate a solution field within the element, C-HiDeNN incorporates nodes from neighboring elements. The support domain of the FEM shape function $N_I$ covers only one layer of elements surrounding the center node $I$, while the support domain of C-HiDeNN shape function $\tilde{N}_I$ covers $(s+1)$ layers of neighboring elements. The treatment of nodal convolution patch centered at the boundary node is detailed in \ref{appenx:BoundaryPatch}.

\subsection{Error estimate for C-HiDeNN}

After setting the reproducing polynomial order $p$ in the convolution patch functions $W_{s,a,p}^{\bm{x}_I} (\bm{x})$, C-HiDeNN interpolants satisfy the $p$-th order reproducing property:
\begin{equation}
    \mathcal{I} (\bm{x}^m) = \sum_{J} \tilde{N}_J(\bm{x}) (\bm{x}_J)^m = \bm{x}^m, m=0, 1,2,\ldots,p.
\end{equation}
To guarantee this reproducing property, we require at least $p+1$ nodes along one dimension in one nodal convolution patch $A_s^{\bm{x}_I}$. This leads to the condition $s \geq p/2$.

When applying C-HiDeNN interpolants to solve elliptic equations, the following error estimate holds \cite{lu2023convolution}:
\begin{equation} \label{eq:C-HiDeNNInterpError}
    \|u-u^h\|_{H^1}\leq C(a,s,p) h^p \|u\|_{H^{p+1}},
\end{equation}
where $C$ is independent of the mesh size $h$, $p$ is the reproducing polynomial order, $u^h$ is the approximate solution using C-HiDeNN, and $u$ is the exact solution to an elliptic equation. The coefficient $C$ depends on hyperparameters $a, s, p$. 

Notably, compared to FEM, C-HiDeNN can achieve higher polynomial orders and accuracy using only linear finite element meshes (e.g., 4-node tetrahedral elements or 8-node brick elements in 3D) without increasing the degrees of freedom (DoFs). A detailed comparison of C-HiDeNN with other numerical methods, including FEM, isogeometric analysis (IGA) \cite{hughes2005isogeometric} and physics-informed neural networks (PINN) \cite{raissi2019physics}, is presented in \ref{appenx:Comparison}. 

\section{ML-VMS with C-HiDeNN interpolations} 
\label{sec:MLVMS}
\subsection{Problem statement and weak form}

Consider an elliptic PDE:
\begin{equation} \label{eq:PDE}
    \mathcal{L} u(\bm{x})  = f(\bm{x}), \bm{x} \in \Omega
\end{equation}
with a Dirichlet boundary condition
\begin{equation} \label{eq:PDE_BC}
    u(\bm{x}) |_{\partial \Omega}=g.
\end{equation}

Let $\Omega \subset \mathbb{R}^d$, where $d \geq 1$ is the number of dimensions, be an open bounded domain with a smooth boundary $\partial \Omega$. The operator $\mathcal{L}$ is a differential operator. $f:\Omega \to \mathbb{R}$ and $g:\partial \Omega \to \mathbb{R}$ are source functions. 

Let $\mathcal{S} \subset H^1(\Omega)$ denote the trial function space, and $\mathcal{V} \subset H^1(\Omega) $ denote the weighting function space. The trial function $u\in\mathcal{S}$ satisfies the Dirichlet boundary condition, and the weighting function $w\in \mathcal{V}$ satisfies zero boundary condition. The weak form of this boundary-value problem is given as follows.

Find $u \in \mathcal{S}$ such that $\forall w\in \mathcal{V}$,
\begin{equation}
    a(w,u)_\Omega = (w,f)_\Omega,
\end{equation}
where $(\cdot,\cdot)_\Omega$ is the $L_2(\Omega)$ inner product, and $a(w,u)_\Omega$ is a bilinear form satisfying
\begin{equation} \label{eq:WeakForm}
    a(w,u)_\Omega=(w,\mathcal{L}u)_\Omega.
\end{equation}
Assuming that the energy norm $\|\cdot\|_{E(\Omega)}=(a(\cdot,\cdot)_\Omega )^{1/2}$ and $H^1$ norm $\|\cdot\|_{H^1 (\Omega)}$ are equivalent, we have
\begin{equation}
    c_1 \|v\|_{H^1 (\Omega)} \leq (a(v,v)_\Omega )^{1/2} \leq c_2 \|v\|_{H^1 (\Omega)}, \forall v \in H^1(\Omega),
\end{equation}
where $c_1$ and $c_2$ are positive constants. For FEM and C-HiDeNN with one single-level refinement over the computational domain, the solution is obtained by solving Eq. (\ref{eq:WeakForm}) directly.

\subsection{Two-level VMS}

In order to acheive high accuracy, a fine-level discretization over the entire computational domain can be extremely time-consuming. Nevertheless, a coarse mesh can result in faster speed, but will lead to inaccurate outcomes. An alternative strategy involves employing a fine mesh in the domains of interest while utilizing a coarse mesh in other parts. This leads to the multilevel method. However, this strategy raises several important questions to answer: How can meshes of different refinement levels be coupled seamlessly? How to determine the mesh sizes for different levels of meshes to balance efficiency and accuracy?

\begin{figure}[htbp]
\centering
\includegraphics[width=3in]{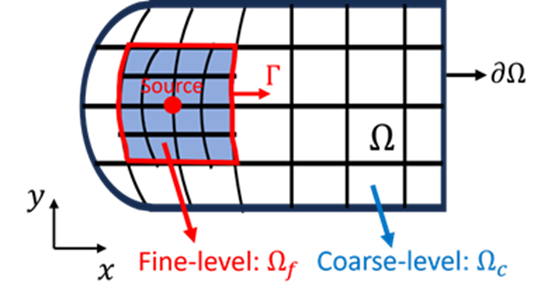}
\caption{The schematic of the two-level VMS method. }
\label{fig:2LevelVMS}
\end{figure}

We first consider the two levels of meshes: a fine-level mesh defined in a localized region and a coarse-level mesh defined throughout the entire domain. As illustrated in Fig. \ref{fig:2LevelVMS}, we define two overlapping domains with different levels of refinement: $\Omega_c$ is the full computational domain (i.e., $\Omega_c=\Omega$) with its boundary the same as $\partial \Omega$; $\Omega_f$ is a subdomain with its boundary $\Gamma$, serving as the interface between different levels. The subdomain $\Omega_f$ covers an integer number of coarse elements in $\Omega_c$, and is then refined with a finer mesh. Each coarse element in the overlapping zone can be represented by an integer number of finer elements from $\Omega_f$ mesh. The element size of the fine mesh is denoted as $h_f$, and the element size of the coarse mesh is $h_c$. For uniform multilevel meshes, the ratio $h_c/h_f$ is an integer denoted by $n$.

Let $\mathcal{I}^{(c)}$ be the $p_c$-th order C-HiDeNN interpolation operator defined on the coarse mesh in $\Omega_c$, given by
\begin{equation}
    \mathcal{I}^{(c)} u = \sum_{I_c} \tilde{N}_{I_c}(\bm{x}) u(\bm{x}_{I_c}),
\end{equation}
where $\tilde{N}_{I_c}$ is the C-HiDeNN shape function centered at the node $\bm{x}_{I_c}$ of the coarse mesh. Let $\mathcal{I}^{(f)}$ be the $p_f$-th order C-HiDeNN interpolation operator defined on the fine mesh in $\Omega_f$, given by
\begin{equation}
    \mathcal{I}^{(f)} u = \sum_{I_f} \tilde{N}_{I_f}(\bm{x}) u(\bm{x}_{I_f}),
\end{equation}
where $\tilde{N}_{I_f}$ is the C-HiDeNN shape function centered at the node $\bm{x}_{I_f}$ of the fine mesh with a smaller element size $h_f$. 

Let $\mathcal{S}_c^h \subset H^1 (\Omega)$ denote the trial function space for coarse interpolations, defined by
\begin{equation}
    \mathcal{S}_c^h = \{ v^h | v^h = \mathcal{I}^{(c)} v, v\in H^{p_c+1}(\Omega_c) \text{ and } v|_{\partial \Omega} = g\},
\end{equation}
associated with a weighting function space $\mathcal{V}_c^h \subset H^1 (\Omega)$ defined by
\begin{equation}
    \mathcal{V}_c^h = \{ w_c | w_c = \mathcal{I}^{(c)} w, w\in H^{p_c+1}(\Omega_c) \text{ and } w|_{\partial \Omega} = 0\}.
\end{equation}
Similarly, let $\mathcal{S}_f^h \subset H^1 (\Omega_f)$ denote the trial function space for fine-level interpolations, defined by
\begin{equation}
    \mathcal{S}_f^h = \{ v^h | v^h = \mathcal{I}^{(f)} v, v\in H^{p_f+1}(\Omega_f) \text{ and } v|_{\Gamma} = u_c\},
\end{equation}
associated with a weighting function space $\mathcal{V}_f^h \subset H^1 (\Omega_f)$ defined by
\begin{equation}
    \mathcal{V}_f^h = \{ w_f | w_f = \mathcal{I}^{(f)} w, w\in H^{p_f+1}(\Omega_f) \text{ and } w|_{\Gamma} = 0\}.
\end{equation}

We define the approximate solution $u^h(\bm{x})$ as follows:
\begin{equation} \label{eq:Two-LevelInterp}
    u^h(\bm{x}) = \left\{
    \begin{array}{cc}
        u_c(\bm{x}) & \bm{x}\in \Omega_c \backslash \Omega_f, \\
        u_f(\bm{x}) & \bm{x}\in \Omega_f,
    \end{array}
    \right.
\end{equation}
where $u_c (\bm{x})\in \mathcal{S}^h_c$ and $u_f (\bm{x})\in \mathcal{S}^h_f$ are interpolation functions defined on the coarse- and fine-level meshes, respectively. To ensure the continuity of $u^h$ across the entire domain, a Dirichlet boundary condition is imposed on the interface $\Gamma$: $u_f |_{\Gamma} = u_c$. In the solution scheme, this condition is enforced by specifying values at the nodes of the fine-mesh, specifically, $u_f(\bm{x}_{I_f}) = u_c(\bm{x}_{I_f}), \forall \bm{x}_{I_f} \in \Gamma$. Note that the approximation function $u_f$ defined in $\Omega_f$ is expressed by the fine-level interpolation operator $\mathcal{I}^{(f)}$ entirely, rather than being decomposed into the sum of a coarse-level approximation function defined by $\mathcal{I}^{(c)}$ and a fine-level approximation function defined by $\mathcal{I}^{(f)}$ like standard VMS \cite{hughes1998variational}.

We use the following weak forms to solve $u^h$:
\begin{eqnarray}
    a(w_c,u^h )_{\Omega_c} = (w_c,f)_{\Omega_c}, &&\forall w_c \in \mathcal{V}_c^h, \label{eq:Two-Level_Weak1}\\
    a(w_f,u^h )_{\Omega_f} =(w_f,f)_{\Omega_f}, &&\forall w_f \in \mathcal{V}_f^h, \label{eq:Two-Level_Weak2}\\
    \text{Boundary conditions: } u_c |_{\partial \Omega} = g,  && u_f |_\Gamma = u_c. \label{eq:Two-Level_BC}
\end{eqnarray}
The functions $u_c$ and $u_f$ are coupled within this system, following the framework of the VMS method \cite{hughes1998variational}, despite the fine-level domain $\Omega_f$ is only a subdomain of the entire domain. We refer to this framework as the two-level VMS. %If $\Omega_f = \Omega_c = \Omega$, the two-level VMS reduces to Hughes' VMS \cite{hughes1998variational}, in which the key assumption $u_f |_{\Gamma}=0$ becomes $u_f |_{\partial \Omega}=0$. 

\subsubsection{Comparison among VMS, GO-MELT and ML-VMS}

We compare different aspects of VMS, GO-MELT and ML-VMS in Table \ref{tab:CompareVMS_GOMELT}. VMS employs a fine mesh across the entire computational domain, whereas GO-MELT and ML-VMS only apply a fine mesh within a small fraction of the entire domain to cover the localized source. In VMS, subgrid terms can be effectively captured to enhance solution accuracy. Building on this foundation, GO-MELT and ML-VMS further achieve higher efficiency by leveraging local fine meshes.

\begin{table}
\begin{center}
\small
\begin{tabular}{|c|c|c|c|}

\hline
Methods & VMS & GO-MELT & Mutli-level VMS \\
\hline

Mesh 
& \begin{minipage}{0.25\textwidth}
    \centering
    \includegraphics[width=\linewidth]{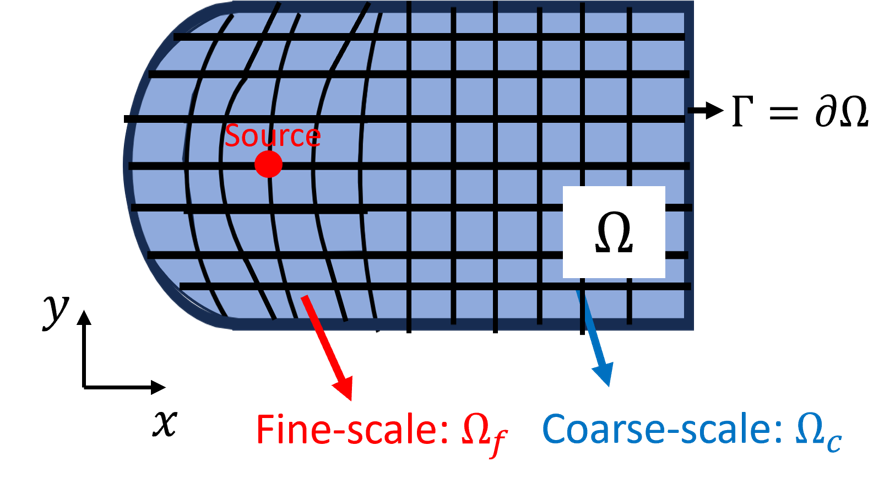}
  \end{minipage}
&
\begin{minipage}{0.25\textwidth}
    \centering
    \includegraphics[width=\linewidth]{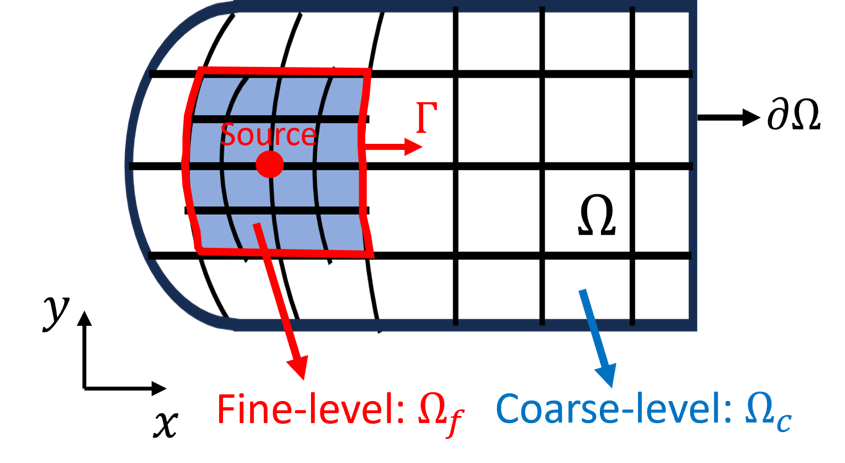}
  \end{minipage}
& \begin{minipage}{0.25\textwidth}
    \centering
    \includegraphics[width=\linewidth]{fig/GO-MELT.png}
  \end{minipage}
\\ 
& $\Omega_f=\Omega_c=\Omega$ & $\Omega_f\subset\Omega_c=\Omega$ & $\Omega_f\subset\Omega_c=\Omega$ \\ \hline
\makecell{Shape \\ functions} & FEM & FEM & \textbf{\makecell{C-HiDeNN \\ (arbitrary order)}} \\ \hline
\makecell{Approximation \\ solutions} 
& $u^h=u_c+u'$ on $\Omega$
& $u^h=\left\{
\begin{array}{cc}
    u_c & \text{on } \Omega \backslash \Omega_f \\
    u_c+u' & \text{on } \Omega_f
\end{array}
\right.$
& $u^h=\left\{
\begin{array}{cc}
    u_c & \text{on } \Omega \backslash \Omega_f \\
    u_f & \text{on } \Omega_f
\end{array}
\right.$

\\ \hline 

\makecell{Interpolation \\ on $\Omega_c$} 
& $u^h=\mathcal{I}^{(c)}u$
& $u^h=\mathcal{I}^{(c)}u$
& $u^h=\mathcal{I}^{(c)}u$
\\ \hline

\makecell{Interpolation \\ on $\Omega_f$} 
& $u^h=\mathcal{I}^{(c)}u+\mathcal{I}^{(f)}(u-\mathcal{I}^{(c)}u)$
& $u^h=\mathcal{I}^{(c)}u+\mathcal{I}^{(f)}(u-\mathcal{I}^{(c)}u)$
& $u^h=\mathcal{I}^{(f)}u$
\\ \hline
\makecell{Interpolation \\ Order} 
& \makecell{Hierarchical $p$-refinement \\ for $\mathcal{I}^{(f)}$ \\ with extra nodes}
& \makecell{Linear FEM \\ interpolation operator \\ $\mathcal{I}^{(c)}$ and  $\mathcal{I}^{(f)}$}
& \makecell{$\mathcal{I}^{(c)}$ \textbf{and} $\mathcal{I}^{(f)}$ \textbf{can be} \\ \textbf{of arbitrary orders} \\\textbf{with linear FE mesh}}
\\ \hline

\makecell{Interface \\ condition} & $u'|_\Gamma=0$ & $u'|_\Gamma=0$ & $u_f|_\Gamma=u_c$ \\ \hline
\makecell{Solution \\ scheme \\ (weak form)} 
& \makecell{$a(w_c, u_c+u')_{\Omega} = (w_c, f)_{\Omega}$ \\ $u'=G*u_c$ \\ ($G$: Green's function) } & \makecell{$a(w_c, u^h)_{\Omega} = (w_c, f)_{\Omega}$ \\ $a(w_f, u^h)_{\Omega_f} = (w_f, f)_{\Omega_f}$ } & \makecell{$a(w_c, u^h)_{\Omega} = (w_c, f)_{\Omega}$ \\ $a(w_f, u_f)_{\Omega_f} = (w_f, f)_{\Omega_f}$ } \\ \hline
Advantage & Accuracy & Accuracy and Efficiency & Accuracy and Efficiency \\ \hline

\end{tabular}
\caption{Comparison among VMS, GO-MELT and ML-VMS. $\mathcal{I}^{(c)}$ and $\mathcal{I}^{(f)}$ are interpolation operators defined in two levels of refinements, respectively. In VMS, the coarse- and fine-scales denoted by $u_c$ and $u'$, respectively, cover the whole domain and standard FEM basis function is used for $\mathcal{I}^{(c)}$ and $\mathcal{I}^{(f)}$. In GO-MELT, a local fine mesh is used to increase the interpolation accuracy in the incremental fine-scale component $u'$ using FEM basis with the same order as coarse-scale interpolators.  In ML-VMS, a local linear fine mesh is used in the fine-scale domain. But C-HiDeNN interpolation operator $\mathcal{I}^{(f)}$ offers fine-grained control on the order of approximation with linear mesh. In the fine mesh, ML-VMS uses a whole fine-level interpolation operator $\mathcal{I}^{(f)}$ to represent the approximation function $u_f$, instead of the sum of coarse- and fine-scale approximations represented by two levels of interpolation operators. }
\label{tab:CompareVMS_GOMELT}
\end{center}
\end{table}

ML-VMS differs from VMS and GO-MELT in two key aspects: the selection of basis functions and the approximation of a fine-scale solution. Both VMS and GO-MELT adopt FEM as their basis functions. However, when using higher-order elements, they require extra degrees of freedom (DoFs). More importantly, the interface between different levels must be carefully designed to ensure that the finer mesh in the overlapping domain fully covers the nodes of the coarser mesh. However, ML-VMS adopts C-HiDeNN as its approximation function, which allows shape functions of arbitrary orders at all levels, without increasing additional degrees of freedom (DoFs) or modifying the mesh. Both VMS and GO-MELT represent the fine solutions as the sum of the coarse-scale $u_c$ and the incremental fine-scale component $u'$. ML-VMS directly represents the fine solution using the fine-scale interpolator. From the perspective of interpolation theory, the interpolation of the exact solution corresponding to VMS and GO-MELT in the fine-level domain is $u^h=\mathcal{I}^{(c)}u+\mathcal{I}^{(f)}(u-\mathcal{I}^{(c)}u)$, while that corresponding to ML-VMS is $u^h=\mathcal{I}^{(f)}u$. If $\mathcal{I}^{(f)}\mathcal{I}^{(c)}=\mathcal{I}^{(c)}$, these two forms are equivalent. Therefore, when FEM is used as the shape function, ML-VMS can be degenerated to GO-MELT.

All three methods are based on a core assumption: the coarse-scale solution at the interface provides boundary conditions for the fine scale. Consequently, VMS and GO-MELT rely on $u'|_\Gamma=0$, while ML-VMS relies on  $u_f|_\Gamma=u_c$.

Finally, the solution schemes are different for these three methods. VMS employs the Green's function method, where the fine-scale solution is expressed in terms of the coarse-scale solution to avoid coupling. Both GO-MELT and ML-VMS use variational forms at different scales: GO-MELT solves for coarse-scale solution $u_c$ and the incremental fine-scale component $u'$, while ML-VMS solves directly for solution for each scale $u_c$ and $u_f$ separately.

\subsubsection{Alternating level algorithm}

% In this subsection, we present the detailed solution scheme for the coupled system (\ref{eq:Two-Level_Weak1}-\ref{eq:Two-Level_BC}).

As the coarse level information $u_c$ and fine-level information $u_f$ are coupled in Eqs. (\ref{eq:Two-Level_Weak1}-\ref{eq:Two-Level_BC}), we propose an alternating-level algorithm to solve this original coupled problem iteratively. As a result, each iteration involves solving a linear problem. Specifically, we solve $u_c$ with $u_f$ given, followed by solving $u_f$ with $u_c$ given. This iterative process continues until the variation of the solution vectors satisfies the convergence criteria. 

To start with, since the fine-level information is unknown at the beginning, so we set $u_f=0$. When solving $u_c$,  We can define a virtual $u_c$ in the overlapping domain $\Omega_f$: $\bar{u}_f = \mathcal{I}^{(c)} u_f$ representing coarse scales, and then define $u'_f=u_f-\bar{u}_f$ representing the difference between fine-level and coarse-level approximation in $\Omega_f$. Note that the difference $u'_f$ may not be reproduced by fine-level interpolations, especially for different orders in two levels, which is different from the linear FEM used in GO-MELT. Based on the above assumption, we solve $u_c$ over the entire domain with given $u'_f=u_f-\mathcal{I}^{(c)}u_f$:
\begin{equation}
    a(w_c,u_c )_{\Omega_c} + a(w_c,u_f-\mathcal{I}^{(c)}u_f )_{\Omega_f} = (w_c,f)_{\Omega_c}, \quad\forall w_c \in \mathcal{V}_c^h
\end{equation}
with boundary conditions $u_c |_{\partial \Omega} = g$.
After $u_c$ is obtained, we solve $u_f$ according to:
\begin{equation}
    a(w_f,u_f )_{\Omega_f} =(w_f,f)_{\Omega_f}, \quad\forall w_f \in \mathcal{V}_f^h
\end{equation}
with Dirichlet boundary conditions $u_f |_{\Gamma} = u_c$.
We continually update the solution $u_c$ and $u_f$ until the variance of the solution falls within the tolerance. The detailed solution scheme is summarized in the Algorithm \ref{alg:2_level}.

\begin{algorithm}[H]
\caption{Alternating level algorithm for 2 levels}
\label{alg:2_level}
\begin{algorithmic}[1]
    \State Set initial fine-level solution $u_f \gets 0$.
    
    \Repeat
        \State \textit{Solve for coarse-level solution $u_c$}
        \State Find $u_c$ with B.C. $u_c |_{\partial \Omega} = g$, such that for all $w_c \in \mathcal{V}_c^h$:
        $$ a(w_c, u_c)_{\Omega_c} + a(w_c, u_f - \mathcal{I}^{(c)}u_f)_{\Omega_f} = (w_c, f)_{\Omega_c} $$
        
        \State \textit{Solve for fine-level solution $u_f$}
        \State Find $u_f$ with Dirichlet B.C. $u_f |_{\Gamma} = u_c$, such that for all $w_f \in \mathcal{V}_f^h$:
        $$ a(w_f, u_f)_{\Omega_f} = (w_f, f)_{\Omega_f} $$
        
    \Until{the variance of the solution is within the desired tolerance.}

\end{algorithmic}
\end{algorithm}

\subsection{$m$-level VMS}

\begin{figure}[htbp]
\centering
\includegraphics[width=2.5in]{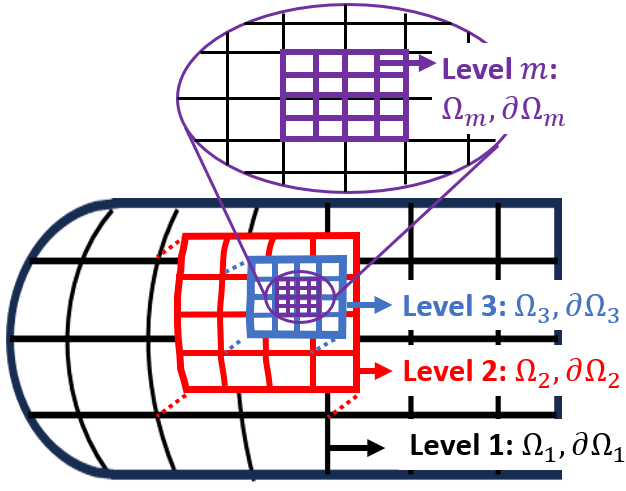}
\caption{The schematic of the $m$-level VMS method. }
\label{fig:nLevelVMS}
\end{figure}

The two-level VMS can be easily extended to an $m$-level VMS. In this extension, we define $m$ levels of overlapping domains that correspond to different refinements, denoted as $\Omega_1,\Omega_2, \ldots,\Omega_m$, as illustrated in Fig. \ref{fig:nLevelVMS}. The $1$st-level domain is the full domain $\Omega$, i.e., $\Omega_1=\Omega$. The domain of each subsequent level is a subdomain of the previous level such that $\Omega_{l+1}\subset \Omega_l, l=1,2,\cdots,m-1$. Each coarse element of the $l$-th level mesh in the overlapping zone $\Omega_{l+1}$ can be represented by an integer number of fine elements from the $(l+1)$-th level mesh. The boundary of $\Omega_l$ is denoted as $\partial \Omega_{l}$, which serves as the interface between the refinements of the $l$-th and $(l+1)$-th level.

The interpolation operator defined in the $l$-th level domain is given by
\begin{equation}
    \mathcal{I}^{(l)} u = \sum_{I_l} \tilde{N}_{I_l}(\bm{x}) u(\bm{x}_{I_l}),
\end{equation}
where $\tilde{N}_{I_l}$ is the C-HiDeNN shape function centered on the node $\bm{x}_{I_l}$ of the $l$-th level mesh with element size $h_l$. 

The $m$-level VMS interpolation is defined by:
\begin{equation}
    u^h(\bm{x})=\left\{
    \begin{array}{cc}
        u_1(\bm{x}), & \bm{x}\in \Omega_1 \backslash \Omega_2, \\
        u_2(\bm{x}), & \bm{x}\in \Omega_2 \backslash \Omega_3, \\
        \cdots & \\
        u_m(\bm{x}), & \bm{x}\in \Omega_m.
    \end{array}
    \right.
\end{equation}
with $u_l (\bm{x})\in\mathcal{S}_l^h, l=1,2,\ldots,m$ being the $l$-th level interpolation with an element size $h_l$. The corresponding trial function space $\mathcal{S}_l^h$ is defined by
\begin{equation}
    \mathcal{S}_l^h = \{ v^h | v^h = \mathcal{I}^{(l)} v, v\in H^{p_l+1}(\Omega_l) \text{ and } v|_{\partial \Omega_l} = u_{l-1}\},
\end{equation}
associated with a weighting function space $\mathcal{V}_l^h \subset H^1 (\Omega_l)$ defined by
\begin{equation}
    \mathcal{V}_l^h = \{ w_l | w_l = \mathcal{I}^{(l)} w, w\in H^{p_l+1}(\Omega_l) \text{ and } w|_{\partial \Omega_l} = 0\}.
\end{equation}
Here, we define $u_0=g$, representing Dirichlet boundary conditions. Note that $u_l (\bm{x})\in\mathcal{S}_l^h, l>1$ is assumed to satisfy $u_l|_{\partial \Omega_l}=u_{l-1}$ on $\partial \Omega_l$ to ensure continuity of $u^h(\bm{x})$ across the interface $\partial \Omega_l$. 
% The $1$st level trial function space is
% \begin{equation}
%     \mathcal{S}_1^h = \{ v^h | v^h = \mathcal{I}^{(1)} v, v\in H^{p_1+1}(\Omega) \text{ and } v|_{\partial \Omega} = g\},
% \end{equation}
% associated with a weighting function space $\mathcal{V}_1^h \subset H^1 (\Omega)$ defined by
% \begin{equation}
%     \mathcal{V}_1^h = \{ w_1 | w_1 = \mathcal{I}^{(1)} w, w\in H^{p_1+1}(\Omega) \text{ and } w|_{\partial \Omega} = 0\}.
% \end{equation}
The weak forms for $m$-level VMS are
\begin{equation} \label{eq:nlevel_weak}
    a(w_l,u^h )_{\Omega_l}=(w_l,f)_{\Omega_l}, \forall w_l\in \mathcal{V}_l^h, l=1,2,\ldots,m
\end{equation}
with boundary conditions
\begin{equation}
    u_1 |_{\partial \Omega} = g, u_l |_{\partial \Omega_l} = u_{l-1}, l=2,\ldots,m.
\end{equation}
Consequently, all $m$ levels of solutions are coupled through the boundary conditions in each level.

\subsubsection{Alternating-level algorithm for $m$-level VMS}

% In this subsection, we present the detailed solution scheme for the coupled system (\ref{eq:Two-Level_Weak1}-\ref{eq:Two-Level_BC}).

We extend the 2-level alternating-level algorithm to the general $m$-level case. Specifically, we solve one level with the other levels given. Initially, we set $u_l=0, l>1$ since the $l$-th ($l>1$) level information is unknown at the beginning. Then we solve for the solution in each different level subsequently. In the overlapping domain $\Omega_{k-1} \backslash \Omega_k (k>l)$, we can define a virtual $u_l$: $u_l = \mathcal{I}^{(l)} u_k$ representing $l$-th level information, and then the rest part $u'_{l,k}=u_k-\mathcal{I}^{l} u_k$ representing the finer-level information. Based on this assumption, we can solve $u_l$ over $\Omega_l$ including overlapping domains:
    \begin{equation}
        a(w_l,u_l)_{\Omega_l} + \sum_{k=l+1}^{m} a(w_l,u_k-\mathcal{I}^{(l)} u_k)_{\Omega_k \backslash \Omega_{k+1}} = (w_l,f)_{\Omega_l}, \quad\forall w_l \in \mathcal{V}_l^h
    \end{equation}
with boundary conditions $u_1 |_{\partial \Omega} = g$ for $l=1$ and $u_l |_{\partial \Omega_l} = u_{l-1}$ for $l>1$. Here, we define $\Omega_{m+1}=\emptyset$, yielding $\Omega_m \backslash \Omega_{m+1}=\Omega_m$. In the implementation, $u_l |_{\partial \Omega_l} = u_{l-1}$ is enforced by specifying values at the nodes of the finer ($l$-th level) mesh.
For the finest level, we solve the following equation:
    \begin{equation}
        a(w_m,u_m )_{\Omega_m} =(w_m,f)_{\Omega_m}, \quad\forall w_m \in \mathcal{V}_m^h
    \end{equation}
    with Dirichlet boundary conditions $u_m |_{\Omega_m} = u_{m-1}$.
We continually update the solution across all levels until the variance of the solution falls within the tolerance. The detailed solution scheme is summarized in the Algorithm \ref{alg:m_level}.

\begin{algorithm}[H]
\caption{Alternating level algorithm for $m$ levels}
\label{alg:m_level}
\begin{algorithmic}[1]
    \State \textit{Initialization}
    \For{$l \gets 2$ to $m$}
        \State Set $u_l \gets 0$.
    \EndFor
    
    \Repeat
        \State \textit{Solve for intermediate levels $l=1, \dots, m-1$}
        \For{$l \gets 1$ to $m-1$}
            \State Find $u_l \in \mathcal{V}_l^h$ by solving for all $w_l \in \mathcal{V}_l^h$:
            \begin{equation*}
                a(w_l, u_l)_{\Omega_l} = (w_l, f)_{\Omega_l} - \sum_{k=l+1}^{m} a(w_l, u_k - \mathcal{I}^{(l)} u_k)_{\Omega_k \setminus \Omega_{k+1}}
            \end{equation*}
            \State Apply boundary conditions:
            \If{$l=1$}
                \State $u_1 |_{\partial \Omega} = g$
            \Else
                \State $u_l |_{\partial \Omega_l} = u_{l-1}$
            \EndIf
        \EndFor
        
        \State \textit{Solve for the finest level $u_m$}
        \State Find $u_m \in \mathcal{V}_m^h$ with Dirichlet B.C. $u_m |_{\Omega_m} = u_{m-1} $ such that for all $w_m \in \mathcal{V}_m^h$:
        $$ a(w_m, u_m)_{\Omega_m} = (w_m, f)_{\Omega_m} $$
        
    \Until{the variance of the solution is within the desired tolerance.}
    
\end{algorithmic}
\end{algorithm}

\section{Theoretical error estimate for elliptic problems}

The C-HiDeNN interpolator with the hyperparameter $p$ can achieve a $p$-th order reproducing property, and a $p$-th order ($H^1$-norm) interpolation estimate. When solving the elliptic PDEs using single-level C-HiDeNN, a $p$-th order (energy- or $H^1$-norm) convergence rate can be observed. When a two-level VMS is applied with each level approximated using the C-HiDeNN basis function, the approximation error is estimated using the following theorem.

\begin{theorem}
\textbf{Error estimate for two-level VMS}. Let the $p_c$-th order interpolation estimate (reproducing property) for any interpolation $u_c \in \mathcal{S}_c^h$ and the $p_f$-th order interpolation estimate (reproducing property) for any interpolation $u_f\in \mathcal{S}_f^h$ hold. Then the following error estimate holds:
\begin{equation} \label{eq:ErrorEst_TwoLevel}
    \|u-u^h \|_{H^1 (\Omega) }\leq C_1  h_c^{p_c}  \|u\|_{H^{p_c+1} (\Omega_c \backslash  \Omega_f ) }+C_2  h_f^{p_f}  \|u\|_{H^{p_f+1} (\Omega_f ) }+ \left(C_3  h_c^{p_c } \|u\|_{H^{p_c+1} (\Gamma) }+C_4 h_f^{p_f}  \|u\|_{H^{p_f+1} (\Gamma ) } \right),
\end{equation}
where $u$ is the exact solution to Eq. (\ref{eq:PDE}) under the Dirichlet boundary condition (\ref{eq:PDE_BC}). $u^h$ is the approximate solution defined by the two-level interpolation in Eq. (\ref{eq:Two-LevelInterp}), which satisfies the weak forms as shown in Eqs. (\ref{eq:Two-Level_Weak1}, \ref{eq:Two-Level_Weak2}) and the boundary conditions Eq. (\ref{eq:Two-Level_BC}). We also assume that $u$ is sufficiently smooth.
\end{theorem}

\begin{corollary}
    If the key assumption $u_f |_\Gamma = u_c$ is satisfied exactly (i.e., the solution can be resolved by coarse-level interpolants on $\Gamma$), the following error estimate holds:
\begin{equation} 
    \|u-u^h \|_{H^1 (\Omega) }\leq C_1  h_c^{p_c}  \|u\|_{H^{p_c+1} (\Omega_c \backslash  \Omega_f ) } + C_2  h_f^{p_f}  \|u\|_{H^{p_f+1} (\Omega_f ) }.
\end{equation}
\end{corollary}

\begin{corollary}
    If $\Omega_f=\Omega_c=\Omega$,  the following error estimate holds:
\begin{equation} \label{eq:ErrorEst_VMS}
    \|u-u^h \|_{H^1 (\Omega) }\leq C_2  h_f^{p_f}  \|u\|_{H^{p_f+1} (\Omega ) }+ \left(C_3  h_c^{p_c } \|u\|_{H^{p_c+1} (\partial \Omega) }+C_4 h_f^{p_f}  \|u\|_{H^{p_f+1} (\partial \Omega ) } \right).
\end{equation}
\end{corollary}

\begin{corollary}
    If $\Omega_f=\Omega_c=\Omega$ and $u_f |_{\partial \Omega} = u_c$ is satisfied exactly (i.e., $g$ can be resolved by coarse-level interpolants on $\partial \Omega$, such as vanishing boundary conditions $u|_{\partial \Omega}=g=0$), the following error estimate holds:
\begin{equation} 
    \|u-u^h \|_{H^1 (\Omega) }\leq C_2  h_f^{p_f}  \|u\|_{H^{p_f+1} (\Omega_f ) }.
\end{equation}
\end{corollary}

On the right-hand side of inequality (\ref{eq:ErrorEst_TwoLevel}), the first term refers to the error contribution from the coarse scale; the second term represents the error contribution from the fine scale; the third term accounts for the assumption $u_f|_{\Gamma}=u_c$ deviating from the fine-scale information of the exact solution at the interface $\Gamma$ between two levels of refinements. 

If the key assumption $u_f |_\Gamma = u_c$ is satisfied exactly, the third term vanishes, leading to the conclusion of Corollary 1. In this case, the solution can be resolved by coarse-level interpolations $u_c$ on $\Gamma$, e.g., the exact solution vanishes on $\Gamma$. We can obtain the error estimate for standard VMS under the condition $\Omega_f=\Omega_c=\Omega$ (Corollary 2) \cite{hughes1998variational}, in which the first error term in the original formulation vanishes. If the assumption $u' |_{\partial \Omega} = 0$  (corresponding to $u_f|_{\Gamma}=u_c$) holds exactly for original VMS (e.g., under vanishing boundary conditions) \cite{hughes1998variational}, the error is solely constrained by the second error term, which indicates a $p_f$-th order convergence rate (Corollary 3).

Furthermore, the error estimate given in Eq. (\ref{eq:ErrorEst_TwoLevel}) can be simplified to
\begin{equation} \label{eq:ErrEst_reduced}
    \|u-u^h \|_{H^1 (\Omega) }\leq C^{(c)} \cdot h_c^{p_c} + C^{(f)} \cdot h_f^{p_f},
\end{equation}
where $C^{(c)}=C_1  \|u\|_{H^{p_c+1} (\Omega_c \backslash  \Omega_f )} + C_3  \|u\|_{H^{p_c+1} (\Gamma) }$ is a function of coarse-level hyperparameters $a_c, s_c, p_c$, and $C^{(f)}=C_2    \|u\|_{H^{p_f+1} (\Omega_f ) }+C_4   \|u\|_{H^{p_f+1} (\Gamma ) }$ is a function of fine-level hyperparameters $a_f, s_f, p_f$. The term $\|u\|_{H^{p_c+1} (\Omega_c \backslash  \Omega_f )}$ in $C^{(c)}$ and the term $\|u\|_{H^{p_f+1} (\Omega_f ) }$ in $C^{(f)}$ indicate that the characteristics of the exact solution $u$ significantly influences these two coefficients and the ratio between them.  Unlike the single-level C-HiDeNN error, which is based solely on single-level mesh size, the error in the two-level VMS with C-HiDeNN is governed by the refinements of both two levels of meshes. As a result, the convergence rate is determined by the dominant term between $C^{(c)} \cdot h_c^{p_c}$ and $C^{(f)} \cdot h_f^{p_f}$. This dependency complicates the selection of parameters to achieve an optimal balance between accuracy and efficiency.

When extending the two-level VMS to an $m$-level VMS, we can derive the following error estimate:

\begin{theorem} \label{theorem2}
\textbf{Error estimate for $m$-level VMS}. Let the $p_l$-th order interpolation estimate (reproducing property) for any $l$-th level interpolation $u_l \in \mathcal{S}_l^h$ ($l=1,2,\ldots,n$) hold. Then the following error estimate holds:
\begin{equation} \label{eq:ErrEst_nLevel}
    \|u-u^h \|_{H^1 (\Omega) }\leq \sum_{l=1}^m C^{(l)} (a_l,s_l,p_l;u,\Omega_l,\Omega_{l+1}, \partial \Omega_l, \partial \Omega_{l+1}) \cdot h_l^{p_l},
\end{equation}
where $u$ is the exact solution to Eq. (\ref{eq:PDE}) under the Dirichlet boundary condition (\ref{eq:PDE_BC}); $u^h$ is the approximate solution. We assume that $u$ is sufficiently smooth.
\end{theorem}

As can be seen from the equation, the error estimate comprises $m$ terms, each corresponding to one of all levels of refinements. The coefficient for each term,  $C^{(l)}$, is a function of $a_l, s_l, p_l$, tailored to a specific problem and geometry for the $m$ levels of meshes ($\Omega_1,\Omega_2,\ldots,\Omega_m$). Similar to the error estimate for two-level VMS, the coefficient $C^{(l)}$ is directly influenced by the local regularity of the solution, as it contains the norm of $u$ in $\Omega_l \backslash  \Omega_{l+1}$, specifically $\|u\|_{H^{p_l+1} (\Omega_l \backslash  \Omega_{l+1} )}$. It also encapsulates the errors arising from assumptions $u_l|_{\Omega_l}=u_{l-1}$, deviating from the true solution. Consequently, the behavior of the error is highly sensitive to the characteristics of the exact solution $u$ across the domain. %We will discuss the optimal element size ratio in the followin subsection, and study the convergence rate numerically in the numerical examples.

\begin{remark}
    \textbf{Estimated optimal element size ratios for given parameters based on Theorem 2.} For given parameters $a_l, s_l , p_l, l=1,2,\cdots,n$, the optimal set of element sizes $\{h_1, h_2, \ldots, h_n \}$ holds for
    \begin{equation}
        C^{(1)} \cdot h_1^{p_1} \approx C^{(2)} \cdot h_2^{p_2} \approx \cdots C^{(n)} \cdot h_n^{p_n},
    \end{equation}
yielding the optimal element size ratio
\begin{equation}
    h_l/h_1 = n_{l, opt} = \left\lceil \left( C^{(l)}/C^{(1)} \right)^{1/p_1} h_l^{(p_l-p_1)/p_1} \right\rceil, l=1,2,\cdots,n.
\end{equation}
\end{remark}

Here,  we apply the ceiling function to obtain an integer-valued element size ratio. Detailed discussions on this remark are presented in \ref{appenx:OptElemRatio}. In the special case of a two-level analysis, this optimal ratio reduces to:
\begin{equation}
    n_{opt} = \left\lceil \left( C^{(f)}/C^{(c)} \right)^{1/p_f} h_c^{(p_f-p_c)/p_f} \right\rceil.
\end{equation}

\section{Multilevel space-time VMS with C-HiDeNN-TD}
\label{sec:CHIDENNTD}
When modeling problems that are large-scale in nature, we can further reduce the computational cost of ML-VMS with C-HiDeNN using model order reduction.  In this paper, we use the C-HiDeNN-Tensor Decomposition (TD) method, which leverages a special form of tensor decomposition called canonical polyadic decomposition \cite{guo2025tensor}, as the reduced-order model within the ML-VMS framework. For convenience, in the following discussions, $m$-level C-HiDeNN and $m$-level C-HiDeNN-TD denote $m$-level VMS with C-HiDeNN and C-HiDeNN-TD serving as interpolation functions, respectively. In C-HiDeNN-TD, we employ a space-time formulation rather than a conventional time-stepping scheme. Consequently, the problem is treated as a four-dimensional (4D) system (3D space and 1D time). In single-level C-HiDeNN-TD,  the space-time solution can be approximated as:
\begin{equation} \label{eq:TDform1}
    \mathcal{I}^{TD} u = \sum_{q=1}^{Q} \bm{\tilde{N}}^{1D}_x(x) \bm{u}^{(q)}_x \cdot \bm{\tilde{N}}^{1D}_y(y) \bm{u}^{(q)}_y \cdot \bm{\tilde{N}}^{1D}_z(z) \bm{u}^{(q)}_z,
\end{equation}
for 3D spatial problems, and
\begin{equation} \label{eq:TDform}
    \mathcal{I}^{TD} u = \sum_{q=1}^{Q} \bm{\tilde{N}}^{1D}_x(x) \bm{u}^{(q)}_x \cdot \bm{\tilde{N}}^{1D}_y(y) \bm{u}^{(q)}_y \cdot \bm{\tilde{N}}^{1D}_z(z) \bm{u}^{(q)}_z \cdot \bm{\tilde{N}}^{1D}_t(t) \bm{u}^{(q)}_t,
\end{equation}
for 3D time-dependent problems,
where $Q$ is the mode number; $\bm{\tilde{N}}^{1D}_d=[\tilde{N}^{1D}_{d1}(d), \tilde{N}^{1D}_{d2}(d), \ldots, \tilde{N}^{1D}_{dn_d}(d)]$ is the 1D convolution shape functions for the $d$-th dimension; $\bm{u}^{(q)}_d=[u^{(q)}_{d1}, u^{(q)}_{d2}, \ldots, u^{(q)}_{dn_d}]$ is the nodal solution vector for the $q$-th mode and $d$-th dimension. 

Due to TD, the original 4D space-time problem is approximated by a finite sum of multiplications of 1D problems. As a result, the total number of unknowns to be solved decreases from $n_x \times n_y \times n_z \times n_t$ to $(n_x+n_y+n_z+n_t)×Q$, where $n_d$ refers to the number of nodes in the $d$-th dimension. For large-scale problems where $n_d$ can be extremely huge, C-HiDeNN-TD can achieve substantial computational acceleration. Note that a space-time mesh is utilized instead of conventional time-stepping schemes for time-dependent problems.

By integrating the C-HiDeNN-TD approximation into the $m$-level VMS, the $l$-th level space-time solution $u_l(\bm{x},t)$ can be approximated as:
\begin{equation}
    u_l(\bm{x},t)=\sum_{q=1}^{Q_l} \bm{\tilde{N}}^{1D,(l)}_x(x) \bm{u}^{(q),(l)}_x \cdot \bm{\tilde{N}}^{1D,(l)}_y(y) \bm{u}^{(q),(l)}_y \cdot \bm{\tilde{N}}^{1D,(l)}_z(z) \bm{u}^{(q),(l)}_z \cdot \bm{\tilde{N}}^{1D,(l)}_t(t) \bm{u}^{(q),(l)}_t.
\end{equation}
In accordance with the assumption $u_l |_{\partial \Omega_l} = u_{l-1}$ for $l>1$ in ML-VMS, we require that mode number for Level $l$ is larger than that for Level $l-1$, i.e., $Q_1<Q_2<Q_3<\cdots<Q_m$. The first $Q_{l-1}$ modes in $u_l, l>1$ are obtained by interpolating $u_{l-1}$ mode by mode:
\begin{equation}
    \bm{u}^{(q),(l)}_d |_{\partial \Omega_l} = \bm{\tilde{N}}^{1D,(l-1)}_d(d) \bm{u}^{(q),(l-1)}_d, q=1,2,\cdots,Q_{l-1}
\end{equation}
for dimension $d$.
The remaining $Q_{l}-Q_{l-1}$ modes are solved and satisfy zero boundary conditions.

\subsection{Error analysis for elliptic problems with C-HiDeNN-TD}

\begin{theorem}
\textbf{Error decomposition for elliptic problems with C-HiDeNN-TD}. Let $u, u^{C-HiDeNN-TD}$,$u^{C-HiDeNN}$ be the exact solution, the numerical solution of C-HiDeNN-TD and C-HiDeNN  to an elliptic problem, respectively. Both C-HiDeNN-TD and C-HiDeNN utilize the same discretization. Therefore, the C-HiDeNN shape function can be constructed as a product of the 1D shape functions from C-HiDeNN-TD in each dimension. Then the following error decomposition for the elliptic problem holds:
\begin{equation}\label{eq:ErrorDecomposition}
    \left(\| u^{\text{C-HiDeNN-TD}}-u \|_E \right)^2 = \left( \| u^{\text{C-HiDeNN}}-u \|_E \right)^2 + \left( \| u^{\text{C-HiDeNN-TD}}-u^{\text{C-HiDeNN}} \|_E \right)^2.
\end{equation}
\end{theorem}

This theorem reveals that the C-HiDeNN-TD error decomposes into two parts: a discretization error and a decomposition error. The former depends solely on the mesh refinement (and can be estimated via C-HiDeNN error analysis), while the latter is determined by the number of modes used in tensor decomposition. 

The error estimate for the $m$-level VMS with C-HiDeNN-TD reduces to that of the standard C-HiDeNN framework shown in Theorem \ref{theorem2} under one key assumption: a sufficiently large mode number $Q$ is used in TD at each level. An upper bound for the required mode number has previously been established in the reference \cite{zhang2022hidenn}. With a sufficiently large mode number, the decomposition error term vanishes, allowing the C-HiDeNN-TD method's accuracy to match that of the standard C-HiDeNN.

The upper bound on the number of modes proved in \cite{zhang2022hidenn} guarantees that the error between the C-HiDeNN-TD and C-HiDeNN solutions vanishes entirely. In practice, however, far fewer modes are often sufficient to achieve a satisfactory approximation (e.g., an error below $10^{-6}$). The primary challenge is that the rate of this error reduction is problem-dependent, which complicates the a priori estimation of a sufficient mode number. We address this issue and discuss methods for estimating the required number of modes in \ref{appenx:estimateModes}.

\section{Numerical examples}
\label{sec:NUMERICAL}
In this section, we first verify our theory and study the effect of parameters ($a, s, p$ in each level and element size ratio) using 2D Poisson's equation and 1D heat equation, and demonstrate the efficiency of multi-level VMS with C-HiDeNN-TD on the high-dimensional problems, specifically, 3D transient heat transfer simulation with moving heat source under the setting of the real LPBF simulation. All the numerical simulation are run using MATLAB with a CPU of Intel Core i7-11700F.

\subsection{Two-level VMS for the elliptic problem}

Consider a 2D Poisson’s equation defined within $\Omega=[0,20]\times[0,20]$:
\begin{equation} \label{eq:Eg1_Eqn}
    \Delta u(x,y) + b(x,y) = 0.
\end{equation}
where $\Delta$ is the Laplace operator. We manufacture an exact solution as follows:
\begin{equation}
    u(x,y)=\sum_{k=1}^7 \exp \left(-\pi(x-8.2-0.2k)^2 - \pi(y-8.2-0.2k)^2\right).
\end{equation}

The source term $b(x,y)$ is then derived from Eq. (\ref{eq:Eg1_Eqn}) using the method of manufactured solutions. The boundary conditions are also determined from the exact solution.

\begin{figure}[htbp]
\centering
\subfigure[Two-level mesh]{\includegraphics[width=3.2in]{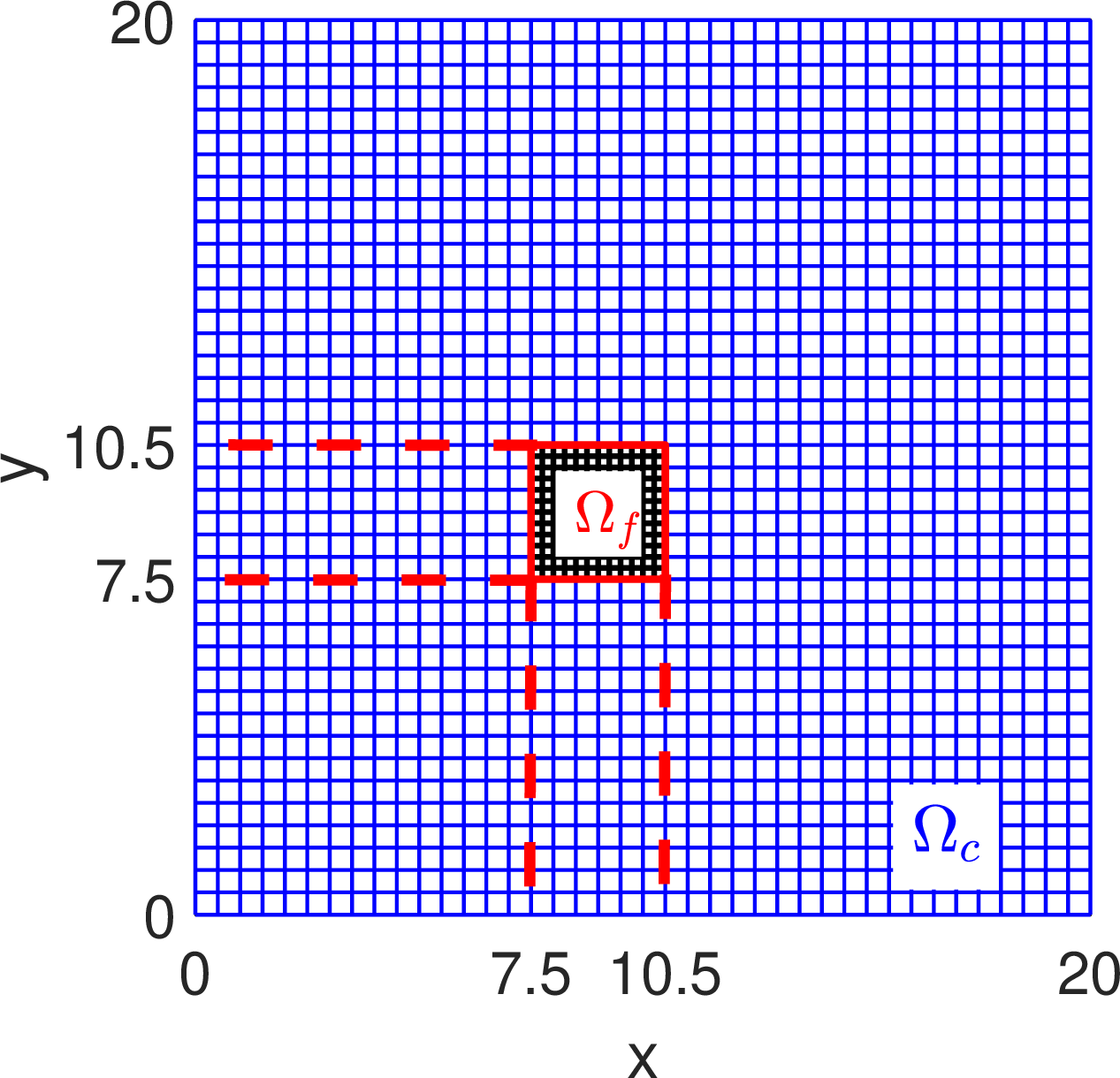}}
\subfigure[Exact solution]{\includegraphics[width=3.2in]{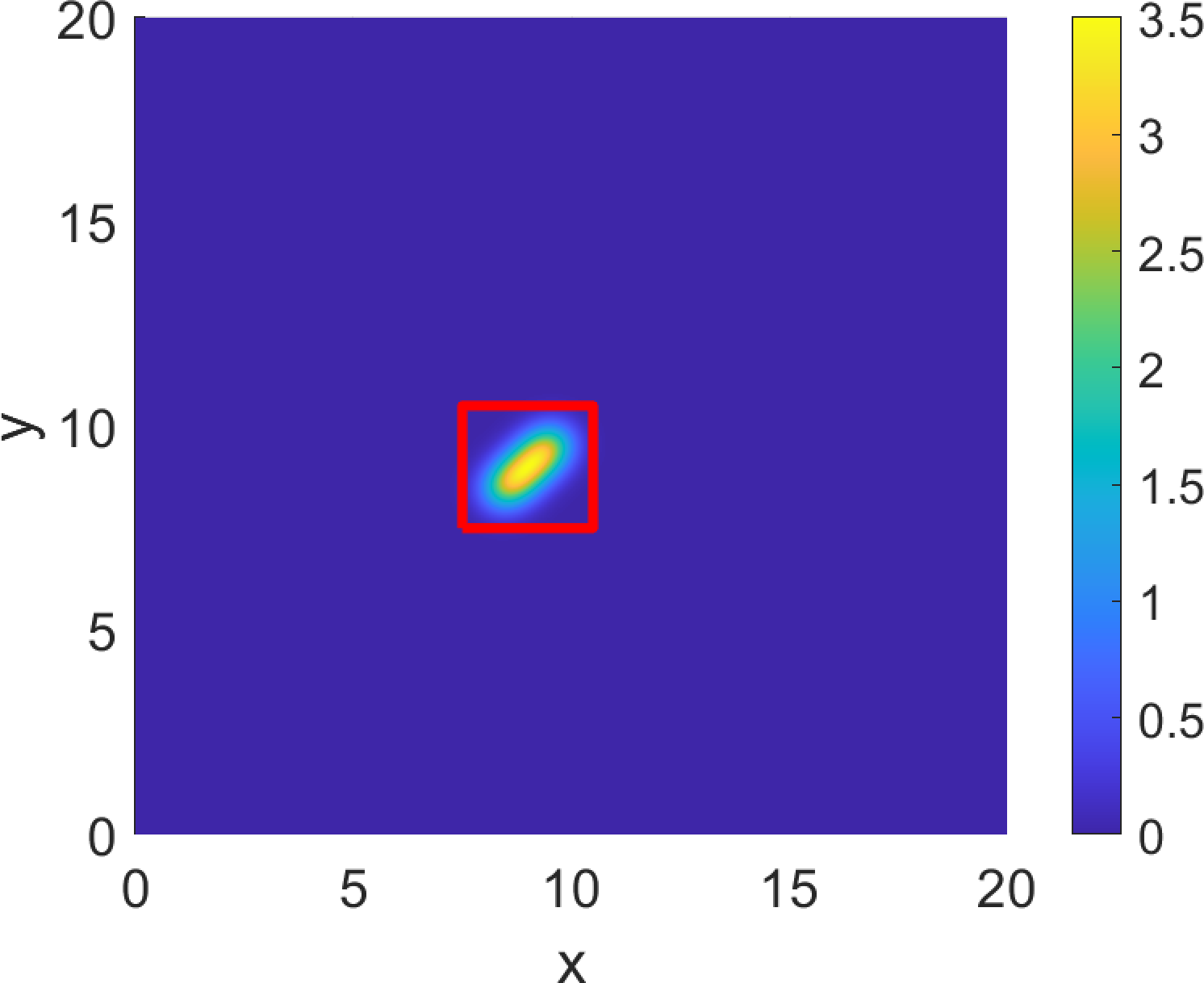}}
\caption{The schematic of the two-level mesh and the exact solution. }
\label{fig:Example1}
\end{figure}

Figure \ref{fig:Example1} shows the 2-level mesh and the exact solution. The fine-level mesh is used in subdomain $\Omega_f=[7.5,10.5]\times[7.5,10.5]$. 

\subsubsection{Two-level VMS with C-HiDeNN}

We apply the two-level VMS with C-HiDeNN (short for two-level C-HiDeNN) to the same problem. We numerically verify our theory for error estimation, and then study the optimal hyperparameters including the mesh size ratio $n=h_c/h_f$ and the polynomial order for each level to balance the tradeoff between efficiency and accuracy. Details refer to \ref{appenx:two-levelVMS_elliptic}.

% The optimal accuracy for various cases at similar costs is closely aligned. The optimal element size ratio is approximately between 2 and 4. 

\begin{figure}[htbp]
\centering
\includegraphics[width=0.6\textwidth]{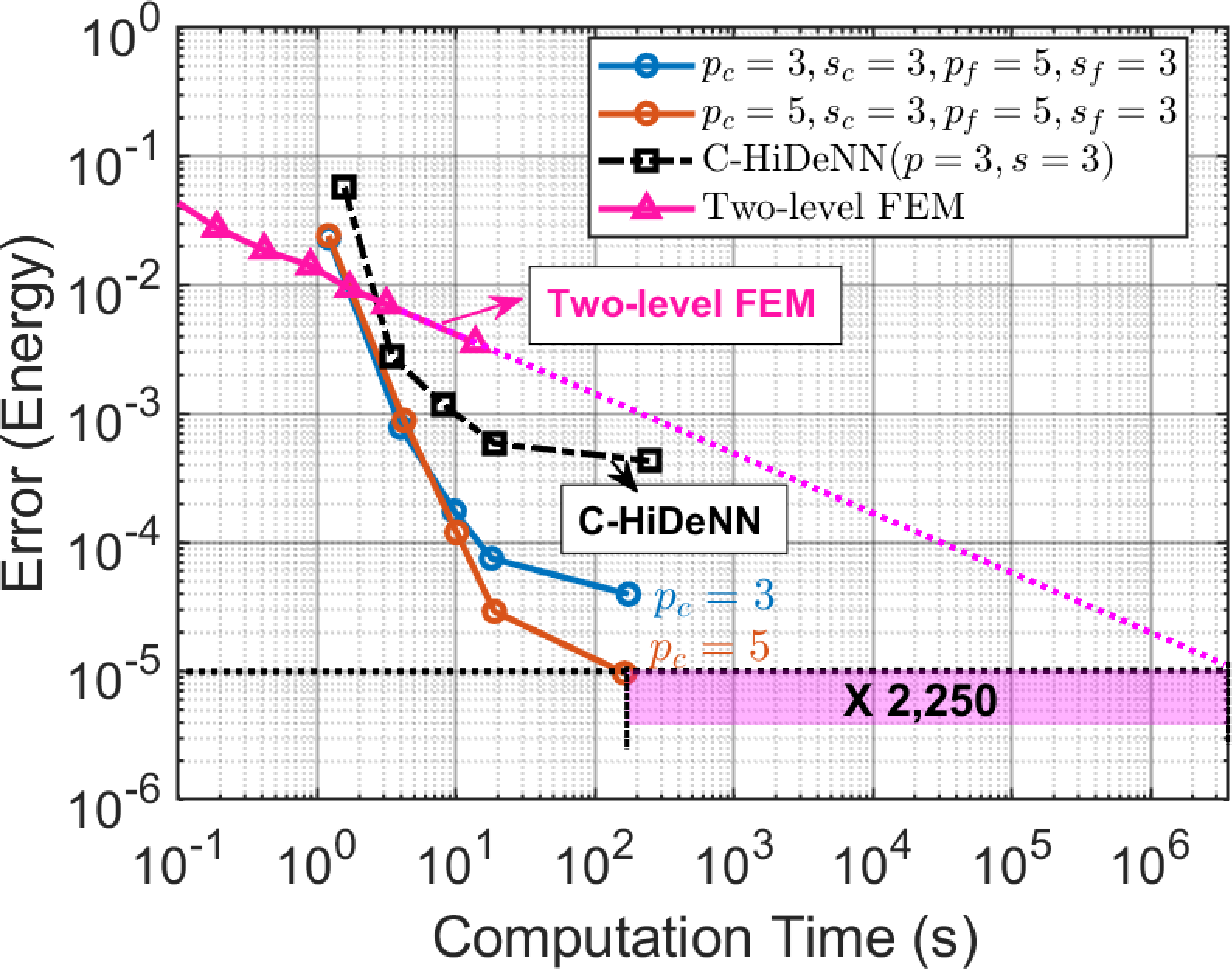}
\caption{Comparisons among two-level C-HiDeNN with optimal element size ratios, two-level FEM and single-level C-HiDeNN. }
\label{fig:Comparison_vms}
\end{figure}

We compare the performance of two-level C-HiDeNN with single-level C-HiDeNN and two-level VMS using FEM interpolations (short for two-level FEM) in Fig. \ref{fig:Comparison_vms} and Table \ref{table:timecomparison_methods}. When the optimal element size ratio is achieved, the two-level C-HiDeNN demonstrates orders of magnitude improved efficiency over the single-level C-HiDeNN and two-level FEM with the same level of accuracy. If we assume that the two-level FEM maintains the current scaling law between accuracy and computational time during mesh refinement, the two-level C-HiDeNN with $p_c=5$ will achieve a 2,250x speedup to reach an error below $10^{-5}$. The comparison among two-level C-HiDeNN curves with different orders indicates that a greater order of interpolations enhance the performance of two-level C-HiDeNN. 

Based on our theoretical and numerical findings, we draw the following observations:

$\bullet$ The method with greater two-level orders $p_c$ and $p_f$ is more efficient; specifically, it requires less computational time while achieving comparable accuracy. 

$\bullet$ The fine-level convolution patch size $s_f$ influence the performance slightly. (See appendix for more details)

$\bullet$ The optimal element size ratio decreases as $p_f$ increases. The optimal element size ratio $n_{opt}$ is the best one that leverages accuracy and efficiency. In general, a small ratio $n=2$ or $4$ is recommended for larger $p_f$.

$\bullet$ When taking the appropriate element size ratio, the two-level C-HiDeNN yields better efficiency than the single-level C-HiDeNN and two-level FEM with comparable accuracy. This demonstrates the efficacy of ML-VMS framework with C-HiDeNN.

\begin{table}[!htb]
\caption{Computation time comparison for different methods ($p=3, s=3$ or $p_c=3, s_c=3$ for C-HiDeNN) when reaching an energy-norm error less than $10^{-3}$. We take the optimal element size ratio for two-level C-HiDeNN ($n=6$ for $p_c=3, s_c=3, p_f=1, s_f=3$; $n=4$ for $p_c=3, s_c=3, p_f=3, s_f=3$;$n=2$ for $p_c=3, s_c=3, p_f=5, s_f=3$).}
\centering
\begin{tabular}{| c | c | c | c | c |}
\hline
Two-level FEM & C-HiDeNN & \makecell[c]{Two-level C-HiDeNN \\ ($p_f=1, s_f=3$)} & \makecell[c]{Two-level C-HiDeNN \\ ($p_f=3, s_f=3$)} & \makecell[c]{Two-level C-HiDeNN \\ ($p_f=5, s_f=3$)} \\ \hline
108.76 & 7.98 & 11.19 & 5.19 & \textbf{3.66} \\ \hline
\end{tabular}
\label{table:timecomparison_methods}
\end{table}

% \begin{figure}[htbp]
% \centering
% \includegraphics[width=0.9\textwidth]{ratio.png}
% \caption{Optimal ratio based on accuracy and efficiency (a) case 1: $p_c=3, s_c=2; p_f=1, s_f=1$; (b) case 2: $p_c=3, s_c=2; p_f=3, s_f=2$}
% \label{fig:ratio}
% \end{figure}

% \begin{figure}[htbp]
% \centering
% \includegraphics[width=1.0\textwidth]{ratio2.png}
% \caption{Accuracy and efficiency for C-HiDeNN with hyperparameters (a) $n=2$ (b) $n=4$ (c) $n=8$}
% \label{fig:ratio2}
% \end{figure}

% \caption{Error vs. computational time for different refinements and parameters.}
% \label{fig:Convergence_Eg1_ErrVsTime}
% \end{figure}

% \subsection{Two-level VMS with C-HiDeNN-TD for elliptic problem}
% % \subsubsection{Two-level VMS with C-HiDeNN-TD}

% We apply reduced-order methods to solve the following multi-modes solution problem: the Poisson's equation (\ref{eq:Eg1_Eqn}) with a seven-modes exact solution
% \begin{equation}
%     u(x,y)=\sum_{q=1}^7 \exp \left(-\dfrac{\pi}{2}(x-3.2-0.2k)^2 - \dfrac{\pi}{2}(y-3.2-0.2k)^2\right).
% \end{equation}
% We first use (single-level) C-HiDeNN-TD to verify Theorem 3, and study the effect of mode number. Next we study convergence rates for two-level VMS with C-HiDeNN-TD.

\subsubsection{Two-level VMS with C-HiDeNN-TD}

We then employ two-level VMS with C-HiDeNN-TD (short for two-level C-HiDeNN-TD), and take $Q_c=8$ modes for coarse level and $Q_f=14$ modes for fine level. Here, we assume that the truncated error less than $10^{-6}$ ($Q \geq 6$) is accurate enough. Detailed discussions on the selection of mode number refer to \ref{appenx:two-levelVMS_elliptic}.

\begin{figure}[htbp]
\centering
\includegraphics[width=0.6\textwidth]{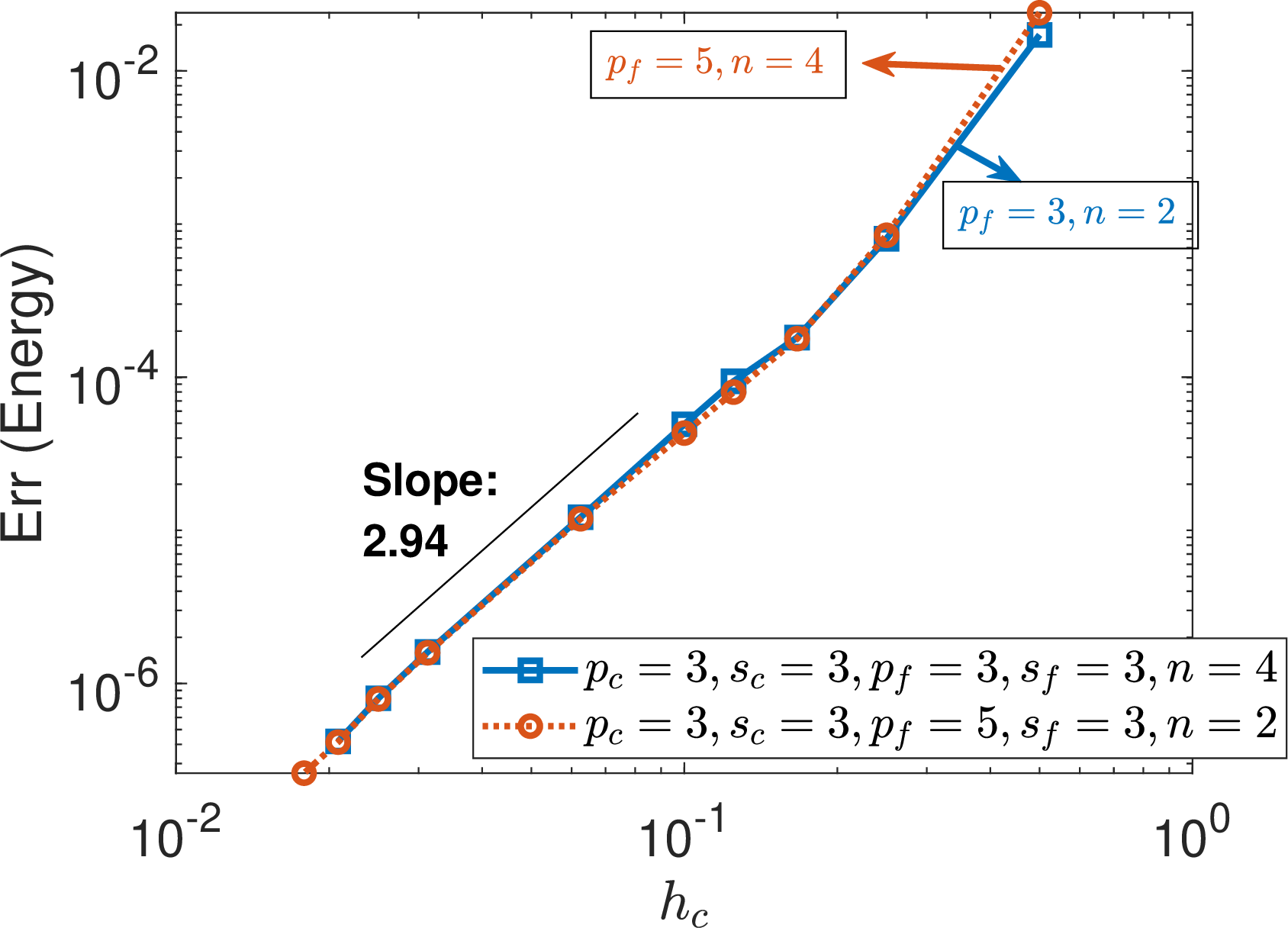}
\caption{Convergence curves for two-level C-HiDeNN-TD with optimal element size ratio: $n=4$ for $p_c=3, s_c=3, p_f=3, s_f=3$ and $n=2$ for $p_c=3, s_c=3, p_f=3, s_f=3$. }
\label{fig:Comparison_TD_0}
\end{figure}

We plot the convergence curve with the optimal element size ratio $n_{opt}$, as shown in Fig. \ref{fig:Comparison_TD_0}. In both cases of $p_f=p_c=3$ and $p_f=5>p_c=3$, the convergence rate is approximately $p_c=3$, and their convergence curves are close to each other. This is because the fine-level error term $C^{(f)}\cdot h_f^{p_f}$ reaches the same level of $C^{(c)}\cdot h_c^{p_c}$, when the optimal element size ratio $n_{opt}$ is taken.

\begin{figure}[htbp]
\centering
\includegraphics[width=0.6\textwidth]{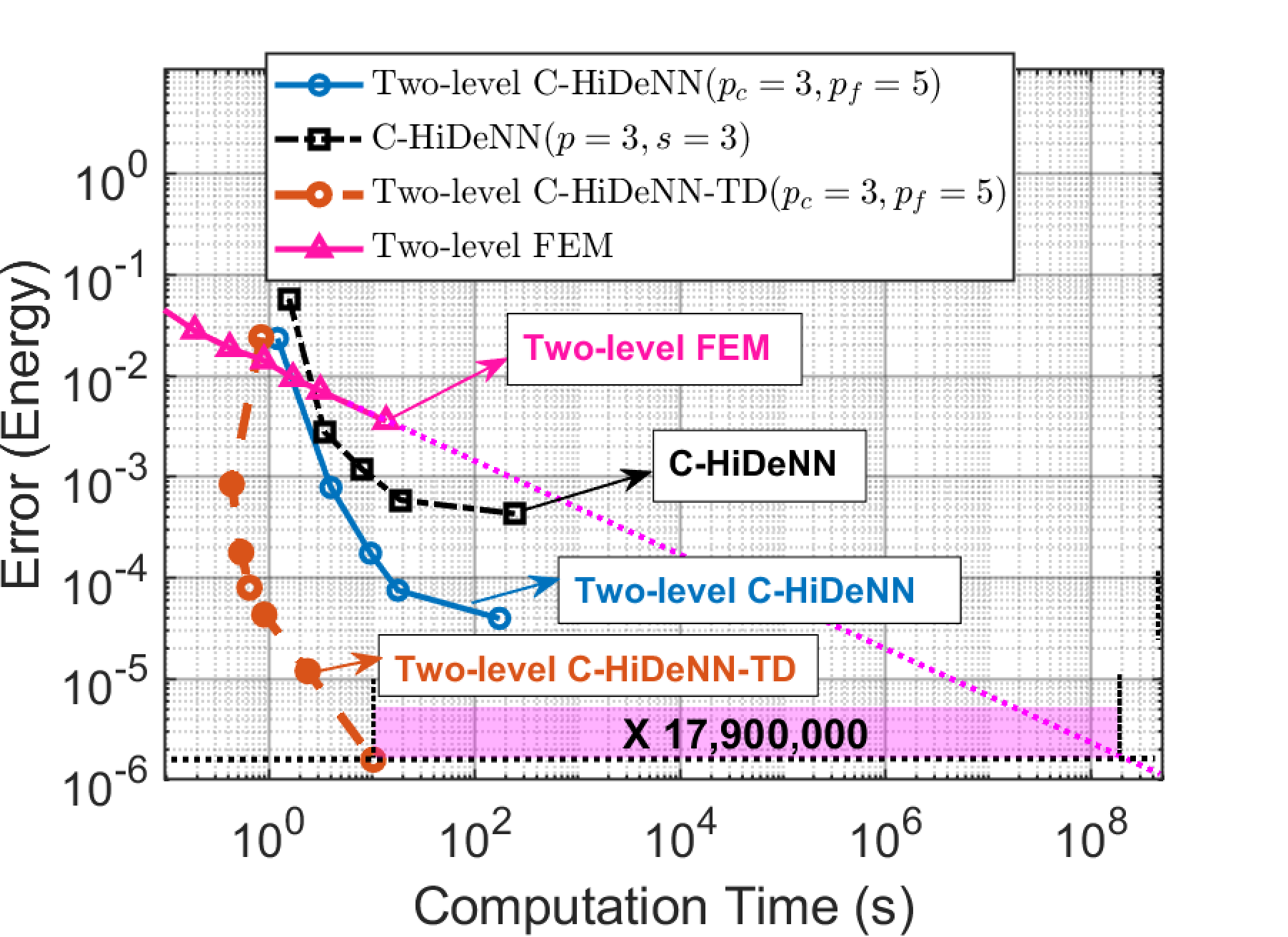}
\caption{Comparisons among two-level C-HiDeNN-TD with optimal element size ratios, two-level C-HiDeNN with optimal element size ratios, two-level FEM and (single-level) C-HiDeNN. }
\label{fig:Comparison_TD2}
\end{figure}

We compare the two-level C-HiDeNN-TD with other methods in Fig. \ref{fig:Comparison_TD2}. The two-level C-HiDeNN-TD is the most efficient method compared to single-level and two-level C-HiDeNN, two-level FEM. The two-level C-HiDeNN with the same parameters ($p_c=3, p_f=3, p_f=5, s_f=3, n=2$), which is the most efficient method in Fig. \ref{fig:Comparison_vms}, requires around 22.7 times more computational time than the two-level C-HiDeNN-TD when achieving an error of $10^{-4}$. Detailed computational time comparisons are listed in Table \ref{table:timecomparison_methods}. The two-level C-HiDeNN-TD takes less than one second to achieve the error less than $10^{-4}$. However, the two-level C-HiDeNN needs 14.54 seconds for the same level of accuracy. The two-level FEM requires more than 100 seconds for the error less than $10^{-3}$. For higher accuracy with an error of $10^{-6}$, the two-level C-HiDeNN-TD achieve $2.84\times 10^7$ times speedup 
compared to the two-level FEM.
% \begin{figure}[htbp]
% \centering
% \includegraphics[width=5in]{Convergence_RefineFine_C-HiDeNN-TD.png}

% \caption{Convergence curves for refining fine-level mesh.}
% \label{fig:Convergence_Eg1_Fine_C-HiDeNN-TD}
% \end{figure}

\begin{table}[!htb]
\caption{Computational time comparison for two-level FEM, C-HiDeNN ($p=3, s=3$), two-level C-HiDeNN ($p_c=3, s_c=3, p_f=5, s_f=3$) and two-level C-HiDeNN-TD ($p_c=3, s_c=3, p_f=5, s_f=3$) when reaching comparable energy-norm errors ($10^{-3}, 10^{-4}, 10^{-5}$ and $10^{-6}$). We take the optimal element size ratio for two-level VMS with C-HiDeNN and C-HiDeNN-TD ($n=2$ for $p_c=3, s_c=3, p_f=5, s_f=3$).}
\centering
\resizebox{1.1\linewidth}{!}{
\begin{tabular}{| c | c | c | c | c | c |}
\hline
Error & 2-level FEM & C-HiDeNN & \makecell[c]{2-level C-HiDeNN \\ ($p_f=5, s_f=3$)} &  \makecell[c]{2-level C-HiDeNN-TD \\ ($p_f=5, s_f=3$)} & \makecell[c]{Speedup \\ over 2-level FEM} \\ \hline
$10^{-3}$ & $108.76$ & $7.98$ & $3.66$ & $0.43$ & $2.53 \times 10^2$ \\ \hline
$10^{-4}$ & $2.73\times 10^4$* & - & $14.54$ & $0.64$ & $4.26\times 10^4$ \\ \hline
$10^{-5}$ & $3.66 \times 10^6$* & - & $ - $ & $3.15$ & $1.16\times 10^6$ \\ \hline
$10^{-6}$ & $4.92 \times 10^8$* & - & $ - $  & $17.34$* & $2.84\times 10^7$ \\ \hline

\end{tabular}}
\raggedright
* Computation time estimated by extrapolation
\label{table:timecomparison_methods}
\end{table}

% \subsection{Two-level VMS with C-HiDeNN-TD for elasticity problem}

\subsection{Space-time two-level VMS with C-HiDeNN and C-HiDeNN-TD for 1D heat problem}

Now we use VMS to solve a parabolic equation. Consider a heat equation within $\Omega \times [0,T]=[-1,1]\times[0,4]$, given by
\begin{equation} \label{eq:Eg3_Eqn}
    u_t(x,t) - u_{xx}(x,t) = f(x,t).
\end{equation}
We manufacture an exact solution:
\begin{equation}
    u^{Ext}(x,t)=\exp \left(-100 x^2\right)\left(1-\exp (-5t)\right),
\end{equation}
yielding boundary conditions
\begin{equation}
    u(-1,t)=\exp(-100)\left(1-\exp (-5t)\right), u(1,t)=\exp(-100)\left(1-\exp (-5t)\right),
\end{equation}
and zero initial conditions
\begin{equation}
    u(x,0)=0.
\end{equation}
The source term $f(x,y)$ is then derived from Eq. (\ref{eq:Eg3_Eqn}) using the exact solution. 

% \subsubsection{C-HiDeNN-TD for single-level space-time mesh}

% We first study the performance of C-HiDeNN-TD and determine the mode number. 

\begin{figure}[htbp]
\centering
\subfigure[Two-level mesh]{\includegraphics[width=2.5in]{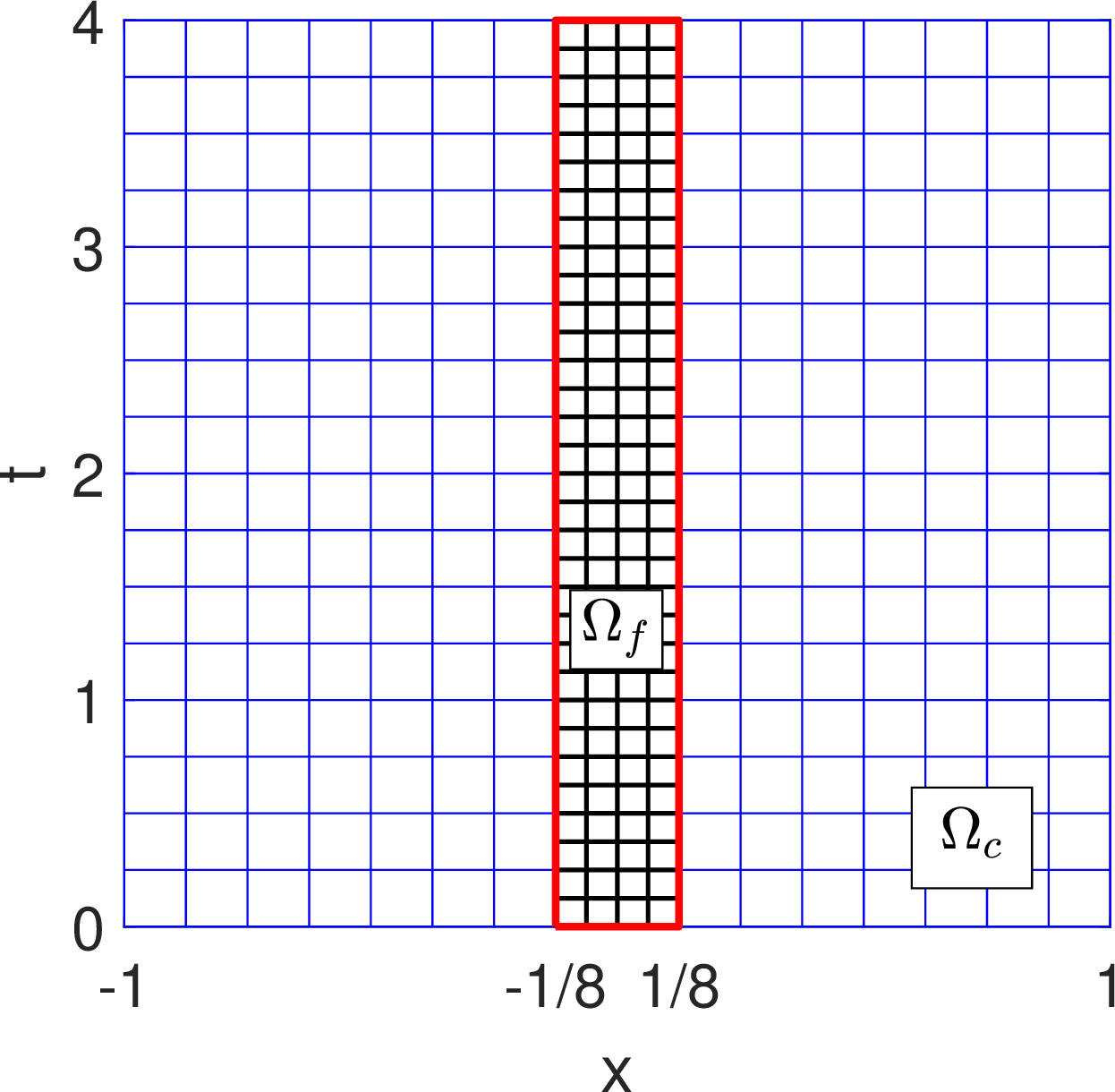}}
\subfigure[Exact solution]{\includegraphics[width=2.5in]{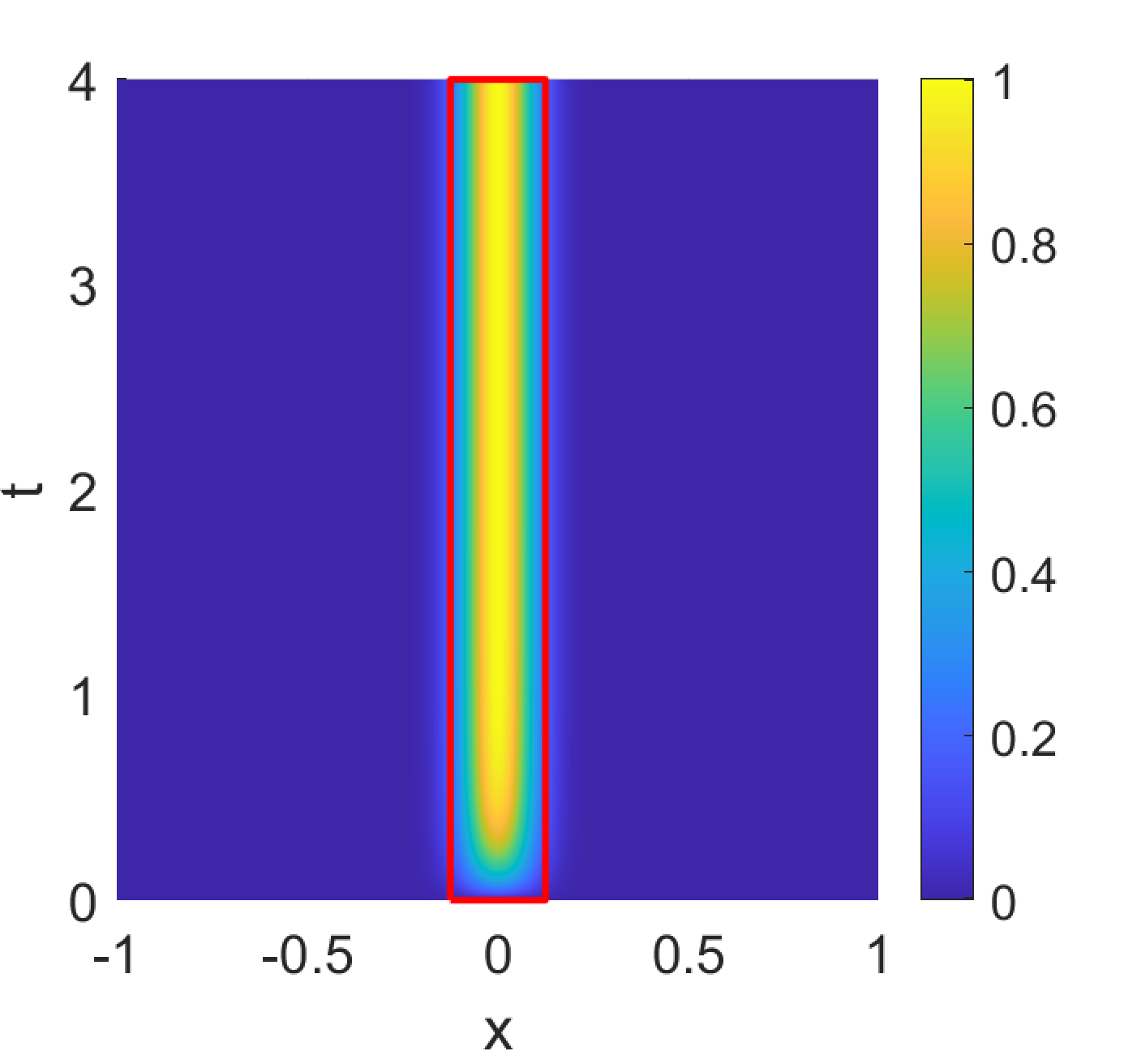}}
\caption{The schematic of the two-level space-time meshes. }
\label{fig:Example3}
\end{figure}

Figure \ref{fig:Example3} illustrates the space-time meshes for two-level VMS. The fine-level mesh is confined into a space-time subdomain $\Omega_f \times [0,T]=[-1/8,1/8]\times[0,4]$. $h_f$ and $\Delta t_f$ denote the fine-level element size and time step size, respectively. $h_c$ and $\Delta t_c$ denote the coarse-level element size and time step size, respectively. 

\begin{figure}[htbp]
\centering
\includegraphics[width=2.5in]{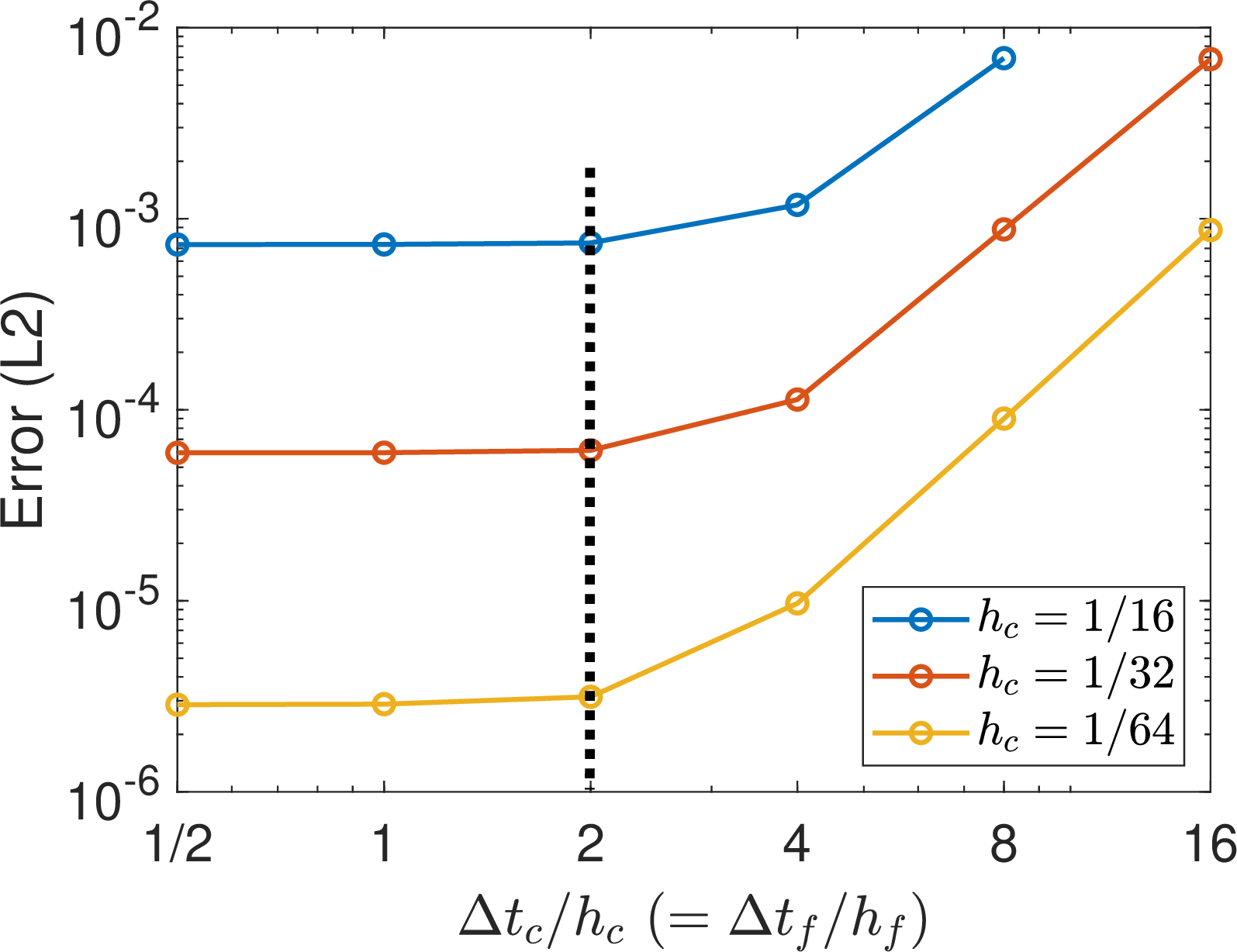}
\caption{Accuracy of two-level C-HiDeNN-TD for different ratios between the time step size and the element size. The parameters are taken as $p_c=3, s_c=3, p_f=3, s_f=3, n=2, Q_c=2, Q_f=4$. When $\Delta t_c / h_c = \Delta_f / h_f \leq 2$, the error remains almost unchanged. Therefore, we take the ratio $\Delta t_c / h_c = \Delta_f / h_f = 2$ in this example.}
\label{fig:STRatio}
\end{figure}

We refine the element size and time step size in each level simultaneously, and take $n=h_c/h_f=\Delta t_c / \Delta t_f$. Table \ref{table:largeTimeStep} lists errors of two-level C-HiDeNN-TD with $(p_c=3, s_c=3, p_f=3, s_f=3, n=h_c/h_f=\Delta t_c/\Delta t_f =2)$. As can be seen, the time step sizes $\Delta t_c$ used for all cases are much larger than the corresponding critical time steps $\Delta t_{Crit}=h_c^2/2$ in the forward Euler scheme. This is due to the fact that the proposed space-time formulation is unconditionally stable. Figure \ref{fig:STRatio} presents the variation in accuracy of the two-level C-HiDeNN-TD with parameters $p_c=3, s_c=3, p_f=3, s_f=3, n=2, Q_c=2, Q_f=4$ with respect to the ratio $\Delta t_c / h_c = \Delta_f / h_f$. When $\Delta t_c / h_c = \Delta_f / h_f \geq 4$, the error increases rapidly as this ratio increases. When $\Delta t_c / h_c = \Delta_f / h_f \leq 2$, it is observed that the error remains almost unchanged, because the error is dominated by the spatial discretization error. Therefore, we take $\Delta t_c / h_c = \Delta_f / h_f = 2$ in this example.

\begin{table}[!htb]
\caption{Performance of two-level C-HiDeNN-TD with $(p_c=3, s_c=3, p_f=3, s_f=3, n=h_c/h_f=\Delta t_c/\Delta t_f =2)$ for large time step sizes. The time step size $\Delta t_c$ is much larger than the critical time step $\Delta t_{Crit}=h_c^2/2$ for the forward Euler scheme, indicating that the space-time two-level C-HiDeNN-TD allows a large time step size. }
\centering
\begin{tabular}{| c | c | c | c | c | c | c |}
\hline
$h_c$ & $\Delta t_c$ & Err ($L^2$) & $\Delta_{Crit}$ for the forward Euler scheme \\ \hline
$1/16$ & $1/2$ & $6.92\times 10^{-3}$ & $1/512$ \\ \hline
$1/32$ & $1/4$ & $8.80\times 10^{-4}$ & $1/2048$ \\ \hline
$1/64$ & $1/8$ & $8.99\times 10^{-5}$ & $1/8192$ \\ \hline

\end{tabular}
\label{table:largeTimeStep}
\end{table}

With the same parameters ($p_c, s_c, p_f, s_f$ and $n$), we compare the efficiency of two-level C-HiDeNN-TD and two-level C-HiDeNN. When error reaches the level of $10^{-4}$, the two-level C-HiDeNN takes around 88 times more computational time than the two-level C-HiDeNN-TD. The two-level C-HiDeNN-TD with greater $p_f$ ($p_f=3$ and $p_f=5$) is much more efficient than that with small $p_f$ ($p_f=1$) for error less than $10^{-4}$. Therefore, we can obtain the similar conclusions with those for the Poission's problem: greater $p_f$ is recommended.

\begin{figure}[htbp]
\centering
\includegraphics[width=0.6\textwidth]{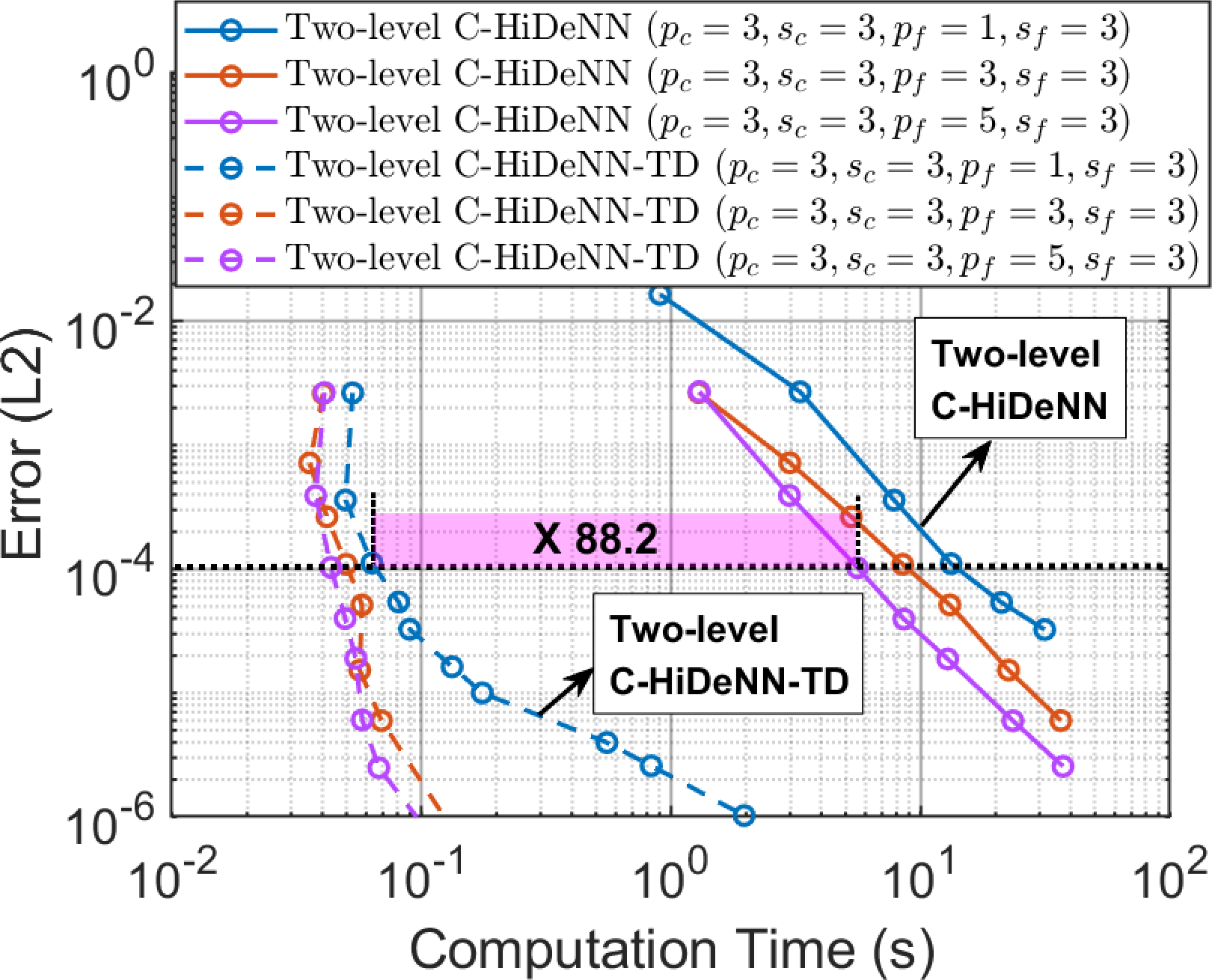}
\caption{Comparisons between two-level C-HiDeNN and two-level C-HiDeNN-TD. }
\label{fig:ST_TD_Comparison}
\end{figure}

\subsection{3D transient heat transfer simulation with moving heat source}

\begin{figure}[htbp]
\centering
\includegraphics[width=0.4\textwidth]{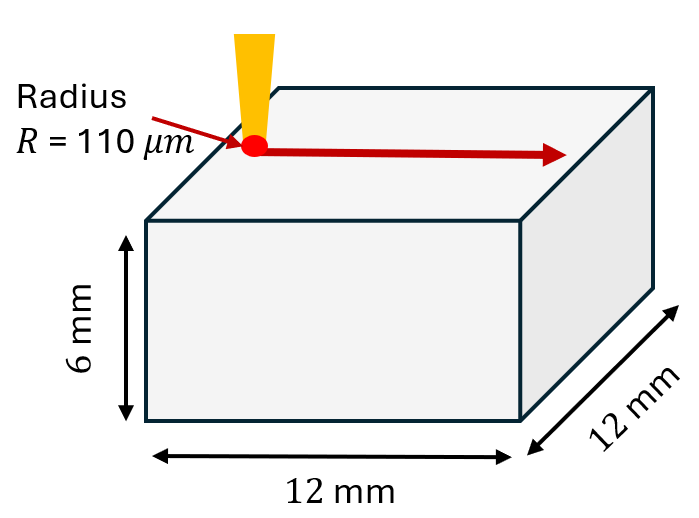}
\caption{LPBF simulation: single track within a 12 $\mathrm{mm^3}$ powder bed domain.}
\label{fig:AM_fig}
\end{figure}

In this section, we focus on modeling the laser powder bed fusion (LPBF) process in additive manufacturing (AM). Here, we demonstrate the capability of ML-VMS C-HiDeNN-TD to simulate single-track LPBF production as illustrated in Fig. \ref{fig:AM_fig}. 

The simulation process is governed by the following 3D heat equation:
\begin{equation} \label{eq:AM}
    \rho c_p u_t - \nabla \cdot (k \nabla u) = f(\bm{x}, t).
\end{equation}
where Ti-6Al-4V is adopted as the printing material. The detailed material and process parameters are defined in Table \ref{table:AM_parameter}. The laser source term $f(\bm{x}, t)$ is described by an ellipsoidal Gaussian intensity profile:
\begin{equation} \label{eq:AM_source}
    f(\bm{x}, t) = A \exp\left( -3\left[\dfrac{(x-x_c(t))^2}{R^2}+\dfrac{(y-y_c)^2}{R^2}+\dfrac{(z-z_c)^2}{D^2} \right] \right),
\end{equation}
where $(x_c(t), y_c, z_c)$ is the center of the laser source moving along the $x$-direction at constant velocity $v$. The laser intensity $A$ is defined by $A=\dfrac{6\sqrt{3}P\eta}{\pi^{3/2}R^2 D}$. 

\begin{table}[!htb]
\caption{Simulation parameters used in LPBF production case. Material parameters are taken from Ti-6Al-
4V literature. }
\centering
\begin{tabular}{| c | c | c | c |}
\hline
Parameter & Variable & Value &  Units \\ \hline
Thermal conductivity & $k$ & 22.0 &  $\mathrm{W \ m^{-1} \ K^{-1}}$ \\ \hline
Density & $\rho$ & 4.27 &  $\mathrm{g \ cm^{3}}$ \\ \hline
Specific heat capacity & $c_p$ & 745 &  $\mathrm{J \ kg^{-1} \ K^{-1}}$ \\ \hline
Laser scan speed & $v$ & 500 &  $\mathrm{mm \ s^{-1}}$ \\ \hline
Laser spot size radius & $R$ & 110 &  $\mathrm{\mu m}$ \\ \hline
Laser penetration depth & $D$ & 50 &  $\mathrm{\mu m}$ \\ \hline
Laser power & $P$ & 200 &  $\mathrm{W}$ \\ \hline
Absorptivity & $\eta$ & 0.25 &  - \\ \hline

\end{tabular}
\label{table:AM_parameter}
\end{table}

% \begin{figure}[htbp]
% \centering
% \includegraphics[width=3in]{MovingHeatSource.png}
% \caption{The schematic of the space-time mesh for a moving heat source. }
% \label{fig:Example1}
% \end{figure}

Before applying the ML-VMS C-HiDeNN-TD to the single-track simulation, we first apply it to a moving heat source problem with an analytical solution to investigate the speed, accuracy, RAM efficiency and disk storage requirement of the method.

We take half of a 12 $\mathrm{mm^3}$ powder bed domain as the computational domain: $\Omega=[-6,6]\times[-6,6]\times[-6,0] \mathrm{mm^3}$. The analytical solution is
\begin{equation}
    u^{\text{Ext}}(\bm{x}, t) = \exp\left[ -3\left(\dfrac{(x-x_c(t))^2}{R^2}+\dfrac{(y-y_c)^2}{R^2}+\dfrac{(z-z_c)^2}{D^2} \right) \right]\cdot \left(1-\exp(5t)\right).
\end{equation}
Its center $(x_c, y_c, z_c)=(-5-vt, 0, 0)$ is located at the top surface $z=0$ and moves along the $x$-direction at a constant velocity $v=500 \mathrm{mm \ s^{-1}}$. The simulation duration is $20 \mathrm{ms}$. The Neumann boundary condition applies to the top surface, denoted by $\Gamma_N$, is
\begin{equation} \label{eq:AM_Neumann}
    u_z|_{\Gamma_N}=0.
\end{equation}
The Dirichlet boundary condition for the remaining surfaces, denoted by $\Gamma_D$, is
\begin{equation}
    u|_{\Gamma_D}=u^{\text{Ext}}|_{\Gamma_D}.
\end{equation}
Here, the value of $u|_{\Gamma_D}$ tends to zero with double-precision accuracy.
The initial condition is 
\begin{equation}
    u(\bm{x}, 0)=0.
\end{equation}

The source term $f(\bm{x},t)$ corresponding to this analytical solution is obtained by methods of manufactured solution. Note that the analytical solution is inherently non-separable due to the exponential term, necessitating a large mode number for C-HiDeNN-TD to achieve the target accuracy. However, since the trajectory of the moving heat source is known a priori, we can employ the following coordinate transformation to effectively reduce the mode number required in TD in the newly transformed coordinate system:
\begin{equation}
    x = \left\{
    \begin{array}{cc}
        \xi + x_c(t) \dfrac{6+\xi}{6-k_s} & \xi\in[-6,-k_s), \\
        \xi+x_c(t) & \xi\in[-k_s,k_s), \\
        \xi+x_c(t)\dfrac{6-\xi}{6-k_s} & \xi\in[k_s,6).
    \end{array}
    \right.
\end{equation}
$k_s$ is the size of a localized domain that covers the core component of the laser source term $f(\bm{x}, t)$. For convenience, in the framework of multi-level VMS, we also take it as the domain size of the finer level.
This transformation maps the physical domain $\Omega(x,y,z,t)$ to the reference domain $\tilde{\Omega}(\xi,y,z,t)$, where the tilde ($^\sim$) denotes quantities in the reference domain. In $\tilde{\Omega}(\xi,y,z,t)$, the transformed moving heat source is fixed at the point $(\xi,y,z)=(0,0,0)$. A general coordinate transformation for an arbitrary velocity vector $\bm{v}$ is provided in \ref{appenx:CoordinateTransformation}. 

Here, we employ a three-level C-HiDeNN-TD to solve this problem. With coordinate transformation, we create the three-level mesh in the reference domain rather than in the physical domain. The level 1 domain $\tilde{\Omega}_1$ covers the whole domain in the reference coordinate. The level 2 domain is $\tilde{\Omega}_2=[-k_s, k_s]\times[-k_s, k_s]\times[-a, 0]$, which is associated with $\Omega_2=[x_c(t)-k_s, x_c(t)+k_s]\times[-k_s, k_s]\times[-k_s, 0]$ in the physical space. The level 3 domain is $\tilde{\Omega}_3=[-k_s/4, k_s/4]\times[-k_s/4, k_s/4]\times[-k_s/4, 0]$, which is associated with $\Omega_3=[x_c(t)-k_s/4, x_c(t)+k_s/4]\times[-k_s/4, k_s/4]\times[-k_s/4, 0]$ in the physical space. The subdomains $\Omega_2$ and $\Omega_3$ for level 2 and level 3 track the moving sources, with sizes of $2k_s$ and $k_s/2$, respectively. Here, $k_s=0.8$ is taken to ensure that the core part of the source is adequately covered.

The approximation function $u_l$ $(l=1,2,3)$ for each level is also defined in the reference domain:
\begin{equation}
    u_l(\xi ,y,z,t)=\sum_{q=1}^{Q_l} \bm{\tilde{N}}^{1D, (l)}_\xi(\xi) \bm{u}^{(q),(l)}_\xi \cdot \bm{\tilde{N}}^{1D, (l)}_y(y) \bm{u}^{(q),(l)}_y \cdot \bm{\tilde{N}}^{1D, (l)}_z(z) \bm{u}^{(q),(l)}_z \cdot \bm{\tilde{N}}^{1D, (l)}_t(t) \bm{u}^{(q),(l)}_t, l=1,2,3
\end{equation}
with the element size $h_l$ and time step size $\Delta t_l$. Similar to the case of the 1D heat equation, the space-time C-HiDeNN-TD allows for large time steps, with a ratio of $\alpha = h_l/\Delta t_l= 0.05\mathrm{mm \ ms^{-1}}$. Through numerical tests, we observed that using smaller time step sizes slightly improved accuracy. For clarity, we define the size of space-time elements as follows: $\tilde{h}_l=\sqrt{(3 h_l^2+\alpha (\Delta t_l)^2)}, l=1,2,3$. Here, we assume that the scaled time dimension is regarded as the fourth spatial dimension. $\tilde{h}_l$ is the diagonal size of scaled spatio-temporal elements. As refining the time step and spatial element simultaneously, we define the ratio $\tilde{n}_{l,k}=\tilde{h}_l/\tilde{h}_k = h_l/h_k = \Delta t_l/ \Delta t_k$.

\begin{figure}[htbp]
\centering
\includegraphics[width=0.6\textwidth]{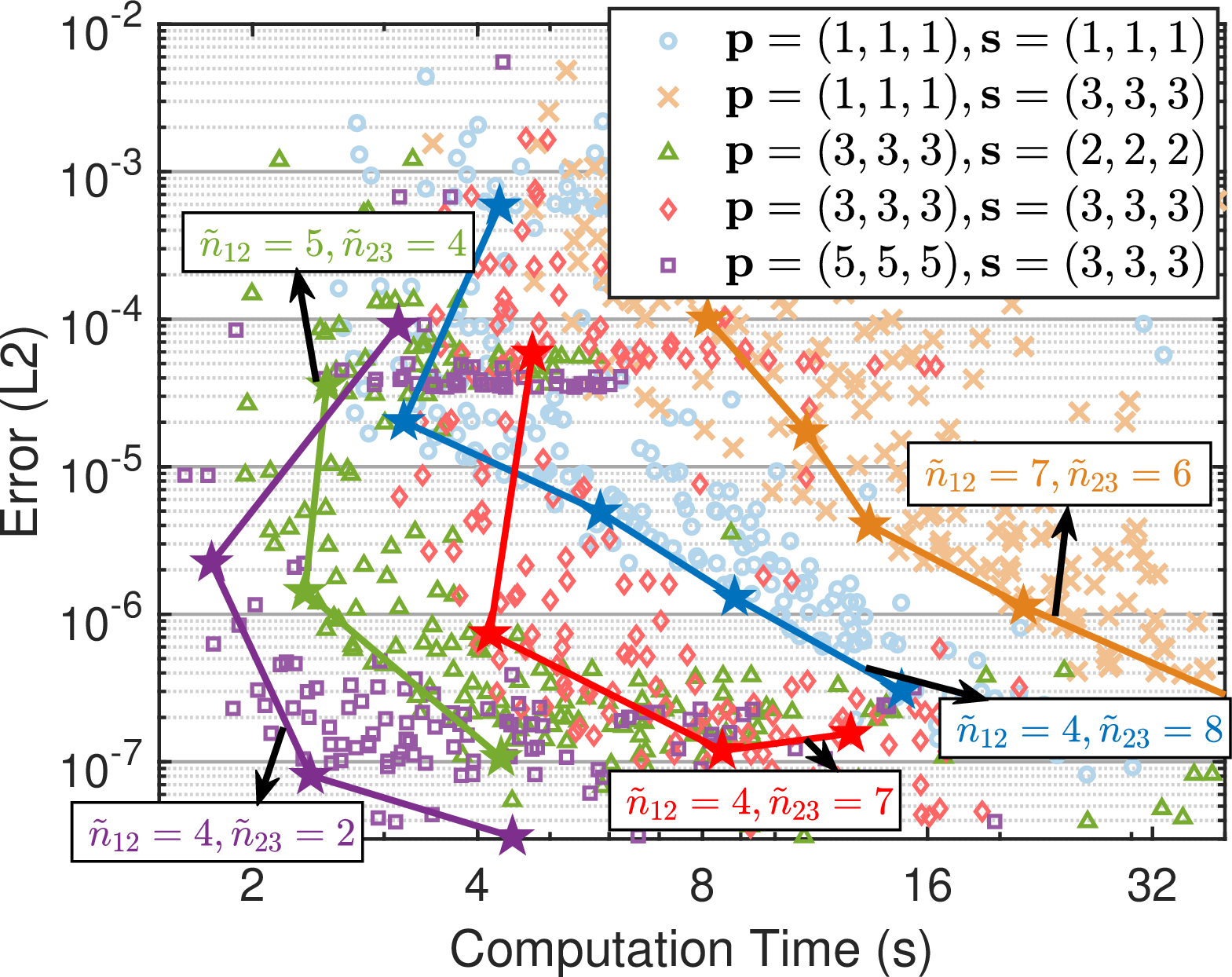}
\caption{Error vs. computational time for various refinements and parameters. Points of different colors represent results for different sets of parameters. The star-marked curves represent results with the optimal combination of element size ratios $\tilde{n}_{12}=\tilde{h}_1/\tilde{h}_2$ and $\tilde{n}_{23}=\tilde{h}_2/\tilde{h}_3$. }
\label{fig:Comparison_TD_3level}
\end{figure}

The numbers of modes for the three levels are $Q_1=2, Q_2=9$ and $Q_3=15$, respectively. The results for various parameters and refinements are illustrated in Fig. \ref{fig:Comparison_TD_3level}. Different colors represent results for distinct sets of hyperparameters $\bm{p}=(p_1, p_2, p_3)$ and $\bm{s}=(s_1, s_2, s_3)$ for different levels. The star-marked curves indicate those with the optimal combination of space-time element size ratios $\tilde{n}_{12}=\tilde{h}_1/\tilde{h}_2$ and $\tilde{n}_{23}=\tilde{h}_2/\tilde{h}_3$ for each set of hyperparameters. It is observed that higher-order approximations require less computational time to achieve comparable accuracy, demonstrating greater efficiency. Furthermore, when using the same orders $\bm{p}=(p_1,p_2,p_3)$, a smaller patch size $\bm{s}$ leads to more efficient results. This is because the computational overhead can decrease with smaller patch size since it directly influences the sparsity of the linear systems.

\begin{figure}[htbp]
\centering
\subfigure{\includegraphics[width=3in]{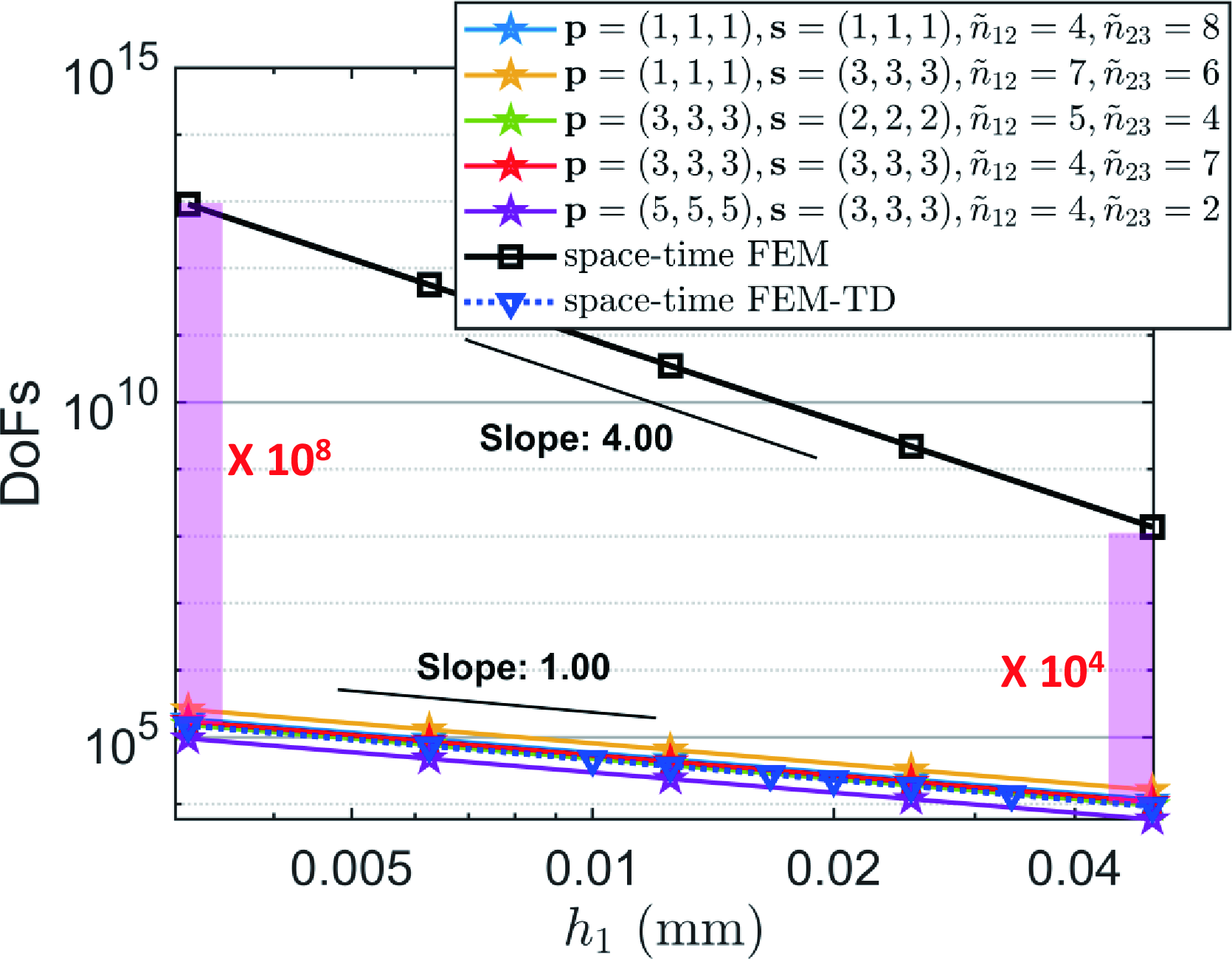}}
\subfigure{\includegraphics[width=3in]{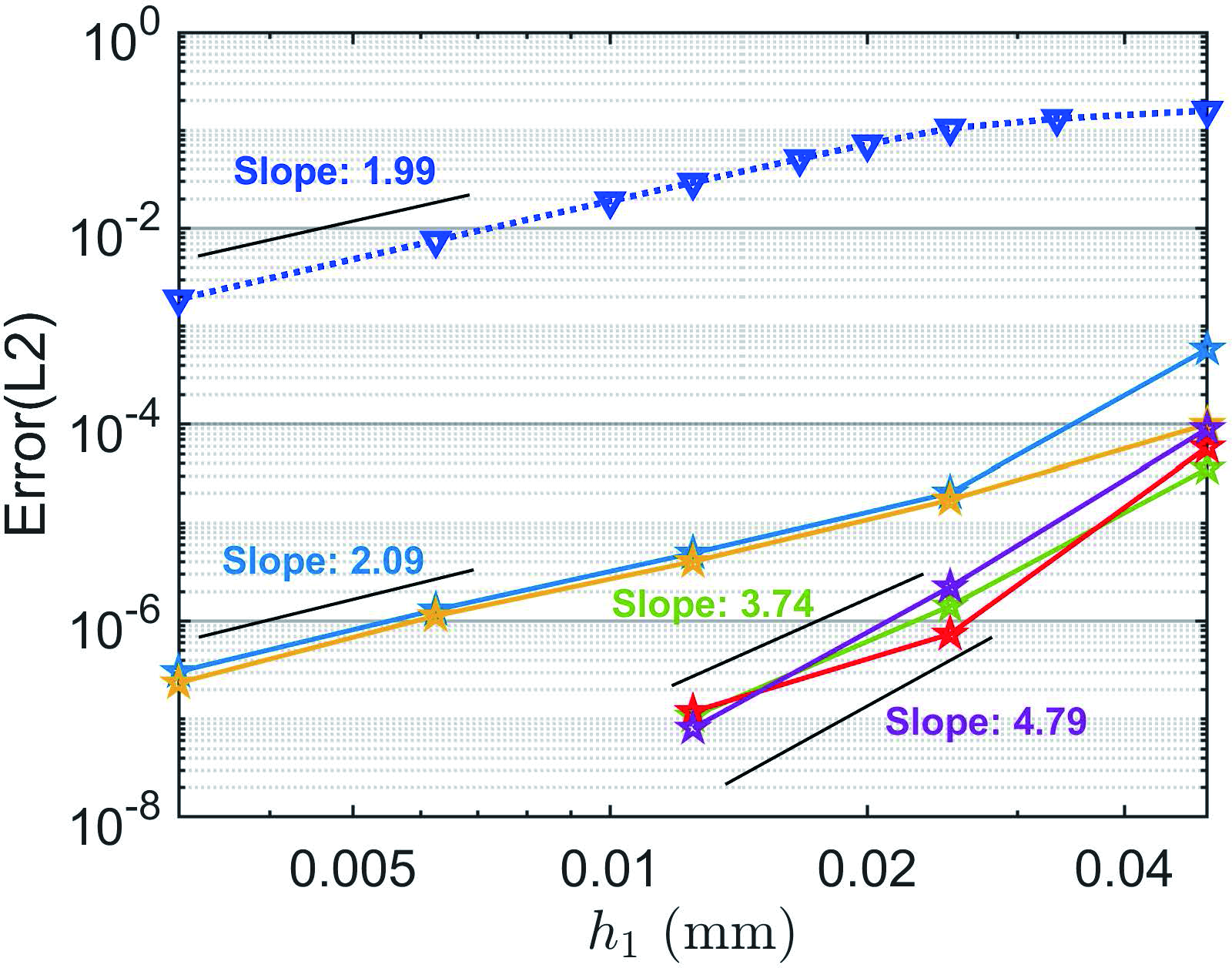}}
\subfigure{\includegraphics[width=3in]{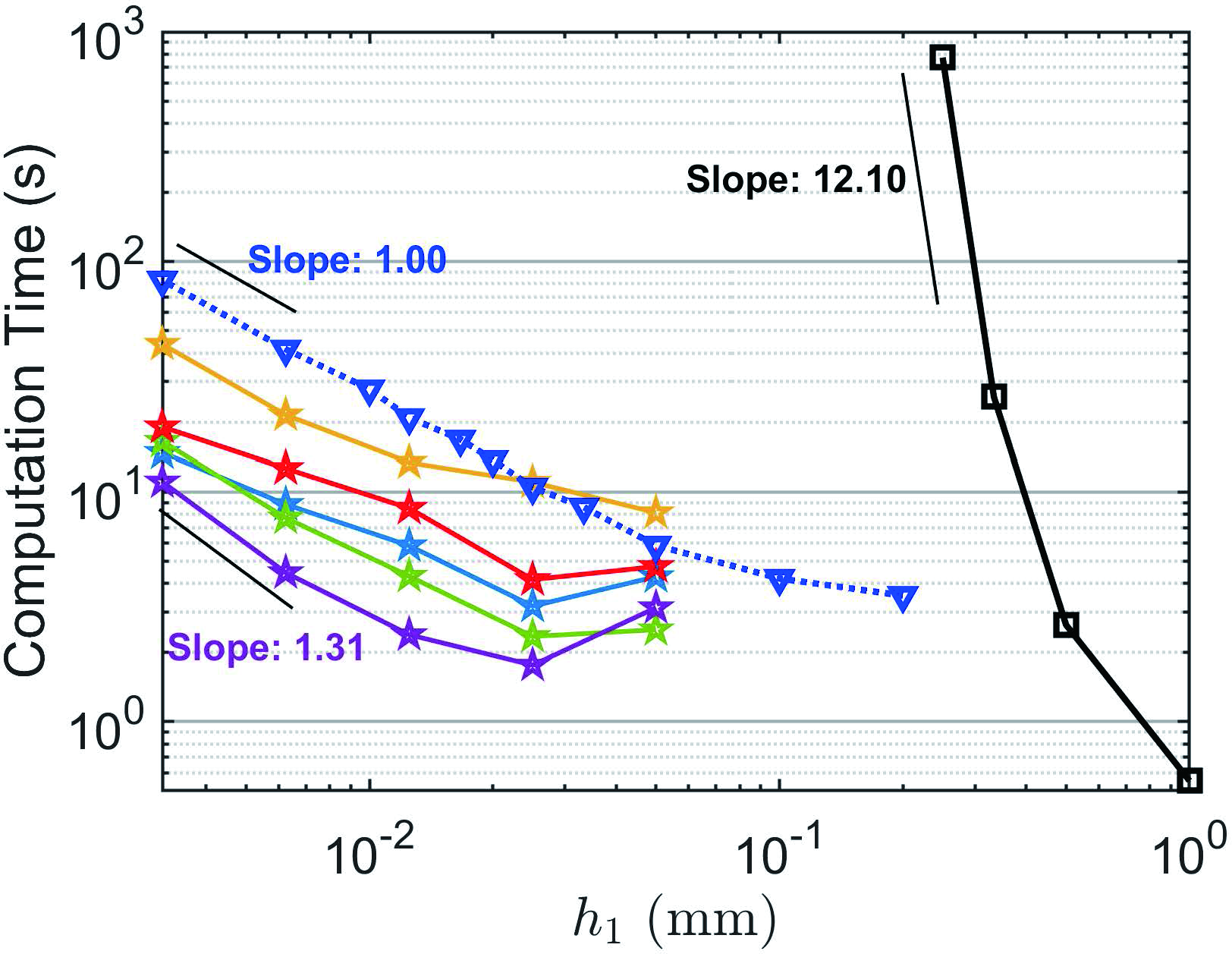}}
\subfigure{\includegraphics[width=3in]{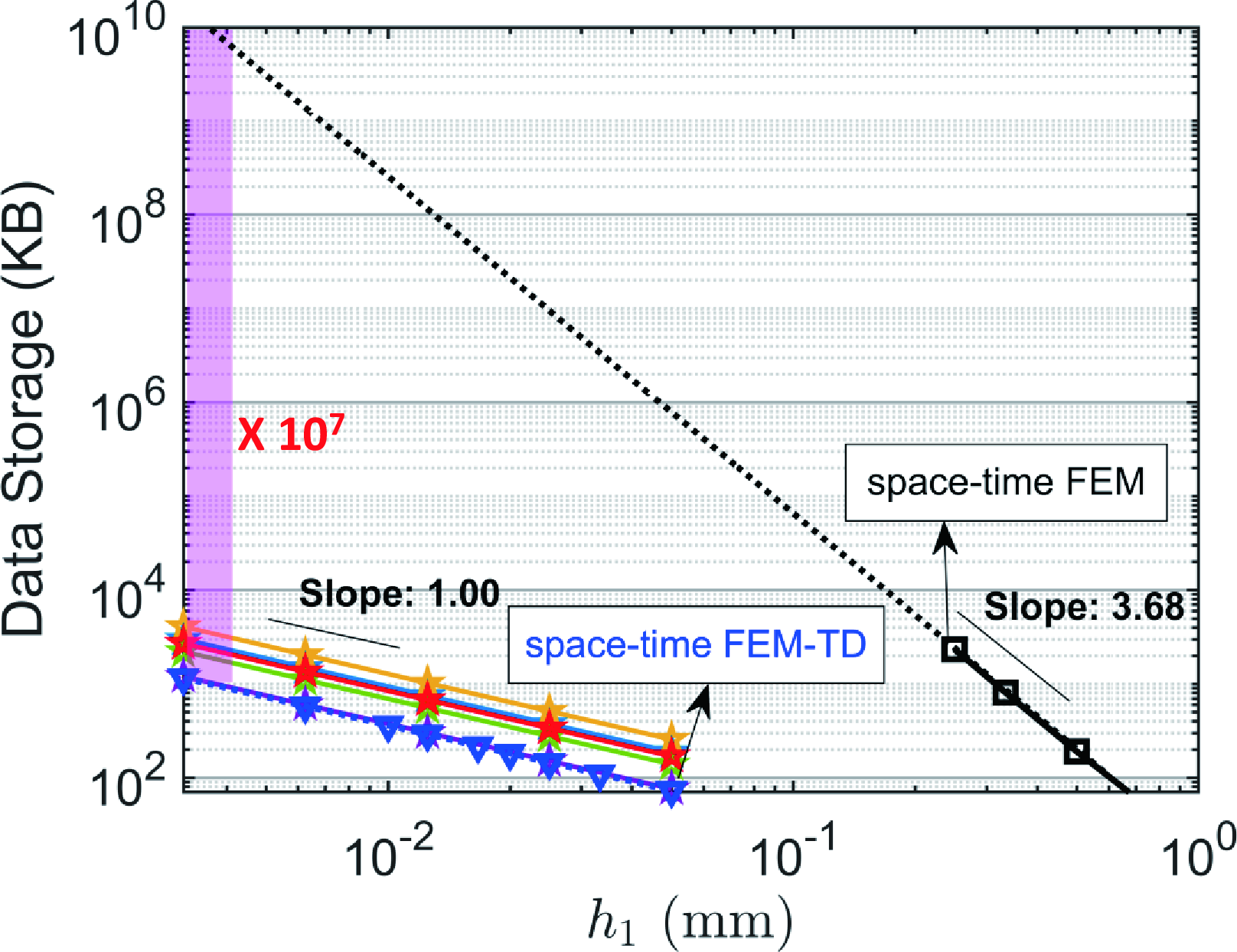}}

%prof liu prefers to use h_tilda instead of h1 because this is space-time
%space-time in multilevel should be in highlight

\caption{Convergence curves for DoFs, error, computational time and data storage with respect to coarsest element size. Each curve represents results for a set of parameters $\bm{p}=(p_1, p_2, p_3)$ and $\bm{s}={s_1, s_2, s_3}$ with optimal element size ratios. The space-time FEM with a uniform mesh and its reduced-order version, the space-time FEM-TD, are taken for comparison. The space-time FEM-TD incorporates tensor decomposition and coordinate transformation techniques; it shares the same space-time mesh as the three-level C-HiDeNN-TD at Level 1, and uses the same mode number, specifically, $Q=15$.}
\label{fig:3LevelCTD}
\end{figure}

To systematically study the performance of ML-VMS C-HiDeNN-TD, we focus on comparing error, computational time, DoFs and data storage with respect to the coarsest element size $h_1$ in Fig. \ref{fig:3LevelCTD}. We take space-time FEM with a uniform mesh and its reduced-order version, the space-time FEM-TD as the baselines for comparison. Here, the space-time FEM employs a continuous Galerkin scheme, using the same spatial element size and time step sizes. The space-time FEM-TD used here adopts the single-scale uniform FE mesh and linear FE basis function but incorporates tensor decomposition and coordinate transformation for model order reduction. The same mode number ($Q=15$) is used as in three-level C-HiDeNN-TD. For fair comparison, all methods adopt the sparse direct solver in solving linear algebraic systems. 

The DoFs of the three-level C-HiDeNN-TD scale linearly with respect to the element size,  whereas the DoFs of the space-time FEM exhibit quartic scaling. For the finest mesh investigated, the space-time FEM requires $10^{13}$ DoFs whereas the three-level C-HiDeNN-TD only requires $10^5$ DoFs. Owing to the technique of tensor decomposition, the space-time FEM-TD also exhibits linear scaling in DoFs. 

In terms of accuracy, the performance of FEM-TD is much lower than that of the three-level C-HiDeNN-TD at comparable resolutions. As expected, the $p$-th order C-HiDeNN-TD achieves the $(p+1)$-th order $L^2$-norm convergence rate. When $p=5$, the convergence rate for three-level-level C-HiDeNN-TD is approximately $4.79$. This suboptimal rate may be attributed to machine precision.

The computational cost of the three-level C-HiDeNN-TD increases with mesh refinement at rates ranging from $0.59$ to $1.31$. As a data-free reduced-order model, the computational cost of the space-time FEM-TD scales linearly with refinement but is still much higher than three-level C-HiDeNN-TD. It is observed that the computational time for $\bm{p}=(3,3,3), \bm{s}=(3,3,3)$ is consistently longer than that for $\bm{p}=(3,3,3), \bm{s}=(2,2,2)$. This is because a larger $\bm{s}$ results in a wider semi-bandwidth of the stiffness matrices, thereby increasing computational time. However, the comparison among different orders become complicated, as the order affects solution-related resolutions (see the upper-right subplot in Fig. 14), and these resolutions can significantly influence the number of iterations in the TD solution scheme. Although higher-order C-HiDeNN-TD methods have larger $\bm{s}$, their superior resolutions lead to fewer iterations, thus exhibiting significantly less computational time under the same refinement level. The cost of the full-order space-time FEM with a uniform mesh is significantly higher than that of all reduced-order methods, and it increases rapidly at a rate of $12.10$ for the finest space-time mesh studied.

Given that data storage is proportional to the number of DoFs, the three-level C-HiDeNN-TD also exhibits linear scaling in its data storage. Its data storage for the finest mesh studied is projected to be at least $10^7$ times lower than that of the space-time FEM. 

We can conclude that higher-order ML-VMS C-HiDeNN-TD leads to greater efficiency and smaller data storage while maintaining comparable accuracy. This conclusion is further supported by Table \ref{table:3level_time}. We list computational time and data storage for various parameters with optimal element size ratios when achieving an accuracy of the order of $10^{-7}$. It is observed that the computational time and data storage for $\bm{p}=(5,5,5)$ are the smallest. Given similar level of accuracy and taking the same orders across different levels, smaller $\bm{s}$ effectively reduces both time and storage costs. The optimal three-level C-HiDeNN-TD achieves $4,854$ times speedup and $533$ times less storage compared to the space-time FEM-TD. This is because the space-time FEM-TD requires much finer meshes to achieve comparable accuracy ($10^{-7}$). The three-level C-HiDeNN-TD is projected to achieve a speedup of over $10^{50}$ and a reduction in storage of $10^{15}$ compared to the full-order space-time FEM. This projection assumes that the full-order model employs linear-order interpolations, uniform space-time meshes, and a direct solver for the resulting linear algebraic systems, while targeting high-resolution solutions. However, it should be noted that in practice, the efficiency of the full-order method could also be improved through the use of high-order interpolations, adaptive meshing, or iterative solvers. 

\begin{table}[!htb]
\caption{Performance for different parameters $\bm{p}=(p_1, p_2, p_3)$ and $\bm{s}=(s_1, s_2, s_3)$ with comparable accuracy ($l^2$-norm error of the order of $10^{-7}$). The optimal element size ratios are taken. }
\centering
\resizebox{1.1\linewidth}{!}{
\begin{tabular}{| c | c | c | c | c | c | c | c |}
\hline
$\bm{p}=(p_1, p_2, p_3)$ & $\bm{s}=(s_1, s_2, s_3)$ & \makecell[c]{Optimal\\ $(\tilde{n}_{12}, \tilde{n}_{23})$} & \makecell[c]{Given\\ Accuracy} & \makecell[c]{Computation \\ Time (s)} & \makecell[c]{Speedup over \\ FEM-TD}  & \makecell[c]{Data \\Storage (MB)} & \makecell[c]{Storage gain \\ over FEM-TD}\\ \hline
\multirow{2}{*}{(1,1,1)} & (1,1,1) & (4,8) & $3.05\times10^{-7}$ & 14.77 & 785 & 2.94 & 54 \\
\cline{2-8}
& (3,3,3) & (7,6) & $2.34\times10^{-7}$ & 43.68 & 266 & 4.00 & 40\\ \hline
\multirow{2}{*}{(3,3,3)} & (2,2,2) & (5,4) & $1.07\times10^{-7}$ & 4.30 & 2698 & 0.54 & 296 \\
\cline{2-8}
& (3,3,3) & (4,7) & $1.19\times10^{-7}$ & 8.50 & 1365 & 1.32 & 121\\ \hline
(5,5,5) & (3,3,3) & (4,2) & $8.08\times10^{-8}$ & 2.39 & \textbf{4854} & 0.30 & \textbf{533} \\ \hline
\multicolumn{3}{|c|}{Space-time FEM-TD (single level)} & $1\times 10^{-7}$ & $1.16\times 10^4$* & - & $160$* & - \\
\hline
\multicolumn{3}{|c|}{Space-time FEM (single level)} & $1\times 10^{-7}$ & $6.20\times 10^{50}$* & - & $1.93 \times 10^{15}$* & - \\
\hline

\end{tabular}}
\\
\raggedright
* These data are estimated via extrapolation, based on the results of the space-time FEM (with a coarse uniform mesh) and space-time FEM-TD. Here, we assume the number of modes are sufficient to ensure that space-time FEM-TD achieves accuracy comparable to that of space-time FEM. We thus estimate the resolution required for space-time FEM to reach the target accuracy.
\label{table:3level_time}
\end{table}

\begin{figure}[htbp]
\centering
\includegraphics[width=0.5\textwidth]{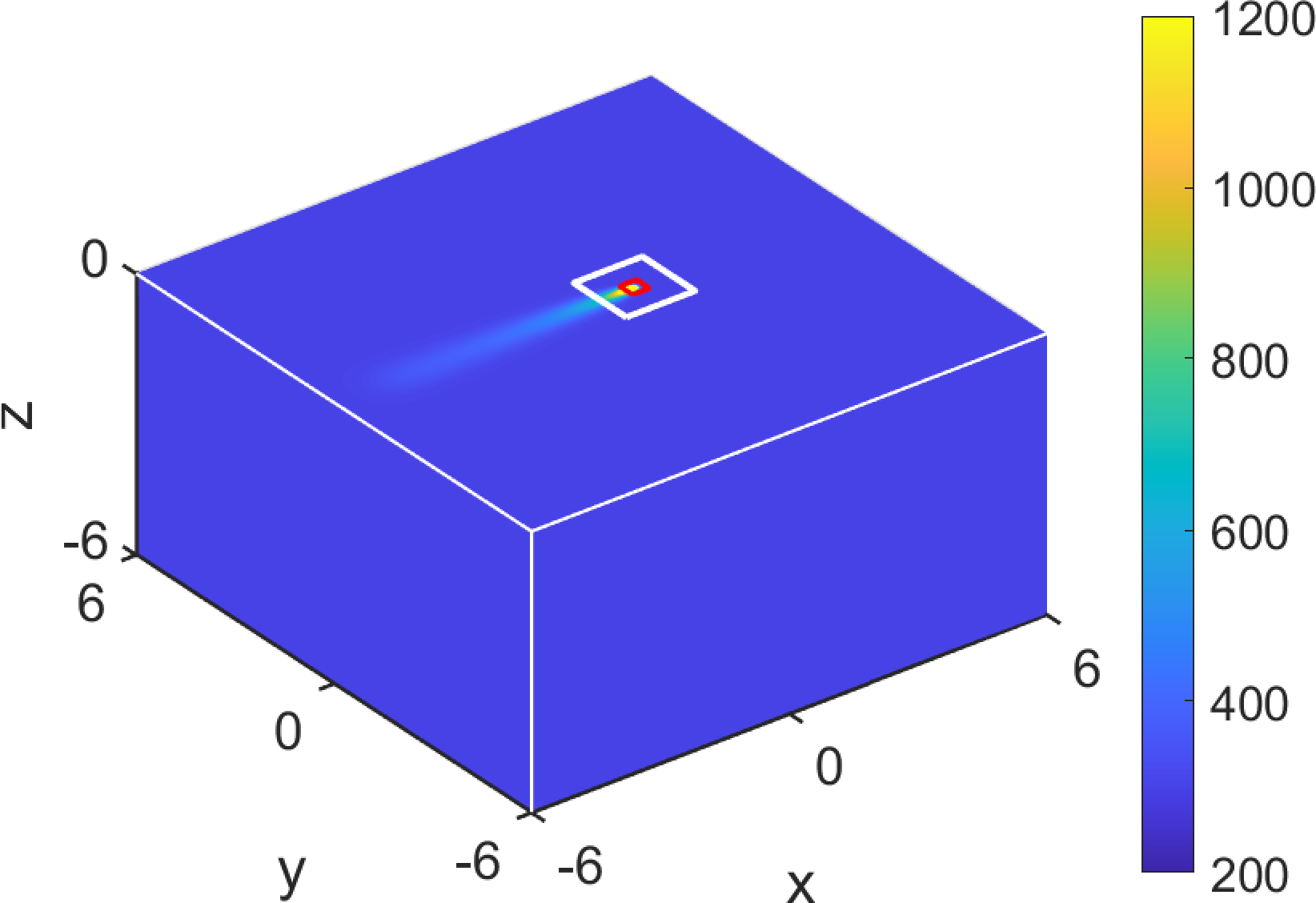}
\caption{Temperature solution for LPBF production at time $t=12 \mathrm{ms}$. The white box represents Level 2, and the red box represents Level 3. These two levels track the moving sources.}
\label{fig:AM_sol}
\end{figure}

Subsequently, we apply the three-level C-HiDeNN-TD to the simulation of single-track LPBF production, which is governed by the 3D heat equation (\ref{eq:AM}) and the laser source term given in Eq. (\ref{eq:AM_source}). The zero-flux Neumann boundary condition Eq. (\ref{eq:AM_Neumann}) is applied to the top surface, and the ambient temperature $T_{amb} = 298.15 \mathrm{K}$ is used for Dirichlet boundary conditions for the remaining surfaces. Based on the previous study for hyperparameters, we choose $\bm{p}=(5,5,5), \bm{s}=(3,3,3)$ and $\tilde{n}_{12}=4, \tilde{n}_{23}=2$ for optimal performance. 

In the simulation, Level 1 is discretized using a $480 \times 480 \times 240$ element mesh with a mesh size equal to $25 \mathrm{\mu m}$. Level 2 and Level 3 are discretized into $256 \times 256 \times 128$ and $128 \times 128 \times 64$ uniform meshes, respectively. The resulting smallest spatial element size is $3.125 \mathrm{\mu m}$, which corresponds to $2.83 \times 10^{10}$ spatial DoFs for a uniform mesh over the whole spatial domain if only single-level mesh is used.
The whole simulation takes about $306$ s on a single Interl Core i7-11700F CPU. The temperature solution at time $t=12$ ms is presented in Fig. \ref{fig:AM_sol}.

%\section{Outlook}

\section{Discussions}
\label{sec:discussions}
The proposed ML-VMS framework is distinguished from existing multilevel methods and GO-MELT by its novel solution construction in overlapping domains. Instead of decomposing the solution into additive coarse and fine level components, this new approach represents the solution directly via C-HiDeNN interpolation with different hyperparameters defined on the linear FE mesh. This allows for arbitrary order of interpolation across all levels without increasing the degrees of freedom (DoFs), which significantly enhances computational efficiency and reduces memory usage. In this paper, we generalize the framework to an $m$-level method and provide a rigorous error analysis that considers both mesh discretization errors at each level and errors arising from interface assumptions. Our numerical results demonstrate that higher-order approximations achieve superior computational efficiency and lower memory requirements while maintaining comparable accuracy.

C-HiDeNN-TD has been leveraged in the ML-VMS framework to further accelerate the computational time while minimizing the requirement of memory and storage requests for high-dimensional large-scale analysis. TD exhibits a linear growth of the number of DoFs instead of the polynomial growth in standard FEM. Contrary to explicit time-stepping schemes used in GO-MELT, the ML-VMS with space-time C-HiDeNN-TD is not constrained by the stability condition. Combined with higher-order temporal C-HiDeNN elements, much larger steps can be used to ensure both efficiency and accuracy. This is particularly beneficial for simulations with long durations as in LPBF. More importantly, in space-time ML-VMS C-HiDeNN-TD, the solution is obtained simultaneously for all time steps. Therefore, this is also distinct from the semi-discretization formulation used in the time-marching formulation as in \cite{2025arXiv251000516G}.

\begin{figure}[htbp]
\centering
\includegraphics[width=5in]{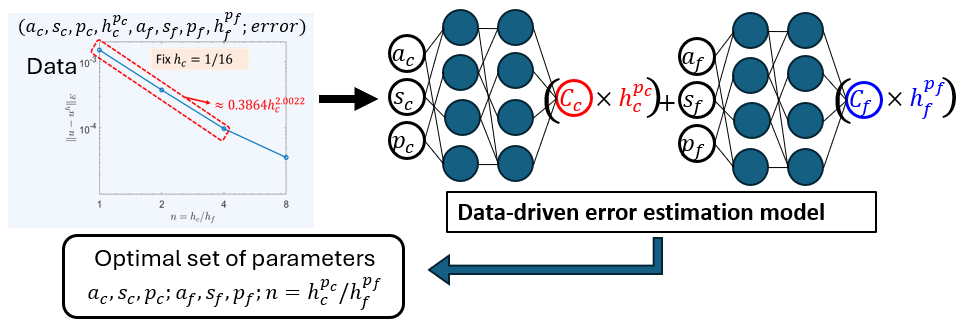}
\caption{The architecture of the data-driven error estimation model for two-level VMS. The model error is the sum of two terms: $C^{(c)} \cdot h_c^{p_c} $ and $ C^{(f)} \cdot h_f^{p_f}$. The coefficients  $C^{(c)}$ and $C^{(f)}$ can be predicted using neural networks.}
\label{fig:DataDrivenErrModel}
\end{figure}

To optimize the hyperparameters in ML-VMS C-HiDeNN-TD ($a_l, s_l, p_l$), we can introduce a data-driven error model where neural networks learn the mapping from a level's parameters to its error contribution (e.g., coefficients $C^{(c)}$ and $C^{(f)}$ in the error estimate $C^{(c)}h_c^{p_c}+C^{(f)}h_f^{p_f}$, see Fig. \ref{fig:DataDrivenErrModel}). The primary advantage of this method is its computational efficiency. Our error estimation theory (Eqs. \ref{eq:ErrEst_reduced} and \ref{eq:ErrEst_nLevel}) justifies training these networks exclusively on coarse-mesh data, as the model is guaranteed to extrapolate effectively to finer meshes. Therefore, this data-driven approach can provide a robust and computationally inexpensive means of optimizing the ML-VMS framework.

%To avoid the complexities of multilevel mesh adaptation, we do not include nodal coordinates as training inputs, although the underlying C-HiDeNN framework inherently supports $r-$adaptivity.

\section{Conclusion and future work}
\label{sec:CONCLUSION}
In this paper, we propose a novel MultiLevel Variational MultiScale (ML-VMS) framework to address a key challenge in engineering simulation: resolving disparate length scales, which can lead to a prohibitively large number of degrees of freedom (DoFs). In the ML-VMS framework, we utilize the C-HiDeNN and C-HiDeNN-TD interpolation functions on linear FE meshes, adjusting the hyperparameters to control the interpolation order at each level. This formulation ensures a balance between computational efficiency and accuracy, and it avoids the mesh incompatibility difficulties that arise across different levels in standard multilevel methods using higher-order FE basis functions.

In the future, the proposed ML-VMS can be easily extended to a multiphysics setup where each level can be used to account for different physics. For example, in high-fidelity LPBF simulation, the Navier-Stokes equation can be used at the finest level to account for complex interaction among laser, fluid and gas in the melt-pool region. Concurrently, the coarse level could model the heat transfer in the surrounding solid material using a transient heat transfer equation to improve computational efficiency.

Moreover, the capability of ML-VMS framework can be further extended to solve parametric PDEs by using the space-parameter-time (S-P-T) C-HiDeNN-TD basis function. In this case, the PDE coefficients such as laser power, absorptivity, and velocity can be treated as input variables \cite{guo2025tensor}. In contrast to standard data-driven surrogate models that require offline data generation—a process that can be prohibitively expensive for large-scale, high-dimensional parametric problems such as LPBF—the S-P-T C-HiDeNN-TD model is inherently data-free. It operates by directly finding the parameterized solution through the minimization of the parameterized PDE's residual. Consequently, it stands as a promising digital twin model for the LPBF process optimization and control.

Finally, despite the versatility of the ML-VMS framework, its numerical implementation for different PDEs and future extension to parameterized systems can still be a cumbersome task due to the growing complexity of the problems. Recent advances in large-language models (LLMs) for mathematical reasoning and coding show great promise in streamlining the automatic derivation, coding, and verification \& validation of complex numerical algorithms \cite{guo2025large}. By learning the fundamentals of ML-VMS, a specialized ``mechanistic LLM" could be developed to greatly alleviate the mathematical derivation and code implementation challenges and make the ML-VMS framework accessible to a wide range of engineering problems such as fatigue analysis and electronic design automation (EDA).

\appendix

\section{Construction of convolution patch functions}
\label{appenx:ConvolutionPatch}

In C-HiDeNN, the convolution patch function $\bm{W}_{s,a,p}^{\bm{x}_I}  (x)$ plays a similar role as the kernel function in convolution neural networks, which connects the neighboring nodes to the centered node $I$. It can be arbitrary function that satisfies the reproducing property and Kronecker delta property in its corresponding nodal convolution patch $A^{\bm{x}_I}_s$ with a patch size $s$. In this work, we take the radial basis functions as the convolution patch functions. As a result, $\bm{W}_{s,a,p}^{\bm{x}_I}  (x)$ is defined in the following matrix form:
\begin{equation}
    \bm{W}_{s,a,p}^{\bm{x}_I}  (\bm{x})=\bm{\Psi}(\bm{x})\bm{A} + \bm{P}(\bm{x}) \bm{K},
\end{equation}
where $\bm{W}_{s,a,p}^{\bm{x}_I}  (x)=\{W_{s,a,p,K}^{\bm{x}_I} (\bm{x})|K\in A_s^{\bm{x}_I} \}$ is a vector of convolution patch functions defined in $A_s^{\bm{x}_I}$, and $\bm{\Psi}(\bm{x})=\{\Psi_J (\bm{x})|J\in A_s^{\bm{x}_I} \}$ is a vector of radial basis functions. A radial basis function $\Psi_J (\bm{x})$ is a function of radial distance between a point of interest $\bm{x}$ and a node $J$ in $A_s^{\bm{x}_I}$. That is, $\Psi_J=\psi(\|\bm{x}-\bm{x}_J\|_2)$, where $\psi(z)$ is taken as a cubic spline function:
\begin{equation} \label{eq:CubicSpline}
    \psi(z)=\begin{cases}
    \frac{2}{3} -4z^2+4z^3 \quad \forall z \in [0, \frac{1}{2}]\\
    \frac{4}{3} -4z+4z^2- \frac{4}{3}z^3 \quad \forall z \in (\frac{1}{2}, 1]\\
    0 \quad \forall z \in (1,+\infty)
   \end{cases}
\end{equation}
$\bm{P}=[P_1(\bm{x}), P_2(\bm{x}), \cdots, P_m(\bm{x})]$ is a vector of functions to be reproduced under reproducing conditions such as polynomials. For example, in 2D, the second order reproducing condition ($p = 2$) requires $m = 6$ polynomial basis
functions: $1, x, y, x^2, xy, y^2$. $\bm{A}$ and $\bm{K}$ are the coefficients to construct $\bm{W}_{s,a,p}^{\bm{x}_I}  (x)$ through linear combinations of $\bm{\Psi}(\bm{x})$ and $\bm{P}(\bm{x})$. They are determined by imposing the Kronecker delta property:
\begin{equation}
    W_{s,a,p,K}^{\bm{x}_I} (\bm{x}_J) = \delta_{JK}, \text{for } J, K\in A_s^{\bm{x}_I}
\end{equation}
and reproducing conditions:
\begin{equation} \label{eq:ReproducingConditions}
    \sum_{J\in A_s^{\bm{x}_I}} W_{s,a,p,J}^{\bm{x}_I}(\bm{x}) P_k(\bm{x}_J) = P_k(\bm{x}), k=1,2,\cdots,m.
\end{equation}

When the dimension of $\bm{P}(\bm{x})$ is smaller than the dimension of $\bm{\Psi}(\bm{x})$, i.e., $n_s>m$, where $n_s$ is the number of nodes in $A_s^{\bm{x}_I}$, the solution for coefficients is not unique. We take the coefficients as
\begin{equation} \label{eq:RadialCoefficients}
    \bm{A}=\bm{R}^{-1}(\bm{I}_{n_s \times n_s}-\bm{QK}), \bm{K}=\left( \bm{Q}^T \bm{R}^{-1} \bm{Q} \right)^{-1} \bm{Q}^T \bm{R}^{-1},
\end{equation}
that satisfies the Kronecker delta and reproducing properties. $\bm{I}_{n_s \times n_s}$ denotes an $n_s$-dimensional identity matrix. The moment matrices $Q$ and $R$ are matrices composed of the nodal values of $\bm{\Psi}$ and $\bm{P}$: $\bm{R}=\{ \Psi_J(\bm{x}_I) \}_{I,J}$ and $\bm{Q}=\{ P_j(\bm{x}_I) \}_{I,j}$.

\section{Lower bound for the number of nodes in the nodal patch}

To ensure the reproducing conditions in Eq. (\ref{eq:ReproducingConditions}), the number of nodes $n_s$ in the nodal patch $A_s^{\bm{x}_I}$ must be equal to or greater than $m$, i.e., $n_s \geq m$. If $\bm{P}(\bm{x})$ is composed of standard $d$-dimensional polynomial basis functions up to the $p$-th order, we have $m=C^d_{p+d}$. If $\bm{P}(\bm{x})$ is composed of the product of one-dimensional polynomial basis functions, we have $m=(p+1)^d$. For regular mesh, there are $n_s=(2s+1)^d$ nodes in the nodal patch. A rough lower bound to ensure the reproducing property is $s \geq p/2$.

In particular, if $n_s = m$, the coefficients $\bm{A}$ vanishes, i.e., the convolution patch functions $\bm{W}_{s,a,p}^{\bm{x}_I}  (\bm{x})$ are independent of $\bm{\Psi}(\bm{x})$. We have the following theorem:

\begin{theorem}
Let the coefficients $\bm{A}$ and $\bm{K}$ be defined by Eq. (\ref{eq:RadialCoefficients}). If $n_s = m$, we have:
\begin{equation}
    \bm{A} = \bm{0}, \bm{K} = \bm{Q}^{-1}.
\end{equation}
\end{theorem}

\begin{proof}
    When $n_s = m$, the matrix $\bm{Q}$ becomes an $m \times m$ invertible matrix. Therefore,
    \begin{equation}
        \bm{K}=\left( \bm{Q}^{-1} \bm{R} \bm{Q}^{-T} \right) \bm{Q}^T \bm{R}^{-1} = \bm{Q}^{-1}.
    \end{equation}
    This yields $\bm{A} = \bm{R}^{-1}\bm{0} = \bm{0}$.
\end{proof}

\begin{corollary}
Let $\bm{P}(\bm{x})$ be polynomial basis and $n_s=m$, the convolution patch functions $\bm{W}_{s,a,p}^{\bm{x}_I}  (\bm{x})$ reduces to Lagrange polynomials.
\end{corollary}

\section{The treatment of boundary nodal convolution patch}
\label{appenx:BoundaryPatch}

\begin{figure}[htbp]
\centering
\includegraphics[width=4in]{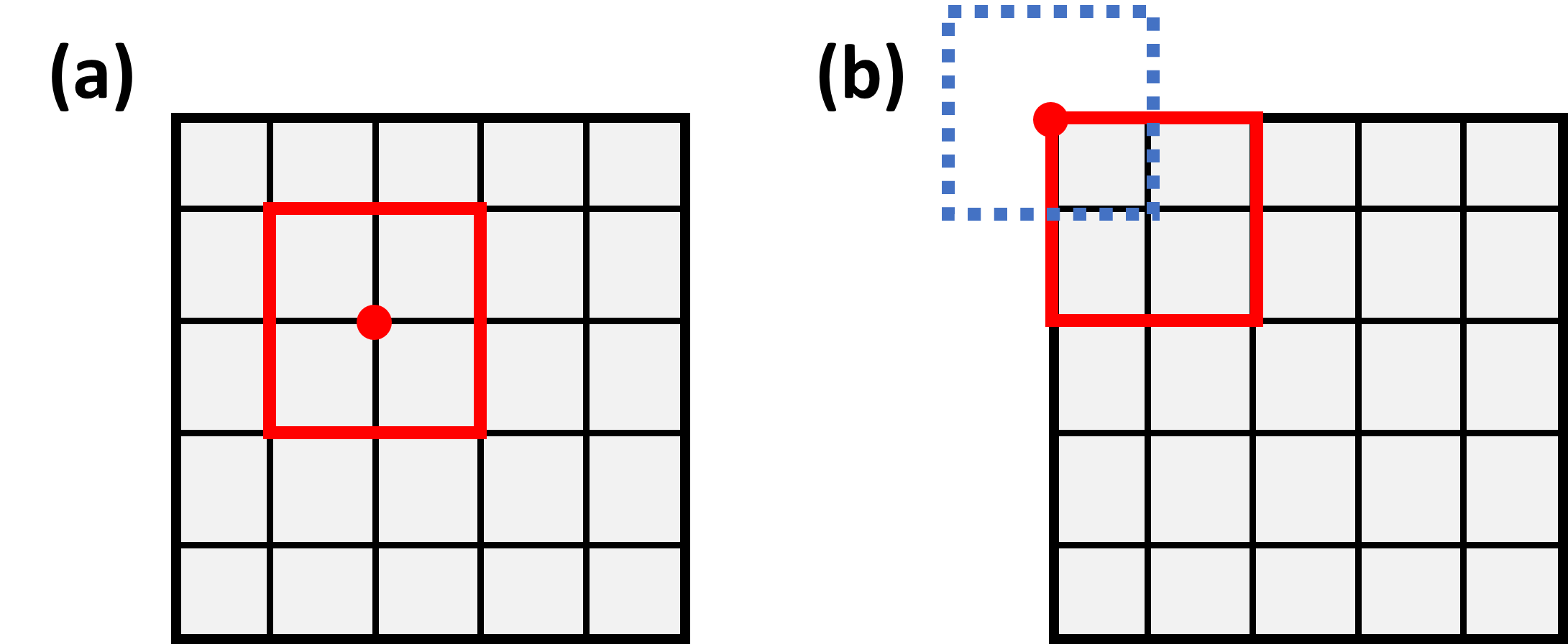}
\caption{Illustration for a nodal convolution patch domain on a 2D mesh: (a) patch domain for an interior node; (b) patch domain for a boundary node.}
\label{fig:BoundaryTreatment}
\end{figure}

In general, the node $I$ is located at the center of nodal convolution patch $A_s^{\bm{x}_I}$, as illustrated in Fig. \ref{fig:BoundaryTreatment} (a). However, this definition of the patch is problematic for elements near the boundary of the computational domain, because there are not enough neighboring nodes for a given patch size (original $A_s^{\bm{x}_I}$ represented by a blue dashed box). To address this, we move the nodal convolution patch into the interior domain, i.e., the red solid box in Fig. \ref{fig:BoundaryTreatment} (b). Therefore, the node $I$ may not be located at the center in this new nodal convolution patch $A_s^{\bm{x}_I}$.

\section{Comparison among FEM, IGA, PINN and C-HiDeNN}
\label{appenx:Comparison}

A comparison of C-HiDeNN with other numerical methods, including FEM, isogeometric analysis (IGA) \cite{hughes2005isogeometric} and physics-informed neural networks (PINN) \cite{raissi2019physics}, is listed in Table \ref{table:Compare_C-HiDeNN}. Compared to IGA and PINN, C-HiDeNN does not require any special treatment for the imposition of exact boundary conditions due to its inherent interpolation properties. By increasing the reproducing polynomial order $p$ and reconstructing the convolution patch functions $W_{s,a,p}^{\bm{x}_I} (\bm{x})$ accordingly, C-HiDeNN can build higher order shape functions using only linear elements (3- or 4-node elements in 2D and 4- or 8-node elements in 3D), achieving faster convergence rates than FEM without increasing DoFs. Although C-HiDeNN naturally admits $r$-adaptivity by optimizing nodal coordinates \cite{liu2023hidenn}, we will not discuss $r$-adaptivity in this work and focus on controlling other hyperparameters, including the convolution patch size $s$ and reproducing polynomial order $p$.

\begin{table}[]
    \centering
    \begin{tabularx}{\textwidth}{|X|X|X|X|X|}
    \hline
         & FEM & IGA & PINN & C-HiDeNN \\
        \hline
        \makecell[l]{DoFs \\ attached to} & nodes & \makecell[l]{control \\ variables} & weights, biases & nodes \\
        \hline
        Interpolation property & Yes & In general, No & No & Yes \\
        \hline
        Reproducing property & polynomials up to order p & exact geometry & No & any nonlinear activation \\
        \hline
        Adaptivity & $r,h,p$ & $h,k$ & neurons, layers, activations & \makecell[l]{$r,h,p,s,a$, \\activations} \\
        \hline
    \end{tabularx}
    \caption{Comparison among FEM, IGA, physics-informed neural network (PINN) and C-HiDeNN.}
    \label{table:Compare_C-HiDeNN}
\end{table}

% \begin{figure}[htbp]
% \centering
% \includegraphics[width=0.9\textwidth]{ratio.png}
% \caption{Optimal ratio based on accuracy and efficiency (a) case 1: $p_c=3, s_c=2; p_f=1, s_f=1$; (b) case 2: $p_c=3, s_c=2; p_f=3, s_f=2$}
% \label{fig:ratio}
% \end{figure}

% \begin{figure}[htbp]
% \centering
% \includegraphics[width=1.0\textwidth]{ratio2.png}
% \caption{Accuracy and efficiency for C-HiDeNN with hyperparameters (a) $n=2$ (b) $n=4$ (c) $n=8$}
% \label{fig:ratio2}
% \end{figure}

\section{Proof of Theorem 1 and Theorem 2}
We first define a two-level interpolation of the real solution $u$, given by
\begin{equation}
    v^h(\bm{x}) = \left\{
    \begin{array}{cc}
        v_c(\bm{x}) & \bm{x}\in \Omega_c \backslash \Omega_f, \\
        v_f(\bm{x}) & \bm{x}\in \Omega_f,
    \end{array}
    \right.
\end{equation}
with 
\begin{equation}
    v_c(\bm{x}) = \mathcal{I}^{(c)} u, v_f(\bm{x}) = \mathcal{I}^{(f)} u.
\end{equation}
Note that $v_f'=v_f-v_c$ might not vanish on the boundary $\Gamma$, and therefore $v^h$ might be discontinuous across $\Gamma$.

\begin{lemma}
There exists an interpolation function $\tilde{v}_f' \in \mathcal{V}^h_f$ such that the deviation between $v_f'$ and $\tilde{v}_f'$ can be estimated by
\begin{equation}
    \| v_f'-\tilde{v}_f' \|_{H^1(\Omega_f)} \leq c_{1,\Gamma} h_c^{p_c} \|u\|_{H^{p_c+1}(\Gamma)} + c_{2,\Gamma} h_f^{p_f} \|u\|_{H^{p_f+1}(\Gamma)},
\end{equation}
where $c_{1,\Gamma}$ and $c_{2,\Gamma}$ are positive constants.
\end{lemma}

\begin{proof}
We take a simple 2D case for illustration first: $\Omega_f$ is a 2D regular domain $[a,b]\times[c,d]$ with a uniform mesh. Denote the coordinates $\bm{x}=(x,y)$ and nodes $(x_i,y_j ),i=1,2,…,n_1,j=1,2,…,n_2$, where $n_1$ and $n_2$ are number of nodes along $x,y$-directions, respectively. Denote the value of $v_f'$ on the boundary is $g$, i.e., $g=v_f' |_{\Gamma}$. We construct the deviation $\delta v_f' = v_f'-\tilde{v}_f'$ directly such that it satisfies the boundary condition $\delta v_f'|_{\Gamma}=g$, and it vanishes away the boundary. The detailed constructions are given as follows.

Near the left boundary $\Gamma_{\text{left}}$, $\delta v_f'$ is constructed by product rule:
\begin{equation}
    g (a,y) \cdot \max(1-(x-a)/(\beta h_f),0).
\end{equation}
$\max(1-(x-a)/(\beta h_f),0)$ is a piece-wise linear function, which equals to 1 on the left boundary, and equals to 0 at $x=a+\beta h_f$. The parameter $\beta$ is an integer such that $\beta h_f\in[l_1,l_2 ],0<l_1<l_2<\min(b-a, d-c)$, where $l_1$ and $l_2$ are constants independent of $h_f$. That is, when the mesh is finer, $\beta$ will increase. Hence, the function $g (a,y) \cdot \max(1-(x-a)/(\beta h_f),0)$ can be reproduced by interpolations in $\Omega_f$. Its $H^1$-norm can be estimated by
\begin{eqnarray}
    && \|g (a,y) \cdot \max(1-\dfrac{x-a}{\beta h_f},0)\|_{H^1(\Omega_f)} \\ \nonumber
    &=& \left[ \int_c^d \left(\dfrac{\beta h_f}{3} (g (a,y))^2 + \dfrac{\beta h_f}{3} (\dfrac{\partial}{\partial y}g (a,y))^2 + \dfrac{1}{\beta h_f}(g (a,y))^2 \right) \mathrm{d}y \right]^{1/2} \\ \nonumber
    &\leq& \left(\dfrac{l_2}{3}+\dfrac{1}{l_1} \right)^{1/2} \|g (a,y)\|_{H^1(\Gamma_{\text{left}})}.
\end{eqnarray}

$\delta v_f'$ near other boundaries is constructed in the same manner. Thus $\delta v_f'$ is expressed by
\begin{eqnarray} \delta v_f' &=& g(a,y) \cdot \max(1-(x-a)/(\beta h_f), 0) + g(b,y) \cdot \max(1-(b-x)/(\beta h_f), 0) \\ \nonumber && + g(x,c) \cdot \max(1-(y-c)/(\beta h_f), 0) + g(x,d) \cdot \max(1-(d-y)/(\beta h_f), 0) \\ \nonumber && - g(a,c) \max(1-(x-a)/(\beta h_f), 0) \cdot \max(1-(y-c)/(\beta h_f), 0) \\ \nonumber && - g(a,d) \max(1-(x-a)/(\beta h_f), 0) \cdot \max(1-(d-y)/(\beta h_f), 0) \\ \nonumber && - g(b,c) \max(1-(b-x)/(\beta h_f), 0) \cdot \max(1-(y-c)/(\beta h_f), 0) \\ \nonumber && - g(b,d) \max(1-(b-x)/(\beta h_f), 0) \cdot \max(1-(d-y)/(\beta h_f), 0). \end{eqnarray}
Notably, the sum of the first four terms is corrected by subtracting the corner parts they share, so that the shared parts are not calculated twice. As the coarse mesh is covered by the fine mesh and C-HiDeNN has the interpolation property, the four verticies are repeating nodes and the value of $g$ equals to zero at these points. The $H^1$-norm of $\delta v_f'$ can be estimated by
\begin{eqnarray}
    \|\delta v_f'\|_{H^1(\Omega_f)} &\leq& \left(\dfrac{l_2}{3}+\dfrac{1}{l_1} \right)^{1/2}\| g\|_{H^1(\Gamma)} \\ \nonumber
    &=& \left(\dfrac{l_2}{3}+\dfrac{1}{l_1} \right)^{1/2}\| \mathcal{I}^{(f)}u-u+u-\mathcal{I}^{(c)}u\|_{H^1(\Gamma)} \\ \nonumber
    &\leq& \left(\dfrac{l_2}{3}+\dfrac{1}{l_1} \right)^{1/2} \left(\| \mathcal{I}^{(f)}u-u\|_{H^1(\Gamma)}+\|u-\mathcal{I}^{(c)}u\|_{H^1(\Gamma)}\right) \\ \nonumber
    &\leq& c_{1,\Gamma} h_c^{p_c} \|u\|_{H^{p_c+1}(\Gamma)} + c_{2,\Gamma} h_f^{p_f} \|u\|_{H^{p_f+1}(\Gamma)}.
\end{eqnarray}
The last inequality comes from the $p_c$-th order interpolation estimate (reproducing property) for coarse interpolations and the $p_f$-th order interpolation estimate (reproducing property) for fine interpolations.

$\tilde{v}_f'$ is obtained by subtracting the constructed $\delta v_f'$ from $v_f'$. Above conclusions can be extended to arbitrary $\Omega_f$ with a non-uniform mesh.

\end{proof}

Then we prove Theorem 1.

\begin{proof}
 According to Lemma 1, there exists an interpolation function $\tilde{v}_f' \in \mathcal{V}^h_f$ such that the deviation between $v_f'$ and $\tilde{v}_f'$ can be estimated by
\begin{equation}
    \| v_f'-\tilde{v}_f' \|_{H^1(\Omega_f)} \leq c_{1,\Gamma} h_c^{p_c} \|u\|_{H^{p_c+1}(\Gamma)} + c_{2,\Gamma} h_f^{p_f} \|u\|_{H^{p_f+1}(\Gamma)},
\end{equation}
where $c_{1,\Gamma}$ and $c_{2,\Gamma}$ are positive constants.

The real solution also satisfies the weak forms Eqs. (\ref{eq:Two-Level_Weak1}, \ref{eq:Two-Level_Weak2}), so leading to
\begin{eqnarray} \label{eq_appenx:weak1}
    a(w_c,u-u^h )_{\Omega_c} = 0, &&\forall w_c \in \mathcal{V}_c^h, \\
    \label{eq_appenx:weak2}
    a(w_f,u-u^h )_{\Omega_f} =0, &&\forall w_f \in \mathcal{V}_f^h.
\end{eqnarray}
This implies that
\begin{equation} \label{eq_appenx:weak3}
    a(w_c,u-u^h )_{\Omega_f}+a(w_f,u-u^h )_{\Omega_f} =0, \forall w_c \in \mathcal{V}_c^h, w_f \in \mathcal{V}_f^h.
\end{equation}
in the overlapping domain $\Omega_f$.

We have
\begin{eqnarray} \label{eq:proof1}
    && a(u-v^h, u-v^h)_{\Omega} \\ \nonumber
    &=& a(u-u^h, u-u^h)_{\Omega} + a(u^h-v^h, u^h-v^h)_{\Omega} + 2a(u-u^h, u^h-v^h)_{\Omega} \\ \nonumber
    &=& a(u-u^h, u-u^h)_{\Omega} + a(u^h-v^h, u^h-v^h)_{\Omega} + 2a(u-u^h, u_c-v_c)_{\Omega} \\ \nonumber
    && + 2a(u-u^h, u_f'-\tilde{v}_f')_{\Omega_f} + 2a(u-u^h, \tilde{v}_f'-v_f')_{\Omega_f}
\end{eqnarray}
Here, we define $u_c = \mathcal{I}^{(c)}u_f$ and $u_f'=u_f - \mathcal{I}^{(c)}u_f$ in the overlapping domain $\Omega_f$.

Due to the arbitrariness of $w_c$ and $w_f$ in Eqs. (\ref{eq_appenx:weak1}, \ref{eq_appenx:weak2}, \ref{eq_appenx:weak3}), we have 
\begin{equation}
    a(u-u^h, u_c-v_c)_{\Omega} +
    a(u-u^h, u_f'-\tilde{v}_f')_{\Omega_f} =0.
\end{equation}
According to Cauchy-Schwarz inequality, the last term in Eq. (\ref{eq:proof1}) is estimated by
\begin{eqnarray}
    && a(u-u^h, \tilde{v}_f'-v_f')_{\Omega_f} \\ \nonumber
    &\geq& - \left(a(u-u^h, u-u^h)_{\Omega_f}\right)^{1/2} \cdot \left(a(\tilde{v}_f'-v_f', \tilde{v}_f'-v_f')_{\Omega_f}\right)^{1/2} \\ \nonumber
    &\geq& - \left(a(u-u^h, u-u^h)_{\Omega_f}\right)^{1/2} \cdot c_2 \|\tilde{v}_f'-v_f' \|_{H^1(\Omega_f)} \\ \nonumber
    &\geq& - \left(a(u-u^h, u-u^h)_{\Omega_f}\right)^{1/2} \cdot c_2 \left[ c_{1,\Gamma} h_c^{p_c} \|u\|_{H^{p_c+1}(\Gamma)} + c_{2,\Gamma} h_f^{p_f} \|u\|_{H^{p_f+1}(\Gamma)} \right]
\end{eqnarray}

Therefore, we obtain
\begin{eqnarray}
    && a(u-v^h, u-v^h)_{\Omega} \\ \nonumber
    &\geq& a(u-u^h, u-u^h)_{\Omega} - 2 \left(a(u-u^h, u-u^h)_{\Omega_f}\right)^{1/2} \cdot c_2 \left[ c_{1,\Gamma} h_c^{p_c} \|u\|_{H^{p_c+1}(\Gamma)} + c_{2,\Gamma} h_f^{p_f} \|u\|_{H^{p_f+1}(\Gamma)} \right].
\end{eqnarray}
This implies that
\begin{eqnarray}
    && \left( a(u-u^h, u-u^h)_{\Omega} \right)^{1/2} \\ \nonumber
    &\leq& \left[ a(u-v^h, u-v^h)_{\Omega} + \left( c_2 c_{1,\Gamma} h_c^{p_c} \|u\|_{H^{p_c+1}(\Gamma)} + c_2 c_{2,\Gamma} h_f^{p_f} \|u\|_{H^{p_f+1}(\Gamma)}\right)^2 \right]^{1/2} \\ \nonumber
    && + \left( c_2 c_{1,\Gamma} h_c^{p_c} \|u\|_{H^{p_c+1}(\Gamma)} + c_2 c_{2,\Gamma} h_f^{p_f} \|u\|_{H^{p_f+1}(\Gamma)}\right) \\ \nonumber
    &\leq& \left( a(u-v^h, u-v^h)_{\Omega} \right)^{1/2} + 2\left( c_2 c_{1,\Gamma} h_c^{p_c} \|u\|_{H^{p_c+1}(\Gamma)} + c_2 c_{2,\Gamma} h_f^{p_f} \|u\|_{H^{p_f+1}(\Gamma)}\right).
\end{eqnarray}

According to the error estimate for C-HiDeNN interpolants (Eq. (\ref{eq:C-HiDeNNInterpError})), the $p_c$-th order interpolation estimate for $u-v_h=u-\mathcal{I}^{(c)}u$ in $\Omega_c \backslash \Omega_f$ and the $p_f$-th order interpolation estimate for $u-v_h=u-\mathcal{I}^{(f)}u$ in $\Omega_f$, we have the following error estimate:
\begin{equation}
    \|u-u^h \|_{H^1 (\Omega) }\leq C_1  h_c^{p_c}  \|u\|_{H^{p_c+1} (\Omega_c \backslash  \Omega_f ) }+C_2  h_f^{p_f}  \|u\|_{H^{p_f+1} (\Omega_f ) }+ \left(C_3  h_c^{p_c } \|u\|_{H^{p_c+1} (\Gamma) }+C_4 h_f^{p_f}  \|u\|_{H^{p_f+1} (\Gamma ) } \right),
\end{equation}
where $C_1, C_2, C_3$ and $C_4$ are constants.

\end{proof}

Finally, we prove Theorem 2.

\begin{proof}
In the same manner, we can define an $m$-level interpolation of the real solution $u$, given by
\begin{equation}
    v^h(\bm{x})=\left\{
    \begin{array}{cc}
        v_1(\bm{x}), & \bm{x}\in \Omega_1 \backslash \Omega_2, \\
        v_2(\bm{x}), & \bm{x}\in \Omega_2 \backslash \Omega_3, \\
        \cdots & \\
        v_m(\bm{x}), & \bm{x}\in \Omega_m.
    \end{array}
    \right.
\end{equation}
with 
\begin{equation}
    v_l(\bm{x}) = \mathcal{I}^{(l)} u, l=1,2,\cdots, m.
\end{equation}
Define $v_l'=v_l-v_{l-1}, l>1$. Similarly, according to Lemma 1, there exists a new interpolation function $\tilde{v}_l' \in \mathcal{S}_l^h$ such that the deviation between $v_l'$ and $\tilde{v}_l'$ can be estimated by
\begin{equation}
    \| v_l'-\tilde{v}_l' \|_{H^1(\Omega_l \backslash \Omega_{l+1})} \leq c_{1,\Omega_l} h_{l-1}^{p_{l-1}} \|u\|_{H^{p_{l-1}+1}(\partial\Omega_l)} + c_{2,\Omega_l} h_l^{p_l} \|u\|_{H^{p_l+1}(\partial\Omega_l)}, l>1
\end{equation}
where $c_{1,\Omega_l}$ and $c_{2,\Omega_l}$ are positive constants. Here, we define $\Omega_{m+1}=\emptyset$, yielding $\Omega_m \backslash \Omega_{m+1}=\Omega_m$. Therefore, we define a modified $m$-level interpolation:
\begin{equation}
    \tilde{v}^h(\bm{x})=\left\{
    \begin{array}{cc}
        v_1(\bm{x}), & \bm{x}\in \Omega_1 \backslash \Omega_2, \\
        v_2(\bm{x})-v_2'+\tilde{v}_2', & \bm{x}\in \Omega_2 \backslash \Omega_3, \\
        \cdots & \\
        v_m(\bm{x})-v_m'+\tilde{v}_m', & \bm{x}\in \Omega_m,
    \end{array}
    \right.
\end{equation}
which is continuous across all the interfaces.

As the real solution satisfies the weak form Eq. (\ref{eq:nlevel_weak}), we have
\begin{equation} \label{eq_appendx:nlevel_weak}
    a(w_l, u-u^h)_{\Omega_l} = 0, \forall w_l \in \mathcal{V}_l^h, l=1,2,\cdots,n.
\end{equation}

We have
\begin{eqnarray} \label{eq:proof1_theorem2}
    && a(u-v^h, u-v^h)_{\Omega} \\ \nonumber
    &=& a(u-u^h, u-u^h)_{\Omega} + a(u^h-v^h, u^h-v^h)_{\Omega} + 2a(u-u^h, u^h-v^h)_{\Omega} \\ \nonumber
    &=& a(u-u^h, u-u^h)_{\Omega} + a(u^h-v^h, u^h-v^h)_{\Omega} + 2a(u-u^h, u^h-\tilde{v}^h)_{\Omega} + 2a(u-u^h, \tilde{v}^h-v^h)_{\Omega}.
\end{eqnarray}

Due to the arbitrariness of weighting functions in Eq. (\ref{eq_appendx:nlevel_weak}), we have 
\begin{equation}
    a(u-u^h, u^h-\tilde{v}^h)_{\Omega} =0.
\end{equation}
According to Cauchy-Schwarz inequality, the last term in Eq. (\ref{eq:proof1_theorem2}) is estimated by
\begin{eqnarray}
    && a(u-u^h, \tilde{v}^h-v^h)_{\Omega_f} \\ \nonumber
    &=& \sum_{l=2}^{m} a(u-u^h, \tilde{v}_l'-v_l')_{\Omega_l \backslash \Omega_{l+1}}\\ \nonumber
    &\geq& - \sum_{l=2}^{m}\left(a(u-u^h, u-u^h)_{\Omega_l \backslash \Omega_{l+1}}\right)^{1/2} \cdot c_2 \|\tilde{v}_l'-v_l' \|_{H^1(\Omega_l \backslash \Omega_{l+1})} \\ \nonumber
    &\geq& - \left(a(u-u^h, u-u^h)_{\Omega}\right)^{1/2} \cdot c_2 \sum_{l=2}^m \left[ c_{1,\Omega_l} h_{l-1}^{p_{l-1}} \|u\|_{H^{p_{l-1}+1}(\partial\Omega_l)} + c_{2,\Omega_l} h_l^{p_l} \|u\|_{H^{p_l+1}(\partial\Omega_l)} \right]
\end{eqnarray}

Therefore, we obtain
\begin{eqnarray}
    && a(u-v^h, u-v^h)_{\Omega} \\ \nonumber
    &\geq& a(u-u^h, u-u^h)_{\Omega} - 2 \left(a(u-u^h, u-u^h)_{\Omega}\right)^{1/2} \cdot c_2 \sum_{l=2}^m \left[ c_{1,\Omega_l} h_{l-1}^{p_{l-1}} \|u\|_{H^{p_{l-1}+1}(\partial\Omega_l)} + c_{2,\Omega_l} h_l^{p_l} \|u\|_{H^{p_l+1}(\partial\Omega_l)} \right],
\end{eqnarray}
yielding
\begin{equation}
    \left( a(u-u^h, u-u^h)_{\Omega} \right)^{1/2} \leq \left( a(u-v^h, u-v^h)_{\Omega} \right)^{1/2} + 2c_2 \left[ c_{1,\Omega_l} h_{l-1}^{p_{l-1}} \|u\|_{H^{p_{l-1}+1}(\partial\Omega_l)} + c_{2,\Omega_l} h_l^{p_l} \|u\|_{H^{p_l+1}(\partial\Omega_l)} \right].
\end{equation}

According to the error estimate for C-HiDeNN interpolants (Eq. (\ref{eq:C-HiDeNNInterpError})), we have the following error estimate:
\begin{eqnarray}
    \|u-u^h \|_{H^1 (\Omega) }&\leq& c_2/c_1 \|u-v^h \|_{H^1 (\Omega) } + 2c_2/c_1 \left[ c_{1,\Omega_l} h_{l-1}^{p_{l-1}} \|u\|_{H^{p_{l-1}+1}(\partial\Omega_l)} + c_{2,\Omega_l} h_l^{p_l} \|u\|_{H^{p_l+1}(\partial\Omega_l)} \right] \\ \nonumber
    &\leq& c_2/c_1 \left[ \sum_{l=1}^{m-1} C(a_l, s_l, p_l) h_l^{p_l} \|u\|_{H^{p_{l}+1}(\Omega_l\backslash \Omega_{l+1})} + C(a_m, s_m, p_m) h_m^{p_m} \|u\|_{H^{p_{m}+1}(\Omega_m)} \right] \\ \nonumber
    && + 2c_2/c_1 \left[ c_{1,\Omega_l} h_{l-1}^{p_{l-1}} \|u\|_{H^{p_{l-1}+1}(\partial\Omega_l)} + c_{2,\Omega_l} h_l^{p_l} \|u\|_{H^{p_l+1}(\partial\Omega_l)} \right] \\ \nonumber
    &\leq& \sum_{l=1}^m C^{(l)} (a_l,s_l,p_l;u,\Omega_l,\Omega_{l+1}, \partial \Omega_l, \partial \Omega_{l+1}) \cdot h_l^{p_l},
\end{eqnarray}
where coefficients $C^{(l)}, l=1,2,\cdots,m$ are independent of element sizes, and depend on hyperparameters $a_l, s_l$ and $p_l$, tailored to a specific problem and geometry for the $m$ levels of meshes ($\Omega_1,\Omega_2,\ldots,\Omega_m$).

\end{proof}

\section{Discussions on optimal element size ratios}
\label{appenx:OptElemRatio}

The error estimate of two-level VMS is influenced by two levels of refinements, specifically the coarse- and fine-level element sizes $h_c$ and $h_f$. The refinement of $h_c$ or $h_f$ can lead to increasing computational cost, but may not improve the accuracy significantly, if the error is dominated by the other level of hyperparameters. For example, if $C^{(c)} \cdot h_c^{p_c} > C^{(f)} \cdot h_f^{p_f}$, i.e., the error is dominated by coarse-level parameters, a coarser $h_f$ will influence the error slightly. Therefore, we should choose the coarsest $h_f$ to $C^{(c)} \cdot h_c^{p_c} \approx C^{(f)} \cdot h_f^{p_f}$ to ensure the optimal  efficiency that the accuracy remains unchanged. In contrast, if $C^{(f)} \cdot h_f^{p_f} > C^{(c)} \cdot h_f^{p_c}$, a coarser $h_c$ is allowed to reduce costs while maintaining accuracy until these two terms are approximately equal. Therefore, to balance accuracy and efficiency, the optimal element size ratio $n=h_c/h_f$ should be achieved such that
\begin{equation}
    C^{(c)} \cdot h_c^{p_c} \approx C^{(f)} \cdot h_f^{p_f}.
\end{equation}
This yields an optimal ratio
\begin{equation}
    n_{opt} = \left\lceil \left( C^{(f)}/C^{(c)} \right)^{1/p_f} h_c^{(p_f-p_c)/p_f} \right\rceil.
\end{equation}
To ensure an integer $n$, we take the ceil function. 

If we fix the coarse-level mesh and refine the fine-level mesh, the error decreases rapidly on the order $p_f$ for $n<n_{opt}$, while it decreases slowly for $n>n_{opt}$, as the error is dominated by coarse-level parameters for $n>n_{opt}$. In contrast, if we fix the fine-level mesh and refine the coarse-level mesh, the error decreases rapidly at the order $p_c$ for $n>n_{opt}$, while it decreases slowly for $n<n_{opt}$.

The optimal element size ratio is the product of $\left( C^{(f)}/C^{(c)} \right)^{1/p_c}$ and $h_f^{(p_f-p_c)/p_c}$. The first term is the ratio of coefficients $C^{(f)}$ and $C^{(c)}$ to the power of $1/p_c$, which is independent of element sizes. The key factors in $C^{(f)}$ and $C^{(c)}$ are the norms of $u$ defined in different subdomains: $\|u\|_{H^{p_f+1} (\Omega_f ) }$ and $\|u\|_{H^{p_c+1} (\Omega_c \backslash  \Omega_f )}$, respectively. $\Omega_f$ is the domain of interest with fine-level refinements, where the solution exhibits rapid variation. $\Omega_c \backslash  \Omega_f$ represents the other domain with coarse-level refinements, where the solution varies slowly. Hence, the ratio $C^{(f)}/C^{(c)}$ is generally much greater than 1, leading to a large optimal ratio $n_{opt}$. The second term is the element size $h_f$ to the power of $(p_f-p_c)/p_c$. If $p_c=p_f$, this term reduces to 1, making the optimal ratio independent of element sizes. If $p_f>p_c$, the optimal ratio $n_{opt}$ decreases as the mesh is refined. Conversely, if $p_c>p_f$, the optimal ratio $n_{opt}$ increases with mesh refinement. 

\begin{table}[]
    \centering
    \begin{tabular}{|l|l|l|}
    \hline
        Case & Change in $n_{opt}$ \\
        \hline
        $p_c=p_f$ & Independent of $h_c$ and $h_f$ \\
        \hline
        $p_f > p_c$ & Decrease as refining meshes \\
        \hline
        $p_f < p_c$ & Increase as refining meshes \\
        \hline
    \end{tabular}
    \caption{Change in $n_{opt}$ with the refinement of meshes.}
    \label{table:paramters_Eg1_1}
\end{table}

% The optimal element size ratio generally depends on the refinements when $p_c \neq p_f$. If $p_c=p_f$, the optimal ratio is independent of $h_c$ and $h_f$, and is influenced by other parameters such as $a_c, s_c, a_f, s_f$. 

In the same manner, the optimal set of element sizes $\{h_1, h_2, \ldots, h_n \}$ can be determined to satisfy
\begin{equation}
    C^{(1)} \cdot h_1^{p_1} \approx C^{(2)} \cdot h_2^{p_2} \approx \cdots C^{(n)} \cdot h_n^{p_n}.
\end{equation}
This results in the optimal element size ratio
\begin{equation}
    h_l/h_1 = n_{l, opt} = \left\lceil \left( C^{(l)}/C^{(1)} \right)^{1/p_1} h_l^{(p_l-p_1)/p_1} \right\rceil, l=1,2,\cdots,n.
\end{equation}
If $p_l \neq p_1$, the optimal element size ratio $n_l = h_l/h_1$ depends on the element size $h_l$.

\section{Comparison between C-HiDeNN-PGD and C-HiDeNN-TD and estimation for required number of modes}
\label{appenx:estimateModes}

There are two distinct strategies for solving the unknowns in  Eq. (\ref{eq:TDform}): a greedy solution strategy referred to C-HiDeNN-PGD, and an all-at-once solution strategy referred to C-HiDeNN-TD. In C-HiDeNN-PGD, modes are added and solved sequentially. After solving previous modes, a new mode is added and solved with the previous modes fixed. This process continues until the norm of the added mode falls below a specified tolerance. However, C-HiDeNN-TD solves all modes simultaneously. 

In general, C-HiDeNN-TD requires fewer modes to achieve the same accuracy as C-HiDeNN-PGD; however, it necessitates prior determination of the mode number. We give two ways to estimate the mode number for C-HiDeNN-TD. The first strategy is a rough estimation based on C-HiDeNN-PGD results. The C-HiDeNN-PGD results can be set as a good initial guess for C-HiDeNN-TD. It is essential to note that different solution strategies can lead to different levels of accuracy, meaning that the estimate derived from C-HiDeNN-PGD may not accurately reflect the required number of modes of C-HiDeNN-TD for a given level of accuracy. The second strategy is to estimate the mode number based on the trend in the deviation between C-HiDeNN and C-HiDeNN-TD vs. number of modes for coarse meshes. Numerical tests indicate that this deviation (i.e., the second term in the error decomposition (\ref{eq:ErrorDecomposition})) depends primarily on the mode number and is slightly influenced by mesh refinements \cite{zhang2022hidenn}. Therefore, one can estimate the required number of modes based on this deviation for coarse meshes.

\section{Two-level VMS with different hyperparameters for the elliptic problem and numerical verification of Theorem 1}
\label{appenx:two-levelVMS_elliptic}

\subsection{Two-level VMS with C-HiDeNN}
\begin{figure}[htbp]
\centering
\subfigure{\includegraphics[height=1.8in]{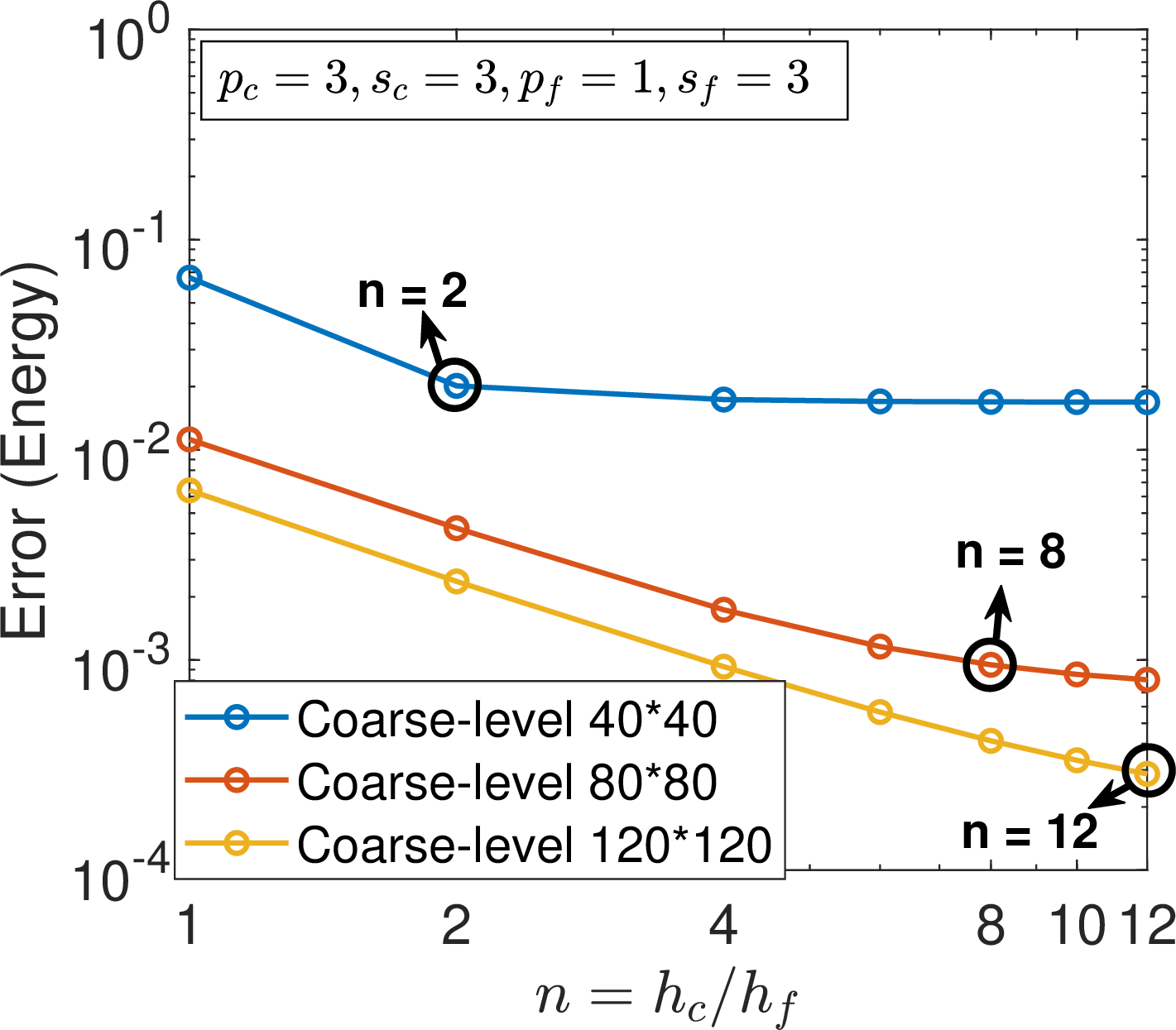}}
\subfigure{\includegraphics[height=1.8in]{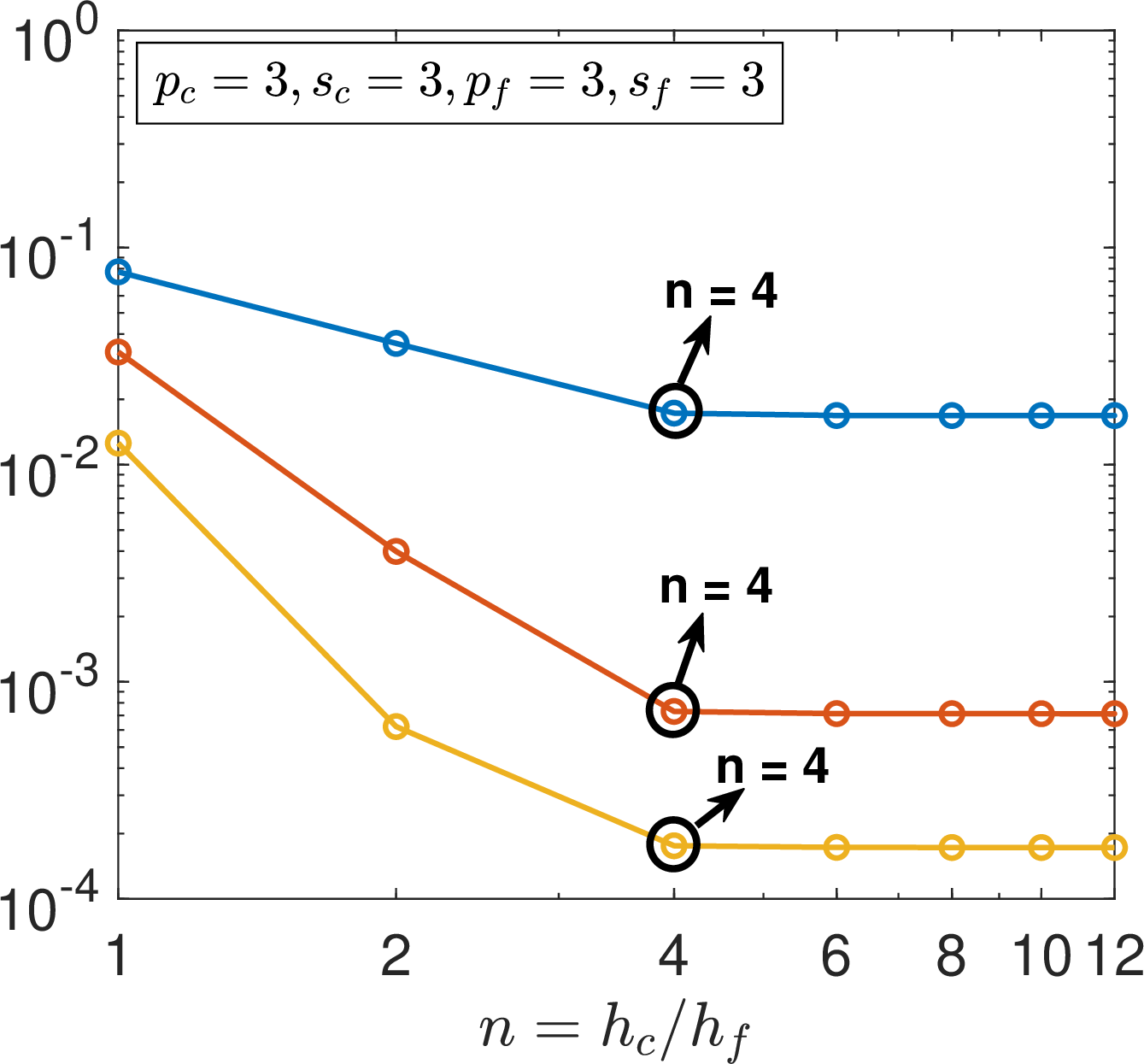}}    
\subfigure{\includegraphics[height=1.8in]{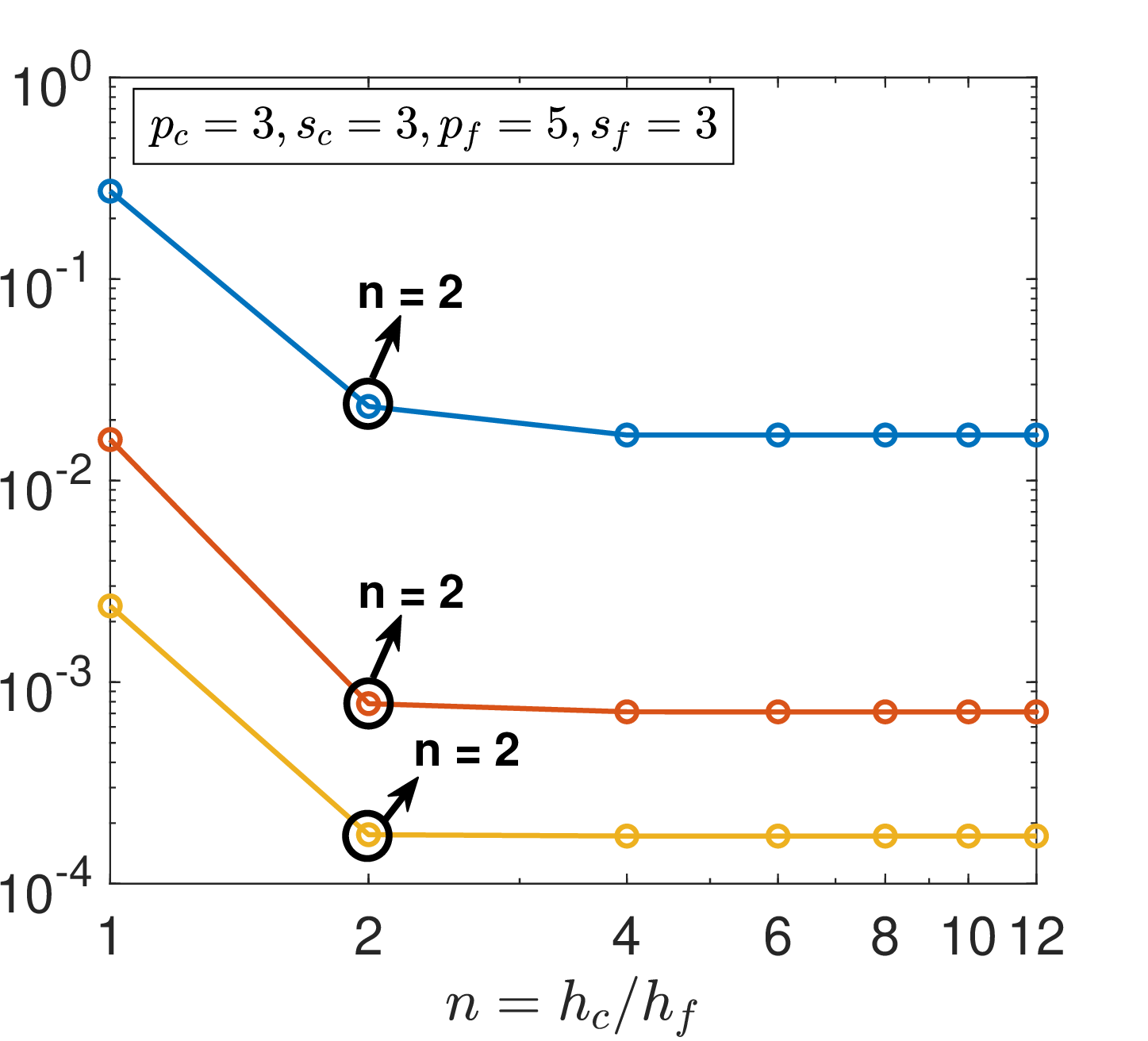}}
\caption{Convergence curves for refining fine-level mesh with fixed coarse-level mesh.}
\label{fig:Convergence_Eg1_Coarse}
\end{figure}

According to the error estimates, two levels of hyperparameters control the error convergence. To verify the theory, we must refine one level of mesh with the other level of mesh fixed. As a result, the error will decrease rapidly initially, and then remain almost constant, as the final error is dominated by the level of approximations for a fixed mesh. We fix the coarse-level mesh, and refine the fine-level mesh. Figure \ref{fig:Convergence_Eg1_Coarse} shows the convergence curves for various coarse-level meshes with three sets of controlling parameters: $(p_c=3, s_c=3, p_f=1, s_f=3)$, $(p_c=3, s_c=3, p_f=3, s_f=3)$ and $(p_c=3, s_c=3, p_f=5, s_f=3)$, corresponding to three cases: $p_f<p_c$, $p_f=p_c$ and $p_f=p_c$. All convergence curves initially exhibit a rapid decrease with the refinement of the fine-level mesh (indicated by a decreasing $n=p_c/p_f$), then followed by a slow decline to constant error levels. However, the turning points are different among three cases. In Fig. \ref{fig:Convergence_Eg1_Coarse} (a), the turning point increases with the refinement of coarse-level mesh, whose trend is consistent with our theoretical optimal element size ratio $n_{opt} = \left\lceil \left( C^{(f)}/C^{(c)} \right)^{1/p_f} h_c^{(p_f-p_c)/p_f} \right\rceil$. In Fig. \ref{fig:Convergence_Eg1_Coarse} (b), the turning points are $n=4$ for all case. In Fig. \ref{fig:Convergence_Eg1_Coarse} (a), the turning points are also consistent. This is because the element size ratio must be an integer and greater or equal to $1$. Since the domain of interest is located in $\Omega_f$, the optimal element size ratio is generally greater than $1$. That is to say, the smallest optimal element size ratio is $n_{opt}=2$ generally. In addition, when $n$ is away from and smaller than turning points, the error decreases with a slope of approximately $p_c$. These observations are consistent with our expectations from the theory of error estimations.

% the convergence curves for different coarse-level meshes are nearly parallel.
% When $n=h_c/h_f<4$, the error decreases with a slope of $p_c=3$. Here, the error is dominated by the fine-level term $C^{f} \cdot h_f^{p_f}$ in the error estimate. When $h_c/h_f>4$, the error decreases slowly. The turning points for various convergence curves are consistent. However, in the left subplot, the turning point increases with the refinement of coarse-level mesh, whose trend is consistent with our theoretical optimal element size ratio $n_{opt} = \left\lceil \left( C^{(f)}/C^{(c)} \right)^{1/p_f} h_c^{(p_f-p_c)/p_f} \right\rceil$. In the right subplot, when $h_c/h_f<2$, the error is dominated by the fine-level term. When $h_c/h_f>2$, the error is almost constant, as it is dominated by the coarse-level parameters.  

% \begin{table}[]
%     \centering
%     \begin{tabular}{|c|c|c|c|}
%     \hline
%         Level & p & s & a \\
%         \hline
%         Coarse & 2 & 3 & 50 \\
%         \hline
%         Fine & 2 & 2& 50 \\
%         \hline
%     \end{tabular}
%     \caption{Parameters for two levels of interpolations.}
%     \label{table:paramters_Eg1}
% \end{table}

% \begin{table}[]
%     \centering
%     \begin{tabular}{|c|c|c|c|}
%     \hline
%         Level & p & s & a \\
%         \hline
%         Coarse & 3 & 2 & 50 \\
%         \hline
%         Fine & 2 & 2& 50 \\
%         \hline
%     \end{tabular}
%     \caption{Parameters for two levels of interpolations.}
%     \label{table:paramters_Eg1_1}
% \end{table}

% Secondly, we fix the coarse-level mesh, and refine the fine-level mesh. The corresponding parameters for interpolations are $p_c=2, s_c=2, p_f=2$ and $s_f=2$. The results for varying $h_f$ are shown in Fig. \ref{fig:Convergence_Eg1_Fine}. The turning points for convergence curves are also located at about $n=2$, which is consistent with the observation from Fig. \ref{fig:Convergence_Eg1_Coarse}. Therefore, the optimal element size ratio for $p_c=p_f=2$ is $n=2$.

% \begin{figure}[htbp]
% \centering
% \includegraphics[width=3in]{RelativeErr_Converge_Refine_h2_p1=2_p2=2.png}

% \caption{Convergence curves for refining fine-level mesh ($p_c=2, p_f=2$).}
% \label{fig:Convergence_Eg1_Fine}
% \end{figure}

% \begin{figure}[htbp]
% \centering
% \includegraphics[width=4in]{ErrVsTime.png}

\begin{figure}[htbp]
\centering
\includegraphics[width=\textwidth]{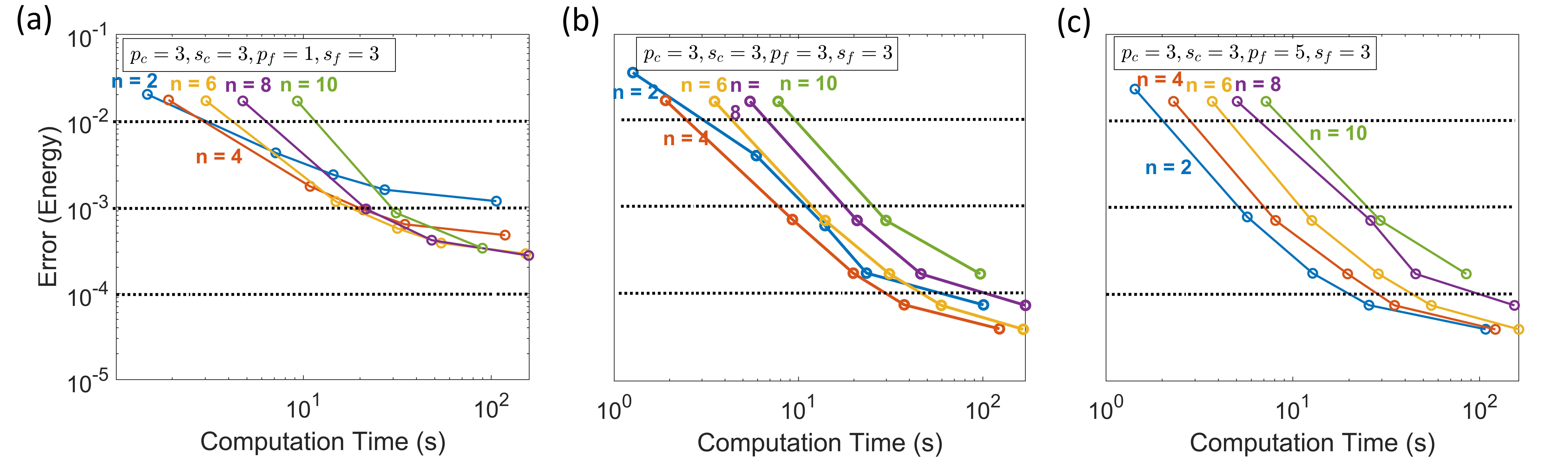}
% \subfigure{\includegraphics[height=1.8in]{ErrCost_3313.png}}
% \subfigure{\includegraphics[height=1.8in]{ErrCost_3333.png}}
% \subfigure{\includegraphics[height=1.8in]{ErrCost_3353.png}}
\caption{Error vs. computational time for different refinements and parameters: (a) $p_c=3, s_c =3, p_f=1, s_f=3$; (b) $p_c=3, s_c=3, p_f=3, s_f=3$; (c) $p_c=3, s_c=3, p_f=5, s_f=3$. Each curve represents the results for specific parameters ($p_c, s_c, p_f, s_f, n$) with refining coarse-level and fine-level meshes.}
\label{fig:Convergence_Eg1_ErrVsTime}
\end{figure}

To study the optimal mesh size ratio $n$ and orders $p_c$ and $p_f$ that leverage accuracy and efficiency, we plot energy-norm errors against computational costs for three cases with various refinements and element size ratios in Fig. \ref{fig:Convergence_Eg1_ErrVsTime}. Each curve in this figure represents the results for specific parameters ($p_c, s_c, p_f, s_f, n$), in which coarse-level and fine-level meshes are refined simultaneously for the remaining $n=h_c/h_f$. In the case of $(p_c=3, s_c=3, p_f=1, s_f=3)$, the shortest computational time for the error $10^{-2}$ is achieved for $n=2$ or $n=4$, while the optimal element size ratio for an error less than $10^{-3}$ is $n=6$ or $n=8$, which aligns with the observations in Fig. \ref{fig:Convergence_Eg1_Coarse} (a). The optimal element size ratios for $p_f=3, s_f=3$ and $p_f=5, s_f=3$ are consistently $n=4$ and $n=2$, respectively, corresponding to the turning points in Fig. \ref{fig:Convergence_Eg1_Coarse}. These observations indicate that the turning point of the convergence curve, i.e., $n_{opt}$, is the best element size ratio that leverages accuracy and efficiency.

\begin{figure}[htbp]
\centering
\subfigure{\includegraphics[height=2.7in]{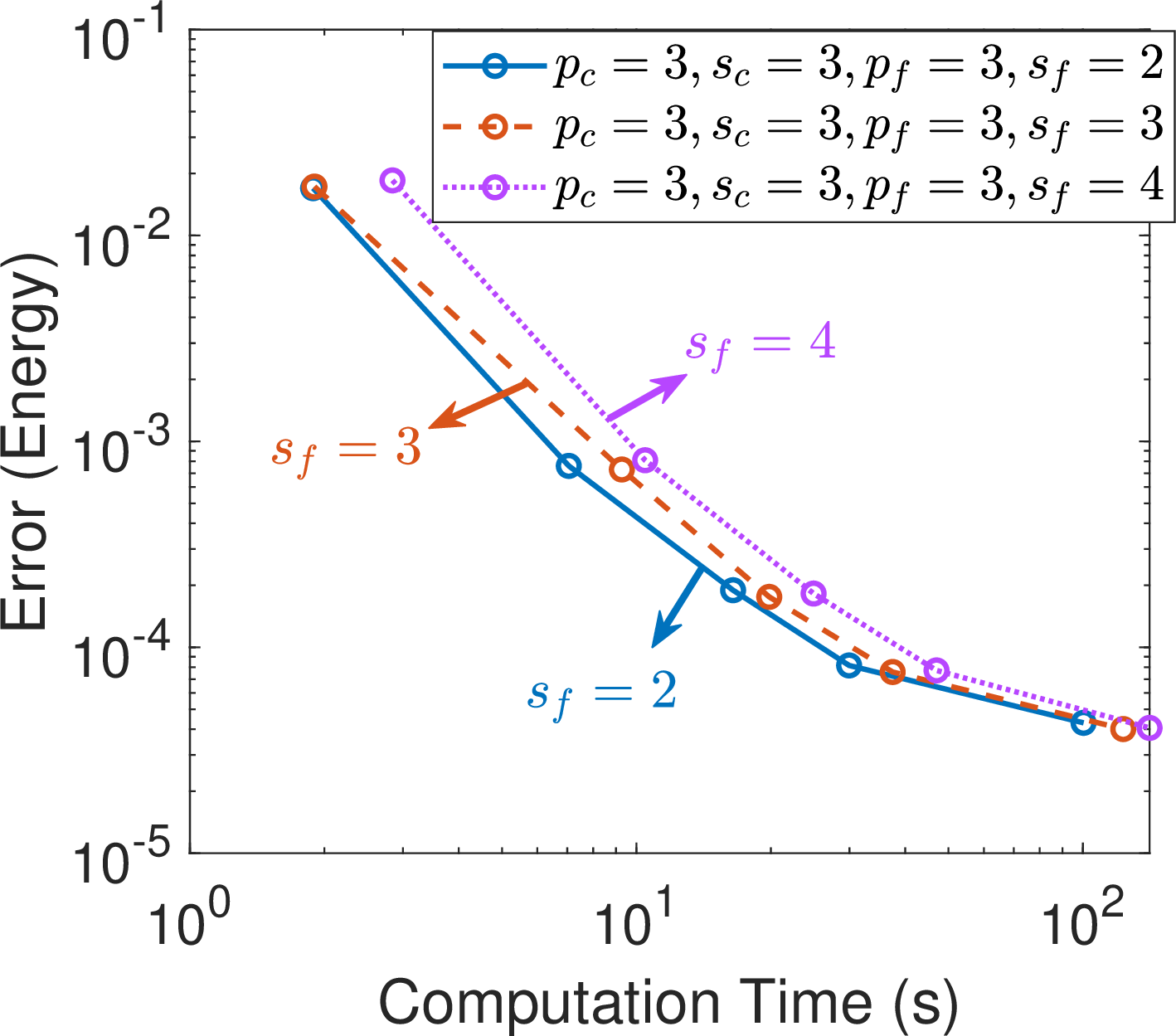}}
\subfigure{\includegraphics[height=2.7in]{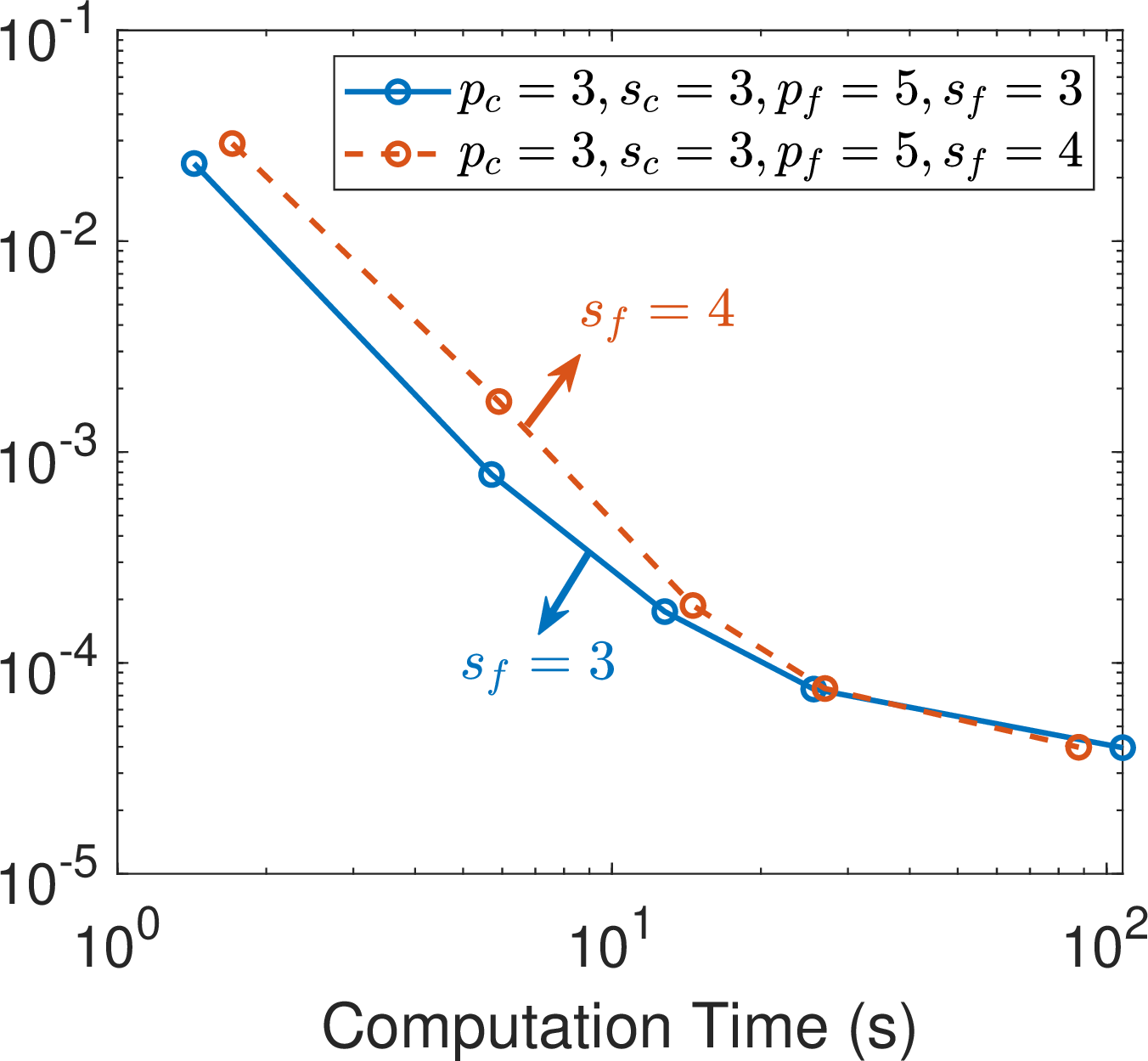}}
\caption{Error vs. computational time for different fine-level patch sizes: (a) $s_f = 2, 3$ and $4$ with other parameters $p_c=3, s_c =3, p_f=1, n=4$ consistent; (b) $s_f = 3$ and $4$ with other parameters $p_c=3, s_c =3, p_f=5, n=2$ consistent. }
\label{fig:Convergence_Eg1_S_ErrVsTime}
\end{figure}

We study the performance for different fine-level patch sizes $s_f$ while the other parameters ($p_c, s_c, p_f, n$) are consistent in Fig. \ref{fig:Convergence_Eg1_S_ErrVsTime}. In the case of $p_c=3, s_c=3, p_f=3$, we take the optimal element size ratio $n=4$. All curves ($s_f=2,3,4$) are close to each other, indicating that $s_f$ slightly influences the efficiency. In the case of $p_c=3, s_c=3, p_f=5$ with the optimal element size ratio $n=2$, the curve for $s_f=3$ is also close to that for $s_f=4$.

\begin{table}[!htb]
\caption{Computation time comparison for different refinements and parameters with comparable accuracy. The smallest computational time for each parameter setting is underlined. The smallest computational time for all cases is bolded.}
\centering
\begin{tabular}{|c| c c c | c c c | c c c |}
\hline
 & \multicolumn{3}{|c}{$p_c=3, s_c=3;$}& \multicolumn{3}{|c}{$p_c=3, s_c=3;$} &\multicolumn{3}{|c|}{$p_c=3, s_c=3;$} \\ 
 & \multicolumn{3}{|c}{$p_f=1, s_f=3$}& \multicolumn{3}{|c}{$p_f=3, s_f=3$} &\multicolumn{3}{|c|}{$p_f=5, s_f=3$} \\ \hline
error & $n=4$ & $n=6$ & $n=8$ & $n=2$ & $n=4$ & $n=6$ & $n=2$ & $n=4$ & $n=6$ \\ \hline
$10^{-2}$ & $\underline{2.20}$ & $3.22$ & $4.73$ & $2.38$ & $\underline{1.95}$ & $2.99$ & $\underline{\bm{1.62}}$ & $1.99$ & $3.10$ \\ \hline
$10^{-3}$ & $15.00$ & $\underline{11.19}$ & $14.48$ & $7.49$ & $\underline{5.19}$ & $7.91$ & $\underline{\bm{3.66}}$ & $5.17$ & $7.94$ \\ \hline
$10^{-4}$ & - & - & - & $78.33$ & $\underline{20.63}$ & $31.90$ & $\underline{\bm{14.54}}$ & $20.63$ & $32.26$ \\ \hline
\end{tabular}
\label{table:timecomparison_p=3}
\end{table}

\begin{figure}[htbp]
\centering
\subfigure{\includegraphics[width=0.6\textwidth]{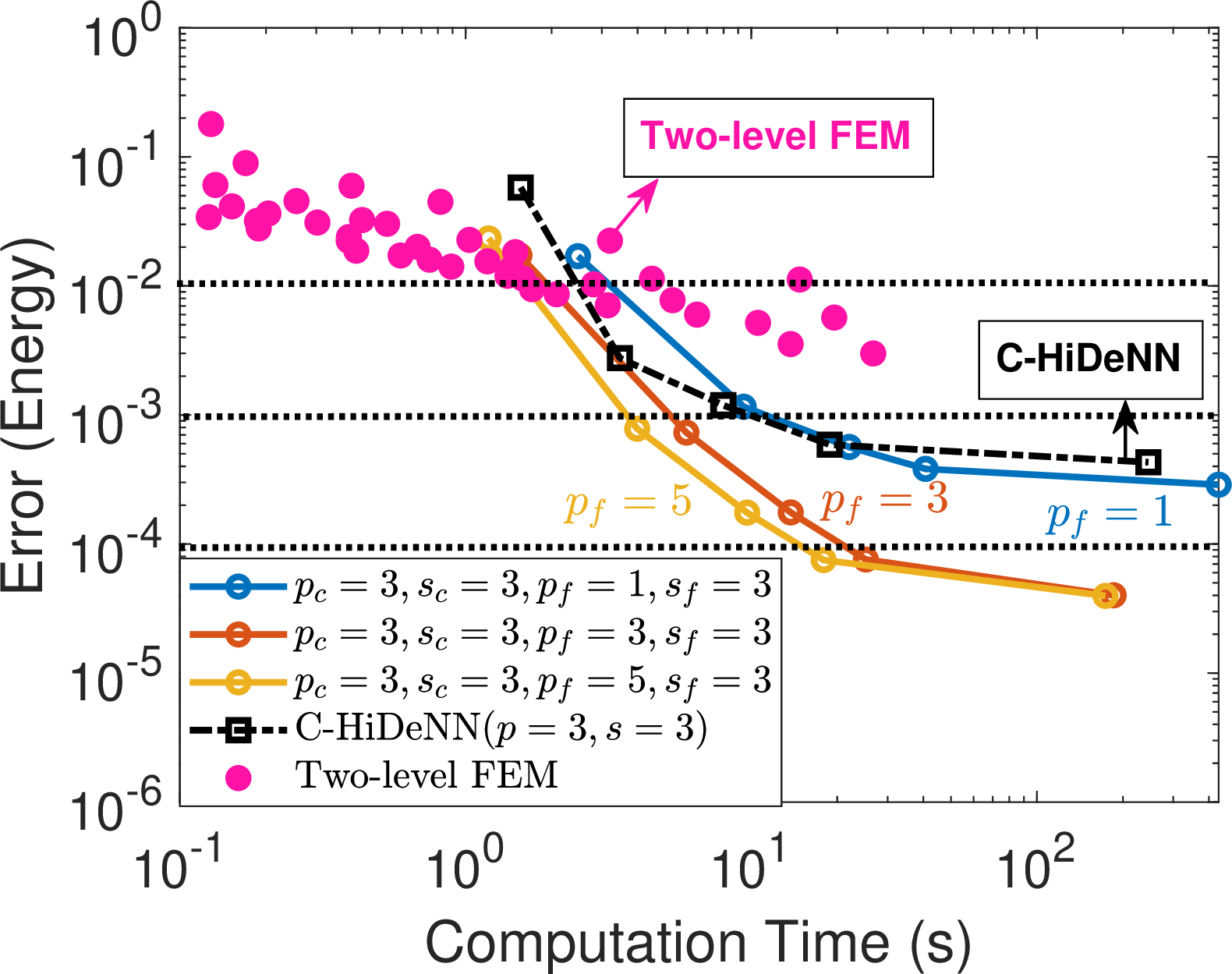}}
\caption{Comparisons among two-level C-HiDeNN with optimal element size ratios, two-level FEM and (single-level) C-HiDeNN. }
\label{fig:Comparison}
\end{figure}

\subsection{Single-level C-HiDeNN-TD and numerical verification of Theorem 3}

Before applying two-level VMS with C-HiDeNN-TD, we first use C-HiDeNN-TD (single-level) to verify Theorem 3, and study the performance with respect to the mode number.

The results for C-HiDeNN-TD and C-HiDeNN with parameters $p=3, s=3$ on the uniform mesh of $240$ by $240$ are listed in Table \ref{table:TD_comparison_p=3}. As expected, C-HiDeNN-TD has a much smaller number of unknowns than C-HiDeNN. The error for C-HiDeNN-TD decreases as the mode number increases. When the mode number is greater than five, C-HiDeNN-TD achieves the same level of accuracy and only requires thousands of DoFs, while C-HiDeNN requires 57121 DoFs. The C-HiDeNN error is the lowest. The energy-norm deviations between C-HiDeNN-TD and C-HiDeNN ($\|u^{\text{C-HiDeNN-TD}}-u^{\text{C-HiDeNN-TD}}\|_E$) are listed in the second column at the end. It is observed that the deviation decreases rapidly with the increase of mode number, and approaches zero. In particular, these deviations are consistently equal to $\sqrt{e^2_{\text{C-HiDeNN-TD}}-e^2_{\text{C-HiDeNN}}}$ except for the last row, which verifies our Theorem 3. The inconsistency in the last row is due to machine precision. 

\begin{table}[!htb]
\caption{Accuracy comparison for C-HiDeNN and C-HiDeNN-TD with parameters $p=3, s=3$ on the 240 by 240 mesh.}
\centering
\begin{tabular}{|c| c | c | c | c | c | c|}
\hline
Mode \# & \multicolumn{2}{c|}{C-HiDeNN-TD}  & \multicolumn{2}{c|}{C-HiDeNN} & Deviations & $\sqrt{e^2_{\text{C-HiDeNN-TD}}-e^2_{\text{C-HiDeNN}}}$ \\
\cline{2-5}
& DoFs & $e_{\text{C-HiDeNN-TD}}$ & DoFs & $e_{\text{C-HiDeNN}}$ & & \\
\hline
1 & 478 & $5.70 \times 10^{-1}$  & 57121 & $1.93 \times 10^{-4}$ & $5.70 \times 10^{-1}$ & $5.70 \times 10^{-1}$ \\
\hline
2 & 956 & $1.46 \times 10^{-1}$ & & & $1.46 \times 10^{-1}$ & $1.46 \times 10^{-1}$ \\
\hline
3 & 1434 & $2.12 \times 10^{-2}$ & & & $2.12 \times 10^{-2}$ & $2.12 \times 10^{-2}$ \\
\hline
4 & 1912 & $1.80 \times 10^{-3}$ & & & $1.79 \times 10^{-3}$ & $1.79 \times 10^{-3}$ \\
\hline
5 & 2390 & $2.13 \times 10^{-4}$ & & & $9.07 \times 10^{-5}$ & $9.07 \times 10^{-5}$ \\
\hline
6 & 2868 & $1.93 \times 10^{-4}$ & & & $2.08 \times 10^{-6}$ & $2.08 \times 10^{-6}$ \\
\hline
7 & 3346 & $1.93 \times 10^{-4}$ & & & $4.33 \times 10^{-8}$ & $3.80 \times 10^{-8}$ \\
\hline

\end{tabular}
\label{table:TD_comparison_p=3}
\end{table}

\begin{figure}[htbp]
\centering
\includegraphics[width=0.6\textwidth]{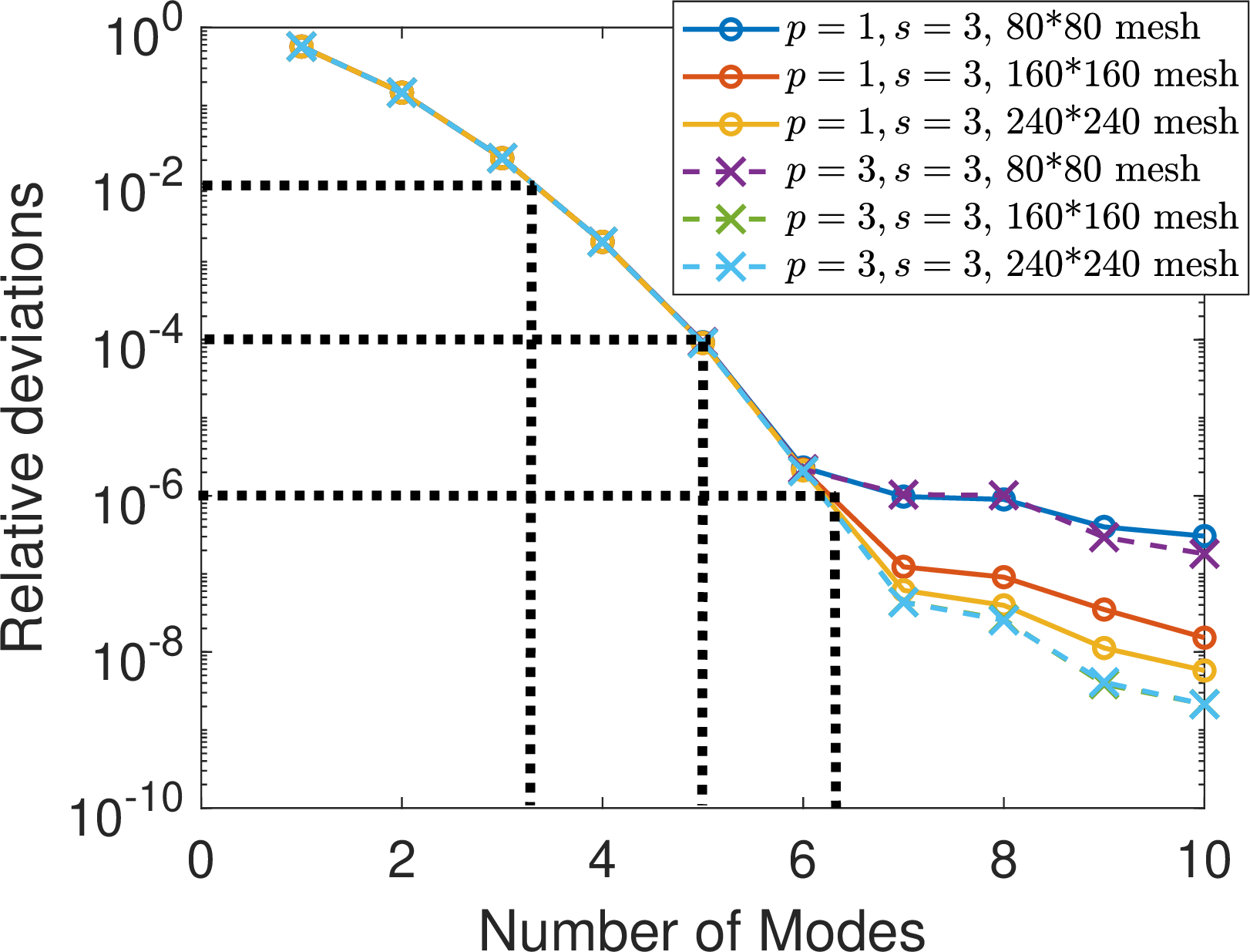}
\caption{Deviations between C-HiDeNN-TD and C-HiDeNN vs. the mode number. Before reaching the precise mode number of the exact solution (7 modes), the trend in deviations with respect to mode number is consistent across various refinements and parameters. This motivates us that we can estimate the mode number required based on the results over coarse meshes at a low computational cost. Here, the mode number greater than five is sufficient.  }
\label{fig:Dev_Q}
\end{figure}

We plot the deviations between C-HiDeNN-TD and C-HiDeNN with respect to the mode number in Fig. \ref{fig:Dev_Q}. The figure shows the trend in deviations across various levels of refinement and reproducing orders $p$. The following observations can be drawn from the figure:

$\bullet$ Initially, all curves nearly coincide, regardless of the refinements and parameter $p$, indicating that the deviations are dominated by the truncated error. 

$\bullet$ We can estimate the mode number based on the results obtained from coarse meshes, because the trend in deviations with respect to the mode number is consistent across different refinements and parameters.

$\bullet$ As the mode number approaches the precise mode number for the exact solution ($Q=7$), the deviations occur, because additional modes are solved to capture the discretization error. If much more modes are used, redundant modes will lead to the singularity in the stiffness matrix.

% The results for (single-level) C-HiDeNN-TD are shown in Fig. \ref{fig:SingleLevel_C-HiDeNN-TD}. As shown in the left subplot, the deviations between C-HiDeNN and C-HiDeNN-TD, which is the second term in the error decomposition Eq. (43), decrease as the mode number $Q$ increasing. This term achieves less than $10^{-7}$ for $Q=6$. Therefore, we take $Q=6$ and study the error vs. the element size $h$. The errors for $p=3$ and $p=5$ both achieve an saturation at about $10^-7$, that is limited by the mode number. The above observations verify our error decomposition for C-HiDeNN-TD presented in Theorem 3.

% \begin{figure}[htbp]
% \centering
% \includegraphics[width=6in]{fig/SingleLevel.png}
% \caption{Convergence curves for (single-level) C-HiDeNN-TD: (a) deviations between C-HiDeNN and C-HiDeNN-TD vs. mode number $Q$; (b) error vs. $h$ for $Q=6$.}
% \label{fig:SingleLevel_C-HiDeNN-TD}
% \end{figure}

\subsection{Two-level VMS with C-HiDeNN-TD}

\begin{figure}[htbp]
\centering
\subfigure{\includegraphics[height=1.8in]{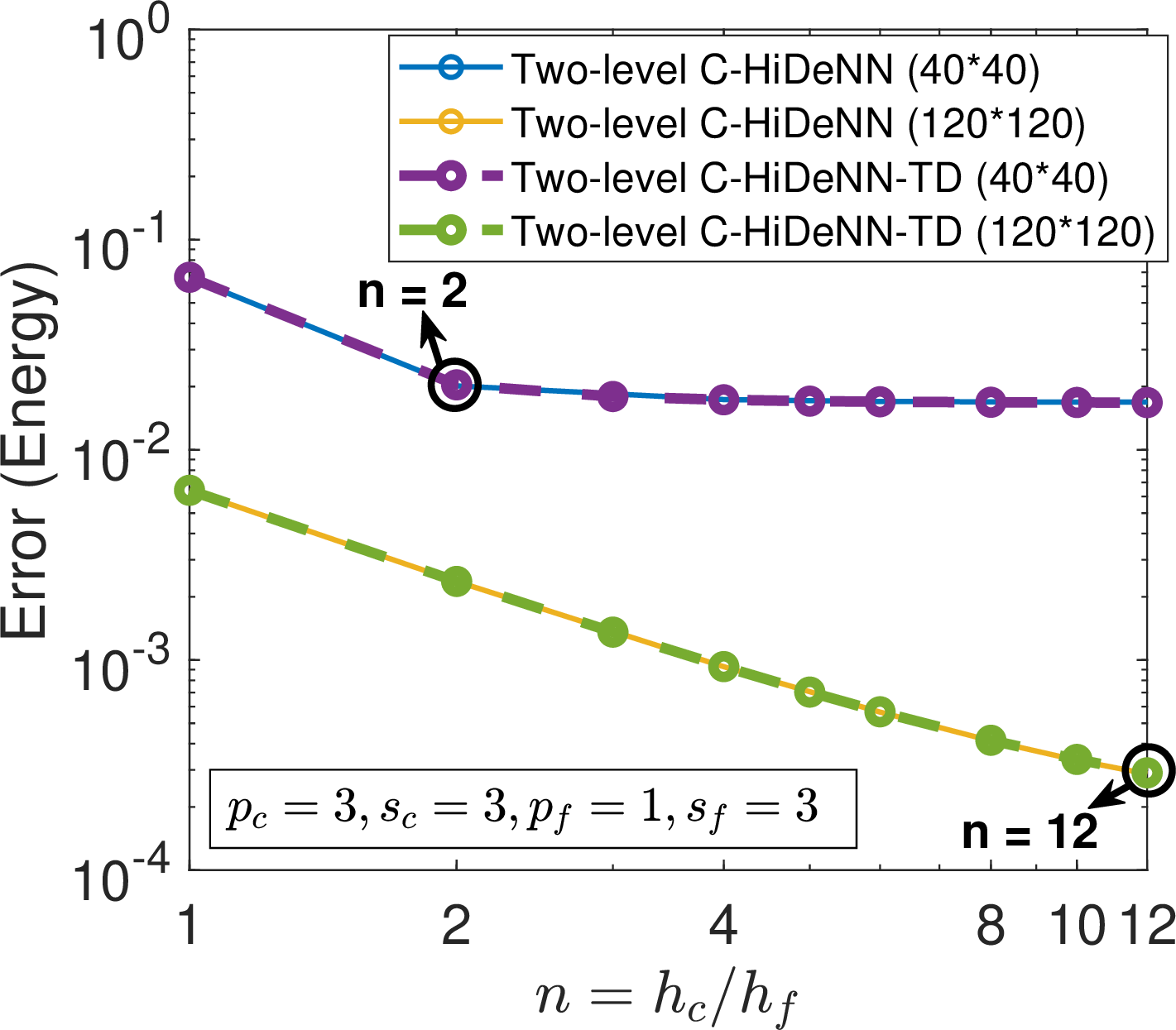}}
\subfigure{\includegraphics[height=1.8in]{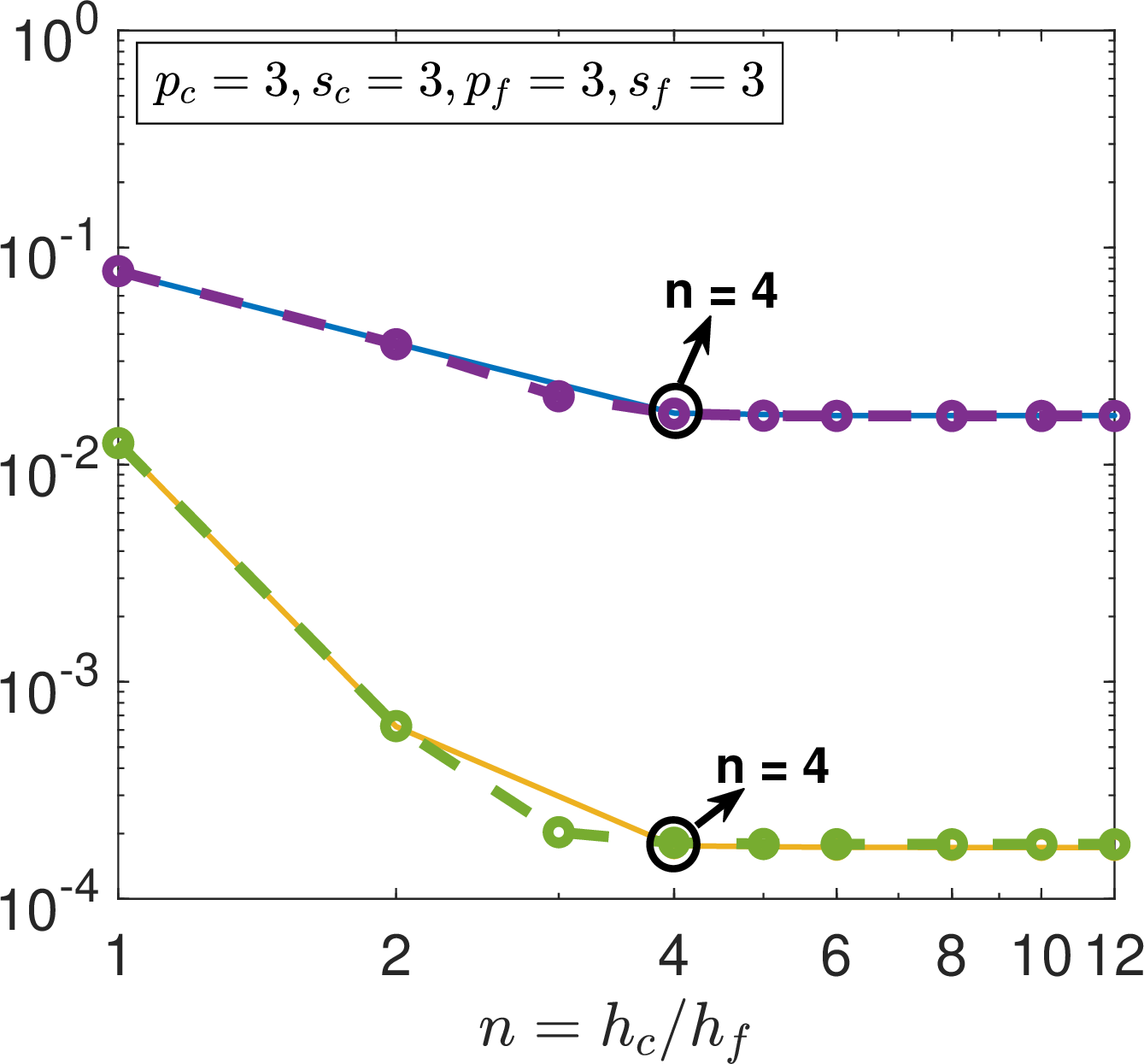}}
\subfigure{\includegraphics[height=1.8in]{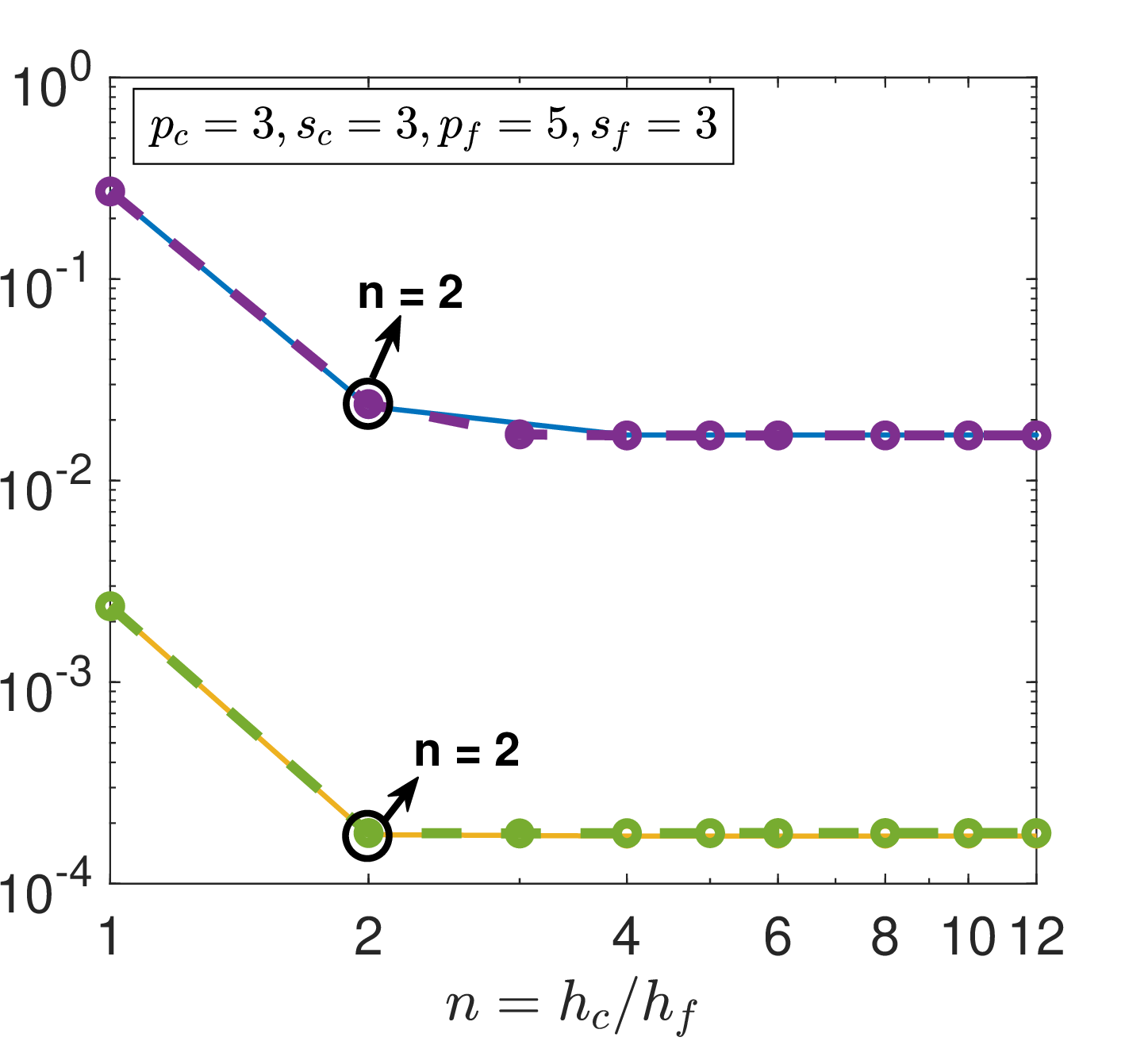}}
\caption{Convergence curves for refining fine-level mesh with fixed coarse-level mesh.}
\label{fig:Convergence_Eg1_Coarse_TD}
\end{figure}

To study the effect of element size ratio in the two-level VMS with C-HiDeNN-TD, we fix the coarse-level mesh and refine the fine-level mesh. As shown in Fig. \ref{fig:Convergence_Eg1_Coarse_TD}, the convergence curves for C-HiDeNN-TD align well with those for C-HiDeNN for various parameters. Hence, we observe similar trends: the errors initially decrease rapidly with refining the fine-level mesh, but maintain constant for very fine fine-level meshes, because the errors are bounded by coarse-level errors. The trend in turning points is also similar to that for two-level C-HiDeNN.

\section{Two-level VMS with different hyperparameters for 1D transient heat problem}

% To study the convergence rates, we first fix the coarse-level space time meshes (i.e., fix $h_c$ and $\Delta t_c$) and refine the fine-level space time meshes (i.e., refine $h_f$ and $\Delta t_f$ simultaneously). We take $6$ modes for both coarse- and fine-levels in C-HiDeNN-TD, as suggested in the analysis for single-level case. As shown in Fig. \ref{fig:Convergence_Eg3_Fine_C-HiDeNN-TD}, the results for C-HiDeNN-TD align well with those for C-HiDeNN, indicating that 6 modes is sufficient. As the fine-level meshes refining, the errors initially decrease at the $p_f$-th order, and then saturate at constants. The final error ($n=32$) for $p_c=5$ is lower than that for $p_c=3$. These observations verify our theories that the error for two-level VMS is controlled by two levels of parameters. 

% \begin{figure}[htbp]
% \centering
% \includegraphics[width=4in]{fig/Convergence_RefineFine_Heat_C-HiDeNN-TD.png}

% \caption{Convergence curves for refining fine-level mesh.}
% \label{fig:Convergence_Eg3_Fine_C-HiDeNN-TD}
% \end{figure}

\begin{figure}[htbp]
\centering
\subfigure{\includegraphics[height=1.8in]{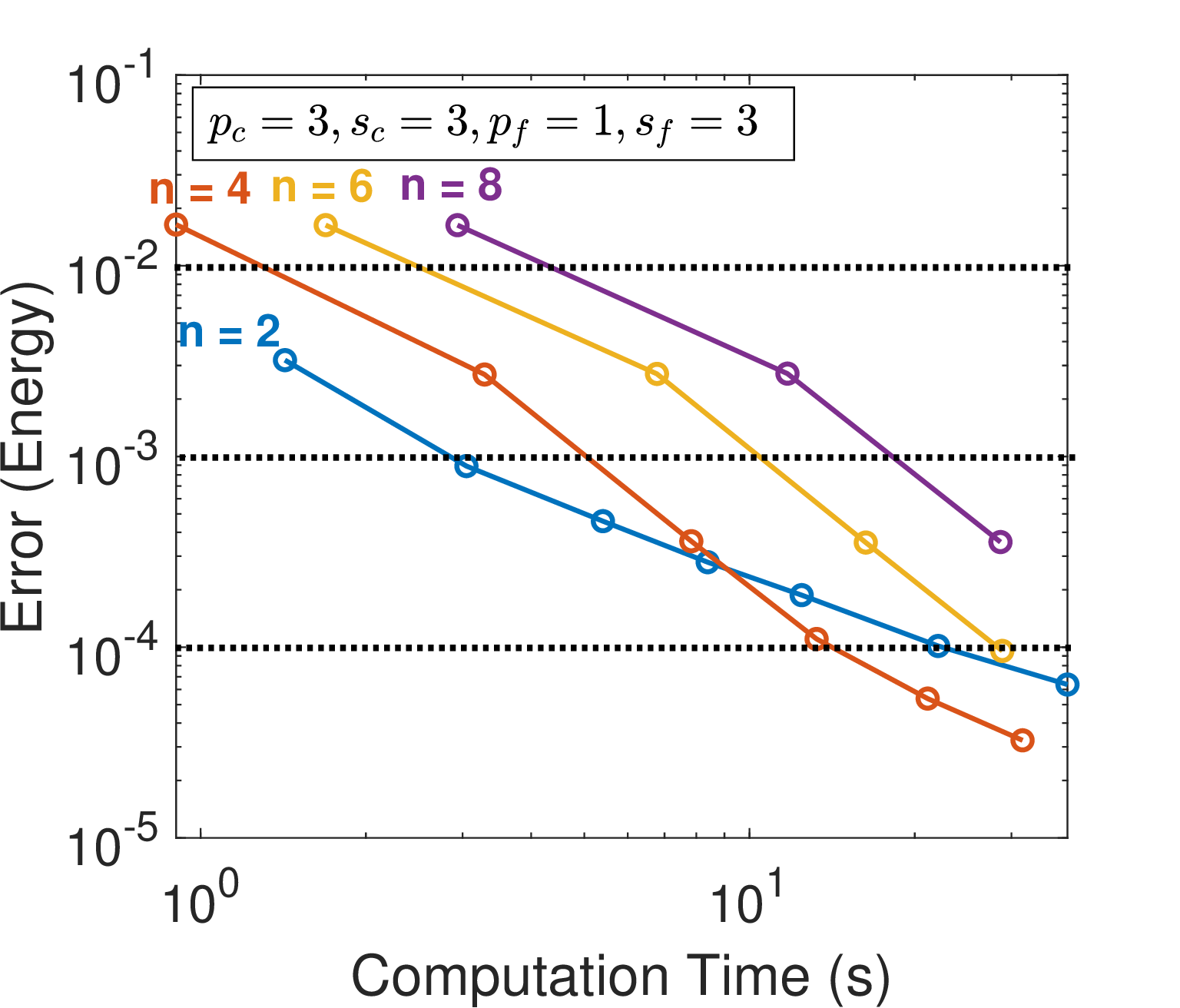}}
\subfigure{\includegraphics[height=1.8in]{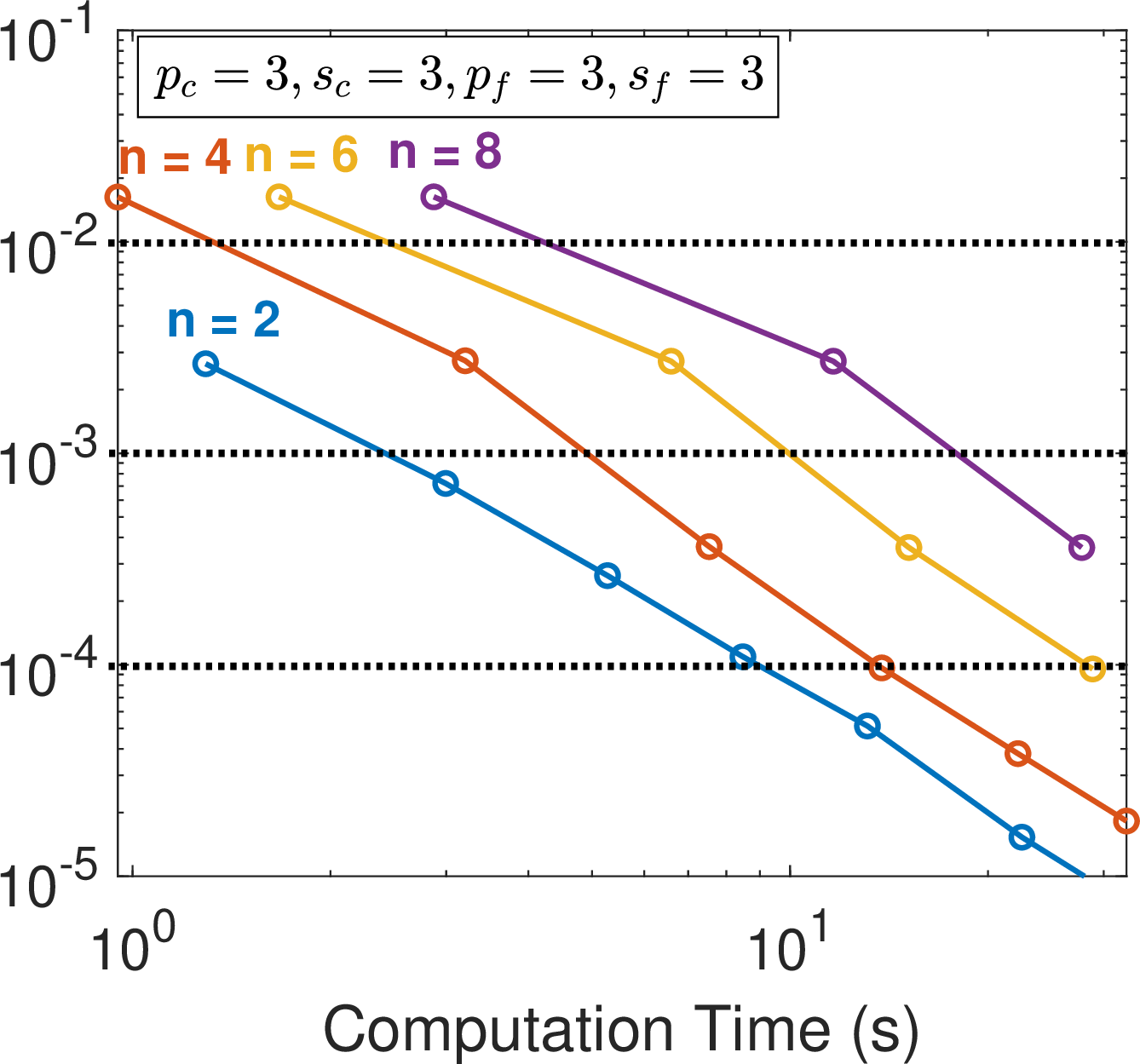}}
\subfigure{\includegraphics[height=1.8in]{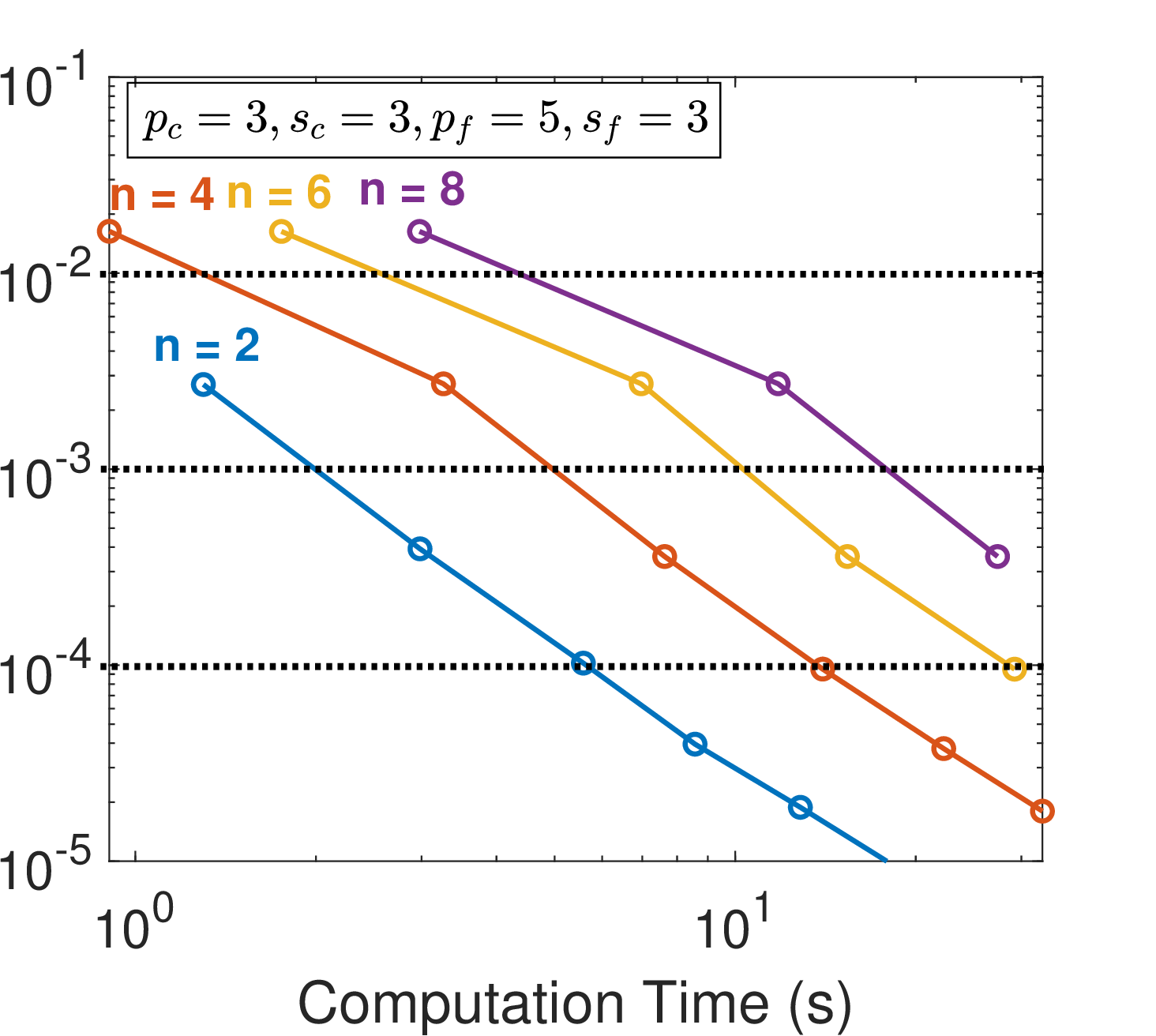}}
\caption{Error vs. computational time for different refinements and parameters: (a) $p_c=3, s_c =3, p_f=1, s_f=3$; (b) $p_c=3, s_c=3, p_f=3, s_f=3$; (c) $p_c=3, s_c=3, p_f=5, s_f=3$. Each curve represents the results for specific parameters ($p_c, s_c, p_f, s_f, n$) with refining coarse-level and fine-level meshes.}
\label{fig:Convergence_ST_ErrVsTime}
\end{figure}

We plot the errors against computational times for three cases (two-level C-HiDeNN-TD with $p_c=3, s_c =3, p_f=1, s_f=3$; $p_c=3, s_c=3, p_f=3, s_f=3$; and $p_c=3, s_c=3, p_f=5, s_f=3$) with various refinements of both two levels in Fig. \ref{fig:Convergence_ST_ErrVsTime}. Each curve represents the results for two-level C-HiDeNN-TD with refining coarse-level and fine-level meshes simultaneously to maintain a consistent element size ratio $n$. Similar to the observations for Poisson's problem, for $p_f=1<p_c=3$, the method for $n=2$ is the most efficient initially, while the method for $n=4$ takes less computational time than that for $n=2$ when error is less than $10^{-4}$. However, in the cases of $p_f=p_c=3$ and $p_f=5>p_c=3$, $n=2$ is always the optimal element size ratio. 

\section{Coordinate transformation for moving heat sources}
\label{appenx:CoordinateTransformation}

In this section, we present a general coordinate transformation algorithem for 3D heat problem with heat sources moving in the $x-y$ plane. Under such coordinate transformation, the problem of moving heat sources can be solved by C-HiDeNN-TD scheme efficiently.

\begin{figure}[htbp]
\centering
\includegraphics[width=5in]{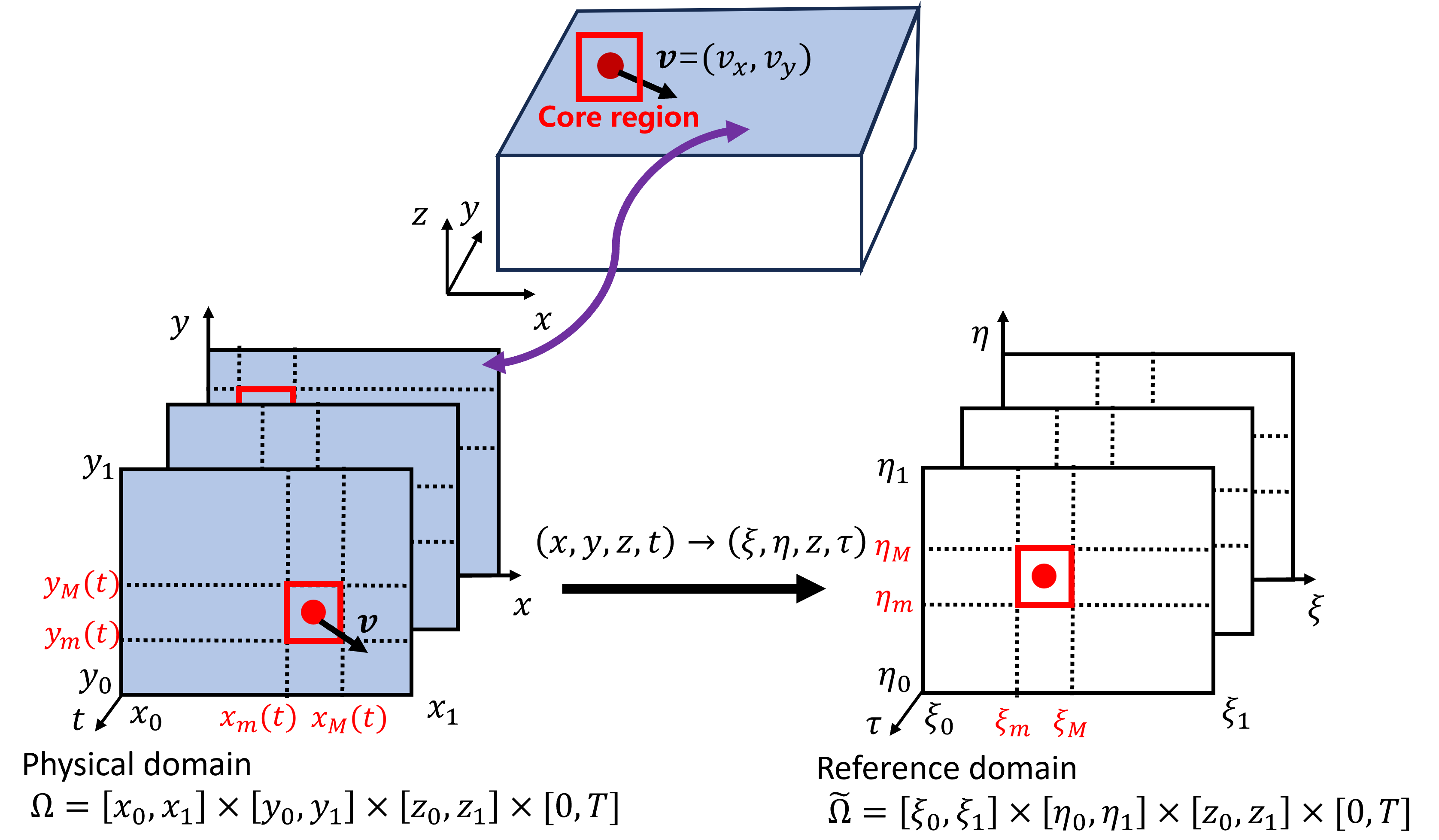}

\caption{Schematic diagram of the coordinate transformation for moving hear sources. The red point is the position of the laser source. }
\label{fig:Mapping_Moving}
\end{figure}

The computational domain $\Omega(\bm{x})$ is $[x_0, x_1]\times[y_0, y_1]\times[z_0, z_1]$. Consider a localized heat source moving at the velocity $\bm{v}=(v_x, v_y)$. We take a box that covers the heat source, associated with the fine-level domain that tracks the source in the framework of ML-VMS. The box is $[x_m(t), x_M(t)]\times[y_m(t), y_M(t)]$ in the $x-y$ plane at time $t$, which also moves at the velocity $\bm{v}=(v_x, v_y)$, as illustrated in Fig. \ref{fig:Mapping_Moving}. The entire $x-y$ plane can be divided into 9 subdomains with the edges of $x=x_m(t), x=x_M(t), y=y_m(t)$ and $y=y_M(t)$. We expect to design a coordinate transformation $(x, y)\to(\eta, \xi)$ to map these subdomains into fixed ones in the transformed space. The domain $[x_m(t), x_M(t)]\times[y_m(t), y_M(t)]$ corresponds to $[\xi_m, \xi_M]\times[\eta_m, \eta_M]$, which is independent of time. The entire $\xi-\eta$ plane is also divided into 9 subdomains with the edges of $\xi=\xi_m, \xi=\xi_M, \eta=\eta_m$ and $y=y_M$, corresponding to 9 subdomains in $x-y$ plane. The corresponding reference domain $\tilde{\Omega}(\bm{\xi})$ in the $\xi-\eta$ plane is $[\xi_0, \xi_1]\times[\eta_0, \eta_1]$.

For clarity, one subdomain in the $x-y$ plane is $[x_a(t), x_b(t)]\times[y_a(t), y_b(t)]$ with
\begin{eqnarray}
    && [x_a(t), x_b(t)]\in \left\{[x_0, x_m(t)], [x_m(t), x_M(t)], [x_M(t), x_1(t)] \right\} \\ \nonumber
    \text{and} &&[y_a(t), y_b(t)]\in \left\{[y_0, y_m(t)], [y_m(t), y_M(t)], [y_M(t), y_1(t)] \right\}.
\end{eqnarray}
 The corresponding subdomain in the $\xi-\eta$ plane is $[\xi_a, \xi_b]\times[\eta_a, \eta_b]$ with $[\xi_a, \xi_b]\in \{[\xi_0, \xi_m], [\xi_m, \xi_M], [\xi_M, \xi_1] \}$ and $[\eta_a, \eta_b]\in \{[\eta_0, \eta_m], [\eta_m, \eta_M], [\eta_M, \eta_1] \}$. A linear mapping is taken to realize such transformation:  
\begin{equation}
    x=x_b(t)\dfrac{\xi-\xi_a}{\xi_b-\xi_a}+x_a(t)\dfrac{\xi_b-\xi}{\xi_b-\xi_a},
    y=y_b(t)\dfrac{\eta-\eta_a}{\eta_b-\eta_a}+y_a(t)\dfrac{\eta_b-\eta}{\eta_b-\eta_a}.
\end{equation}
The transformations in other dimensions ($z$-dimension and $t$-dimension) are identical. In the transformed space, the heat source is fixed within the subdomain $[\xi_m, \xi_M]\times[\eta_m, \eta_M]$.

The Jacobian matrix for the transformation is:
\begin{equation}
    \dfrac{\partial (x,y,t)}{\partial (\xi,\eta,\tau)}=\left[\begin{array}{ccc}
       \dfrac{x_b(t)-x_a(t)}{\xi_b-\xi_a}  & 0 & 0 \\
       0  & \dfrac{y_b(t)-y_a(t)}{\eta_b-\eta_a} & 0 \\
       \dfrac{\mathrm{d} x_b(t)}{\mathrm{d} t} \dfrac{\xi-\xi_a}{\xi_b-\xi_a}+\dfrac{\mathrm{d} x_a(t)}{\mathrm{d} t} \dfrac{\xi_b-\xi}{\xi_b-\xi_a} & \dfrac{\mathrm{d} y_b(t)}{\mathrm{d} t} \dfrac{\eta-\eta_a}{\eta_b-\eta_a}+\dfrac{\mathrm{d} y_a(t)}{\mathrm{d} t} \dfrac{\eta_b-\eta}{\eta_b-\eta_a} &  1
    \end{array}\right].
\end{equation}
We take the transformation $\tau=t$ in the time dimension.
The $z$-dimension is omitted, because the transformation in this dimension is identical and independent of transformations in other dimensions. The inverse of this Jacobian matrix is
\begin{equation}
    \dfrac{\partial (\xi,\eta,\tau)}{\partial (x,y,t)}=\left[\begin{array}{ccc}
       A(t)  & 0 & 0 \\
       0  & B(t) & 0 \\
       C(\xi,t) & D(\eta, t) &  1
    \end{array}\right]
\end{equation}
with
\begin{eqnarray}
    &&A(t) = \dfrac{\xi_b-\xi_a}{x_b(t)-x_a(t)},
    B(t) = \dfrac{\eta_b-\eta_a}{y_b(t)-y_a(t)}, \\ \nonumber
    &&C(\xi,t) = -\left( \dfrac{\mathrm{d} x_b(t)}{\mathrm{d} t} \dfrac{\xi-\xi_a}{x_b(t)-x_a(t)}+\dfrac{\mathrm{d} x_a(t)}{\mathrm{d} t} \dfrac{\xi_b-\xi}{x_b(t)-x_a(t)} \right), \\ \nonumber
    && D(\eta, t) = -\left( \dfrac{\mathrm{d} y_b(t)}{\mathrm{d} t} \dfrac{\eta-\eta_a}{y_b(t)-y_a(t)}+\dfrac{\mathrm{d} y_a(t)}{\mathrm{d} t} \dfrac{\eta_b-\eta}{y_b(t)-y_a(t)} \right).
\end{eqnarray}
which yields the transformations between the derivative operators:
\begin{equation}
    \dfrac{\partial }{\partial x} = A(t)\dfrac{\partial }{\partial \xi},
    \dfrac{\partial }{\partial y} = B(t)\dfrac{\partial }{\partial \eta},
    \dfrac{\partial }{\partial t} = C(\xi,t)\dfrac{\partial }{\partial \xi}+D(\eta,t)\dfrac{\partial }{\partial \eta}+\dfrac{\partial }{\partial \tau},
\end{equation}

As a result, the weak form in the transformed coordinate system can be expressed as:
\begin{equation}
    \int_\Omega \int_0^T (w \cdot u_t + k \nabla w \cdot \nabla u ) \mathrm{d} \bm{x} \mathrm{d} t
    = \int_{\tilde{\Omega}} \int_0^T (w (Cu_{\xi}+Du_{\eta}+u_{\tau}) + k A^2 w_\xi u_{\xi} + k B^2 w_\eta u_\eta + k w_z u_z) \mathrm{det}(\bm{J}) \mathrm{d} \bm{\xi} \mathrm{d} \tau.
\end{equation}
The Jacobian determinant is $\mathrm{det}(\bm{J})=1/(A(t)B(t))$. In particular, $A(t), B(t), C(\xi,t), D(\eta, t)$ and $\mathrm{det}(\bm{J})$ can be written in simple TD form. As a result, this 4D integral can be written as the product of a 1D integral directly, when taking TD weighting functions $w$ and solution $u$, achieving great efficiency.

\bibliography{reference}

% \end{Large}

\end{document}